\documentclass[11pt]{article}

\usepackage{amsmath, amsfonts, amssymb,latexsym,wasysym,graphicx}
\input epsf.sty

\newtheorem{Theorem}{Theorem}[section]
\newtheorem{Lemma}{Lemma}[section]
\newtheorem{Proposition}[Lemma]{Proposition}
\newtheorem{Corollary}[Lemma]{Corollary}
\newtheorem{Definition}[Lemma]{Definition}

\newtheorem{Hypothesis}[Lemma]{Hypothesis}

\newcommand{\BEQ}{\begin{equation}}     
\newcommand{\BEA}{\begin{eqnarray}}
\newcommand{\BD}{\begin{displaymath}}
\newcommand{\EEQ}{\end{equation}}       
\newcommand{\EEA}{\end{eqnarray}}
\newcommand{\ED}{\end{displaymath}}

\newcommand{\del}{\delta}
\newcommand{\Del}{\Delta}
\newcommand{\eps}{\varepsilon}          

\newcommand{\supp}{{\mathrm{supp}}}
\newcommand{\SFLF}{{\mathrm{SFLF}}}

\newcommand{\Region}{{\mathrm{Region}}}

\newcommand{\R}{\mathbb{R}}
\newcommand{\C}{\mathbb{C}}
\newcommand{\Z}{\mathbb{Z}}
\newcommand{\N}{\mathbb{N}}
\newcommand{\D}{\mathbb{D}}
\renewcommand{\P}{\mathbb{P}}

\newcommand{\Id}{{\mathrm{Id}}}

\def\proba{{\mathbb{P}}}
\def\esper{{\mathbb{E}}}
\def\T{{\mathbb{T}}}
\def\F{{\mathbb{F}}}

\newcommand{\eop}{\hfill $\Box$}        

\newcommand{\II}{{\rm i}}               
\newcommand{\half}{{1\over 2}}          
\renewcommand{\vec}[1]{\boldsymbol{#1}} 

\catcode`\@=11
\def\numberbysection{\@addtoreset{equation}{section}
        \def\theequation{\thesection.\arabic{equation}}}
\numberbysection

\begin{document}

\vspace*{1.5cm}
\begin{center}
{\Large \bf From constructive field theory to fractional stochastic calculus. (I) The L\'evy area of fractional
Brownian motion with Hurst index $\alpha\in(\frac{1}{8},\frac{1}{4})$}

\end{center}

\vspace{2mm}
\begin{center}
{\bf  Jacques Magnen and J\'er\'emie Unterberger}
\end{center}

\vspace{2mm}
\begin{quote}

\renewcommand{\baselinestretch}{1.0}
\footnotesize
{Let $B=(B_1(t),\ldots,B_d(t))$ be a $d$-dimensional fractional Brownian motion
with Hurst index $\alpha<1/4$. Defining properly iterated integrals of $B$ is a difficult
task because of the low H\"older regularity index of its paths. Yet rough path theory shows it is
the key to the construction of a stochastic calculus with respect to $B$, or to solving differential
equations driven by $B$.

We show in this paper how to obtain second-order iterated integrals as the limit when the
ultra-violet cut-off goes to infinity of iterated integrals
of weakly interacting fields defined using the tools of constructive field theory, in particular, cluster expansion and renormalization. The construction extends to a large class of Gaussian fields with the same short-distance behaviour, called multi-scale Gaussian fields. Previous constructions \cite{Unt-Holder,Unt-fBm} were of algebraic
nature and did not provide such a limiting procedure.
}
\end{quote}

\vspace{4mm}


\tableofcontents



\section{Introduction}


A major achievement of the probabilistic school since the middle of the 20th century is the
study of diffusion equations, in connection with Brownian motion or more generally
Markov processes -- and also with partial differential equations, through the Feynman-Kac
formula -- with many applications in physics and chemistry \cite{VK}. One of the main tools
is stochastic calculus with respect to semi-martingales $M$. An adapted integral such as $\int_s^t X(u)
dM(u)$ may be understood as a limit in some sense to be defined. Classically one uses
piecewise linear interpolations, $\sum_{s\le t_1<\ldots<t_N\le t} X(t_i)(M(t_{i+1})-
M(t_i))$ or $\sum_{s\le t_1<\ldots<t_N\le t} \frac{X(t_i)+X(t_{i+1})}{2} (M(t_{i+1})-
M(t_i))$; these approximations define in the limit $N\to\infty$
the It\^o, resp. Stratonovich integral. The latter one is actually obtained e.g. 
if $M=W$ is Brownian motion and $X(t)=f(W_t)$ with $f$ smooth   as the limit
$\lim_{\eps\to 0}\int_s^t f(W_{\eps}(u))dW_{\eps}(u)$ for any smooth approximation
$(W_{\eps})_{\eps>0}$ of $W$ converging a.s. to $W$  (see \cite{WZ}, or
\cite{KS} p. 169). The Stratonovich integral
$\int_s^t X(u)d^{Strato} M(u)$ has an advantage over the It\^o integral in that it agrees
with the fundamental theorem of calculus, namely, $F(M(t))=F(M(s))+\int_s^t F'(M(u)) d^{Strato}M(u)$.

\medskip

The semi-martingale approach fails altogether when considering stochastic processes with
lower regularity. Brownian motion, and more generally semi-martingales (up to time reparametrization), are $(1/2)^-$-H\"older,
i.e. $\alpha$-H\"older for any $\alpha<1/2$ \footnote{Recall that a continuous path $X:[0,T]\to\R$ is $\alpha$-H\"older, $\alpha\in(0,1)$, if $\sup_{s,t\in[0,T]} \frac{|X_t-X_s|}{|t-s|^{\alpha}}<\infty$.}.
 Processes with $\alpha$-H\"older paths, where
$\alpha\ll 1/2$, are maybe less common in nature but still deserve interest.
 Among these, the family of multifractional Gaussian processes is perhaps the most
widely studied \cite{PLV}, but one may also cite diffusions on fractals \cite{HL}, sub- or superdiffusions
in porous media \cite{Goldenfeld,Lesne} and the
fascinating multi-fractal random measures/walks in connection with turbulence and
two-dimensional quantum gravity \cite{BM,DS}. Although some of these processes appear in applications
rather as a random landscape in $D\ge 1$ space dimensions
(which raises related but different problems), we view them
here -- restricting ourselves to $D=1$ --
 as some colored noise integrated in time, and wish to define stochastic integration
with respect to them or to solve differential equations driven by them.

\medskip

We concentrate in this article on {\em multiscale Gaussian processes}
 (the terminology is ours) with {\em scaling dimension} or
more or less equivalently {\em H\"older regularity} $\alpha\in(0,1/2)$, the best-known example of
which being {\em fractional Brownian motion} (fBm for short) with Hurst index $\alpha$, $B^{\alpha}(t)$ or simply $B(t)$ \footnote{It is (up to a constant) the unique self-similar Gaussian process with stationary increments. The last property implies that its derivative is a {\em stationary field}.}.
We consider more precisely  a {\em two-dimensional fBm}, $B(t)=(B_1(t),B_2(t))$, with independent, identically distributed components \footnote{The
one-dimensional case is very different and much simpler,  and has been treated in \cite{GNRV}.}.
 The simplest non-trivial stochastic integral is then
${\cal A}(s,t):=\int_s^t dB_1(t_1)\int_s^{t_1} dB_2(t_2)=\int_s^t (B_2(u)-B_2(s))dB_1(u) $,
a twice iterated integral, where $B=(B_1(t),B_2(t))$ is a two-component fBm. Since
$\int_s^t dB_1(t_1)\int_s^{t_1} dB_2(t_2)+\int_s^t dB_2(t_2)
\int_s^{t_2} dB_1(t_1)=(B_1(t)-B_1(s))(B_2(t)-B_2(s))$, one is mainly interested in
the antisymmetrized quantity (measuring an area, as follows from the Green-Riemann
formula),
\BEA  {\cal LA}(s,t)&:=&\int_s^t dB_1(t_1)\int_s^{t_1} dB_2(t_2)-\int_s^t dB_2(t_2)
\int_s^{t_2} dB_1(t_1) \nonumber\\
&=&\int_s^t (B_2(u)-B_2(s))dB_1(u)- (B_1(u)-B_1(s))dB_2(u), \nonumber\\ \EEA
called {\em L\'evy area}.  The
corresponding Stratonovich integral, obtained as a limit either by linear interpolation or
by more refined Gaussian approximations \cite{CQ02,Nua,Unt08,Unt08b}, has been shown to {\em diverge} as soon as
$\alpha\le 1/4$.

This seemingly no-go theorem, although clear and derived by straightforward computations
that we reproduce in short in section 1, appears to be a puzzle when put in front of the results of {\em rough path theory} \cite{Lyo98,LyoQia02,Gu,Lej,Lej-bis,FV}. We shall not enter into the details of this fascinating theory, with its deep connections with nilpotent groups, Hopf algebras and
sub-Riemannian geometry, leaving this to a forthcoming paper \cite{MagUnt-suite} where these structures
will be instrumental to define higher-order integrals. Let us only state that rough path
theory allows one to define integration with respect to an $\alpha$-H\"older path $\Gamma$
in terms of {\em substitutes of iterated integrals} of $\Gamma$ up to order $\lfloor 1/\alpha
\rfloor$, $\lfloor \ .\ \rfloor$=entire part, called {\em rough path over $\Gamma$}, satisfying some H\"older regularity property and
also two algebraic properties,  called respectively  {\em Chen} and {\em shuffle} property. The former
one is of geometric nature and may be seen as an additivity
property for the areas or volumes generated by double or multiple iterated
integrals. In particular, it implies the following compatibility property between
the twice iterated integral  $\cal A$ and the original process $B$,
\BEQ {\cal A}(s,t)={\cal A}(s,u)+{\cal A}(u,t)+(B_2(u)-B_2(s))(B_1(t)-B_1(u)). \label{eq:Chen2}\EEQ
The shuffle property, on the other hand, ensures once again that the fundamental theorem
of calculus holds. Whatever the exact formulation of these two properties, they hold true
trivially for the usual iterated integrals of $\Gamma$ if $\Gamma$ is smooth, and (being
of algebraic nature), also hold true trivially when the rough path over $\Gamma$
is constructed out of true iterated integrals of smooth paths by some limiting procedure
as explained above.

The main point is now the following. General theorems show (i) that rough paths
always exist ({\em existence theorem}); (ii) that any rough path over $\Gamma$ may be
constructed by some rather abstract limiting procedure \footnote{using pieces of horizontal
sub-Riemannian geodesics.} ({\em approximation theorem}). Unfortunately, these theorems
are not constructive; also, they are of pathwise nature, and thus not necessarily appropriate in the probabilistic setting which is our starting point. Recent results using
the underlying Hopf algebra structures \cite{Unt-Holder,Unt-fBm,Unt-resume,FoiUnt,Unt-BPHZ} give an explicit,
general solution to the existence
problem, well-suited in particular for Gaussian processes, but no one knows how to obtain
these rough paths by an explicit limiting procedure.

\bigskip

Our project in this series of papers is to define an explicit rough path over fBm with arbitrary Hurst index, or
more generally {\em multiscale Gaussian fields} (to be defined below) by an explicit,
probabilistically meaningful
limiting procedure, thus solving at last the problem of constructing a full-fledged,
Stratonovich-like integration with respect to fBm.

Roughly speaking, it is obtained by making $B=(B(1),B(2))$ interact through a
weak but singular quartic, non-local interaction, which plays the r\^ole of a squared
kinetic momentum, associated to the {\em rotation of the path}, and makes its L\'evy area {\em finite}. Following the common use
of quantum field theory, this is implemented by multiplying the Gaussian measure by
the exponential weight $e^{-\half c'_{\alpha} \int\int {\cal L}_{int}(\phi_1,\phi_2)(t_1,t_2) |t_1-t_2|^{-4\alpha}
dt_1 dt_2}$ \footnote{The unessential constant $c'_{\alpha}$ is fixed e.g. by demanding that the Fourier transform of
 the kernel $c'_{\alpha}|t_1-t_2|^{-4\alpha}$ is the function $|\xi|^{4\alpha-1}$.} , with
\BEQ {\cal L}_{int}(\phi_1,\phi_2)(t_1,t_2) =\lambda^2 \left\{
(\partial {\cal A}^+)(t_1)  (\partial {\cal A}^+)(t_2)+ (\partial {\cal A}^-)(t_1)(\partial {\cal A}^-)(t_2) \right\},  \EEQ
where:  $\lambda$ (the coupling parameter) is a small, positive constant; $\phi_1,\phi_2$ are the (infra-red divergent) stationary fields associated to $B_1,B_2$, with covariance kernel as in eq. (\ref{eq:cov-phi}),
and similarly,  $ {\cal A}^{\pm}$ are stationary {\em left- and right-turning fields}, built
out of $\phi_1$, $\phi_2$, and representing the singular part of the L\'evy area (see section 1 for details). By assumption $\alpha<1/4$, so that the kernel $|t_1-t_2|^{-4\alpha}$ is locally integrable.
The statistical weight is maximal when $\partial
 {\cal A}^+=
\partial {\cal A}^-=0$, i.e. for sample paths which are ``essentially''  straight lines. Another way to
motivate this interaction is to understand that the divergence of the L\'evy area is due to the accumulation in a small region of space of small loops \cite{Lej}; the statistical
weight is unfavorable to such an accumulation. On the other hand, the quantities in the first-order Gaussian chaos (in other words, the $n$-point functions $\langle
B_{i_1}(x_1)\ldots B_{i_n}(x_n)\rangle_{\lambda}=\frac{1}{Z}\esper\left[B_{i_1}(x_1)\ldots B_{i_n}(x_n)
e^{-\half c'_{\alpha}  \int\int  {\cal L}_{int}(\phi_1,\phi_2)(t_1,t_2) |t_1-t_2|^{-4\alpha} dt_1 dt_2 } \right]$,
$i_1,\ldots,i_n=1,2$, where $Z:=
\esper\left[
e^{-\half c'_{\alpha}  \int\int  {\cal L}_{int}(\phi_1,\phi_2)(t_1,t_2) |t_1-t_2|^{-4\alpha} dt_1 dt_2 } \right]$
is a normalization constant called {\em partition function})
  -- associated to pure {\em translation} -- are insensitive to the interaction. Thus the compatibility property
(\ref{eq:Chen2}) is satisfied, and one has constructed a L\'evy area for fBm.

It would be  very interesting to find an equivalent pathwise description of this
interaction in terms of a {\em constrained}, long-memory self-interacting walk, see
section 1 for more comments. The interaction ${\cal L}_{int}$
 would thus play the r\^ole of a Ginzburg-Landau functional,
whose usual, local version  is conjectured  to give e.g. a field-theoretic description of the
Ising model. Proving such an equivalence should be much easier here (despite the
unusual non-local features) since it is a one-dimensional model.

Starting from the above field-theoretic description, the proof of finiteness and H\"older regularity of the L\'evy area for $\lambda>0$ small enough follows the well-established
scheme of constructive field theory. The main theorem may be stated as follows. As a rule, we denote in this article by $\esper[...]$ the Gaussian expectation and
by $\langle ...\rangle_{\lambda,\rho}$
 the expectation with respect to the $\lambda$-weighted interaction measure with scale $\rho$ ultraviolet
cut-off, so that in particular $\esper[...]=\langle ...\rangle_{0,\infty}$.

\begin{Theorem} \label{th:0.1}

Assume $\lambda\in(\frac{1}{8},\frac{1}{4})$. Consider for $\lambda>0$ small enough the family of probability measures
 (also called: {\em $(\phi,\partial\phi,\sigma)$-model})
\BEQ\proba_{\lambda}^{\to\rho}(\phi_1,\phi_2)=
e^{-\half c'_{\alpha} \int\int dt_1 dt_2 |t_1-t_2|^{-4\alpha} {\cal L}_{int}(\phi_1,\phi_2)(t_1,t_2)} d\mu^{\to\rho}(\phi_1)
d\mu^{\to\rho}(\phi_2),\EEQ
where $d\mu^{\to\rho}(\phi_i)=d\mu(\phi_i^{\to\rho})$ is a Gaussian measure obtained  by an ultra-violet cut-off
at Fourier momentum $|\xi|\approx M^{\rho}$ ($M>1$), see Definition \ref{def:Cjphi} and
 section 5. Then $(\P_{\lambda}^{\to\rho})_{\rho}$
 converges in law when $\rho\to\infty$ to some   measure
$\P_{\lambda}$, and the associated twice iterated integral
 $\int_s^t d\phi^{\to\rho}_1(t_1)\int_s^{t_1}d\phi_2^{\to\rho}(t_2)$
 converges   to a well-defined quantity ${\cal A}(s,t)$ satisfying
the Chen and shuffle properties, with  variance
\BEQ \langle {\cal A}(s,t)^2\rangle_{\lambda}:=
\lim_{\rho\to\infty} \langle {\cal A}(s,t)^2\rangle_{\lambda,\rho}=(\frac{4}{\lambda^2} K_1+K_2)|t-s|^{4\alpha}. \EEQ

 Furthermore, the connected moments of ${\cal A}(s,t)$ in the interacting theory -- except its variance -- coincide with those of the free, Gaussian theory, allowing for the computation of moments of ${\cal A}(s,t)$ in terms of explicit Gaussian integrals. In particular, one finds:
\BEQ \langle {\cal A}(s,t)^{2n} \rangle_{\lambda}\le C_n|t-s|^{4n\alpha}, \quad n\ge 0.\EEQ
As an immediate corollary following from the Kolmogorov-Centsov lemma \cite{RevYor}
\footnote{This lemma states that a process $X$ such that $\esper |X(t)-X(s)|^{p}\le C |t-s|^{1+p\beta}$ has almost surely
$\beta^-$-H\"older paths. Classical hypercontractivity arguments  implies the same conclusion for
 a {\em  Gaussian} process provided simply  $\esper|X(t)-X(s)|^2< C|t-s|^{2\beta}$. }, ${\cal A}(s,t)$ has a.s. $(2\alpha)^-$-H\"older
paths.
\end{Theorem}

The only difficulty is to prove the convergence of the moments of finite distributions of the L\'evy area when $\rho\to\infty$;
let us see how to deduce the assertions of the theorem from this.

The weak convergence of (some subsequence of)  the family of laws $(\P_{\lambda}^{\to\rho})_{\rho}$
is actually a consequence  -- using standard arguments including Prohorov's theorem, \cite{Bil}, chap. 6, 8, 12 -- of the following estimates which imply that $(\P_{\lambda}^{\to\rho})_{\rho}$ is tight,

\BEQ \langle |\phi(t)-\phi(s)|^2\rangle_{\lambda,\rho}\le K|t-s|^{2\alpha},\EEQ
$K$ being a constant uniform in $\rho$.
This last estimate is proved in a stronger Fourier form \footnote{Namely,
$|\phi(t)-\phi(s)|^2\le \int_{|\xi|>\frac{1}{|t-s|}} \frac{K}{|\xi|^{1+2\alpha}} d\xi+\int_{|\xi|<\frac{1}{|t-s|}}
K \frac{|t-s|^2 |\xi|^2}{|\xi|^{1+2\alpha}} d\xi\le K'|t-s|^{2\alpha}.$},
\BEQ \langle |{\cal \phi}(\xi)|^2\rangle_{\lambda,\rho}\le \frac{K}{|\xi|^{1+2\alpha}}.\EEQ

Finally, standard multi-scale estimates, sketched in section 1, show that the generating function of connected moments
of finite distributions of the L\'evy area, $(u_1,\ldots,u_n)\mapsto \log \langle \exp \II \left( u_1 {\cal A}(s_1,t_1)+
\ldots+u_n {\cal A}(s_n,t_n) \right)\rangle_{\lambda,\infty}$,
 has a non-zero radius of convergence, hence (see \cite{Fel2}, Chap. XV, \S 4) the moments of the finite
distributions of the L\'evy area determine its law, and {\em every} subsequence of $(\P^{\to\rho}_{\lambda})_{\rho}$ converges
to the same limit. It would also be possible (by a simple extension of the results of section 6) to prove a speed of
convergence for these moments, e.g.  $\langle \phi^j(s) \phi^j(t)\rangle_{\lambda,\rho}-\langle \phi^j(s)\phi^j(t)\rangle_{\lambda,\rho'}
=K_r \frac{M^{-2\alpha j}}{(1+M^j|t-s|)^r} O(M^{-(8\alpha-1)(\rho'-j)})$ (for some constant $K_r$ depending
on $r=1,2,\ldots$) if $j<\rho'<\rho$.

\bigskip

The result is not difficult to understand using the non-rigorous perturbation
theory. First, by a  trick explained in section 1, one replaces the non-local
interaction ${\cal L}(\phi_1,\phi_2)(t_1,t_2)|t_1-t_2|^{-4\alpha}$
 with a local interaction ${\cal L}(\phi_1,\phi_2,\sigma)(t)$ depending
on a two-component exchange particle field $\sigma=(\sigma_+(t),\sigma_-(t))$.
Then a Schwinger-Dyson identity (a functional integration by parts) relates the moments of $\cal A$ to those of $\sigma$. Simple power-counting arguments show that a connected $2n$-point function of $\sigma$ alone is superficially divergent if and only if $1-4n\alpha\ge 0$. Thus, restricting to $\alpha>1/8$, one only needs to renormalize the two-point function. Since the renormalized propagator of $\sigma$ is screened by an infinite mass term (as we shall see in section 1),
 the theory is free once one has integrated out the $\sigma$-field, hence one retrieves the underlying Gaussian theory $(\phi_1,\phi_2)$. The Schwinger-Dyson identity then shows that the two-point functions of $\cal A$ have been made finite.

\bigskip

Since constructive field theory \cite{Mas,Riv} is not a standard tool in probability, we shall spend
some time explaining in full details and generality its objects, in particular the
{\em cluster expansions} which lie at its heart.  Cluster
expansions are based on a simplified wavelet decomposition
$(\psi^j_{\Del})$ of the fields $\psi=\phi_1,\phi_2,\sigma_{\pm}$ where $j$ is a {\em vertical} (Fourier)
scale index, and $\Del$ a {\em horizontal} (i.e. in direct space)
 interval of size $M^{-j}$ around the center of the wavelet
component. Each $\psi^j_{\Del}$ is to be seen as a {\em degree of freedom} of the theory,
relatively independent from the others, so that the interaction may be expressed as a divergent
  doubly-infinite vertical and horizontal sum. The horizontal (H) and vertical (V)  cluster expansions allow to rewrite
the partition function $Z_V^{\to\rho}$ over a finite volume, with ultraviolet truncation at scale $\rho$,
 as a sum,
 \BEQ Z_V^{\to\rho}=\sum_n \frac{1}{n!} \sum_{\P_1,\ldots,\P_n
\ {\mathrm{ non-overlapping}}} F_{HV}(\P_1)\ldots F_{HV}(\P_n), \label{eq:0.9} \EEQ
where:

\noindent -- $\P_1,\ldots,\P_n$ are disjoint {\em polymers}, i.e. sets of intervals $\Del$ connected
by vertical and horizontal links;  during the course of the expansion, the Gaussian measure
has been modified so that the field components
belonging to different polymers have become independent;

\noindent -- $F_{HV}(\P)$, $\P=\P_1,\ldots,\P_n$ is the $\lambda$-weighted expectation value, $F_{HV}(\P)=
\langle f_{HV}(\P)\rangle_{\lambda}$,
of some function $f_{HV}$ depending only on the field components located in the support of $\P$.

The fundamental idea is that (i) the  {\em polymer evaluation function} $F(\P)$ is all the smaller as the polymer $\P$ is large, due to the polynomial decrease of correlation at large distances (for the horizontal direction), and
to power-counting arguments developed in section 4 for the vertical direction, leading to the image of horizontal islands
maintained together by vertical {\em springs};
 (ii) the horizontal and vertical links in $\P$ (once {\em one} interval belonging  to $\P$ has been fixed)  suppress the invariance by translation,  which normally leads to a divergence when $|V|\to\infty$.  A classical combinatorial trick, called {\em Mayer expansion}, allows one to rewrite
eq. (\ref{eq:0.9}) as a similar sum over {\em trees of polymers}, also called {\em Mayer-extended
polymers} and   denoted by the same letter $\P$, but {\em without
non-overlap conditions}, $Z_V^{\to\rho}=\sum_n \frac{1}{n!} \sum_{\P_1,\ldots,\P_n} F(\P_1)\ldots F(\P_n)$, where $F=F_{HVM}$ is the {\em Mayer-extended polymer evaluation function}, so
that $\ln Z_V^{\to\rho}=\sum_{\P} F(\P).$  In the process, {\em local parts of diverging graphs} have been resummed into
an exponential, leading to a {\em counterterm in the  interaction}; this is the essence of {\em renormalization}.
  On the whole, one finds that in the limit $|V|,
\rho\to \infty$, the free energy
$\ln Z_V^{\to\rho}$ is a sum over each scale of scale-dependent extensive quantities,
 i.e. $\ln Z_V^{\to\rho}=|V| \sum_{j=0}^{\to\rho} M^j f^{j\to\rho}_V$, where $f^{j\to\rho}_V$ converges when
 $|V|\to\infty$ to a finite quantity of order $O(\lambda)$. One retrieves the idea that each interval
 of scale $j$ contains one degree of freedom.
Finally, $n$-point functions are computed as derivatives of an external-field dependent version of the free energy.

\bigskip

Here is the outline of the article. Section 1 recalls the  classical results on the divergence
of the L\'evy area for $\alpha\le 1/4$ and recasts them into notations which are
more appropriate for cluster expansions. It also contains a heuristic section using
{\em perturbative}  field theory and giving some intuition on why adding the interaction
term ${\cal L}_{int}$ briefly described above should make the L\'evy area convergent.
Section 2 contains the definition of multiscale Gaussian fields, including fBm, and
many useful notations and general results concerning the scale (wavelet) decompositions. Section 3 is on cluster expansions, section 4 on renormalization.  Sections 2 to 4 are extremely
general, valid also for fields living on $\R^D$, in the hope that they may serve as a basis for future work, possibly also for $D>1$-models. Section 5, on the
contrary, concentrates on the definition of our specific $(\phi,\partial\phi,\sigma)$-model.
The proof of finiteness of $n$-point functions of the L\'evy area and freeness of the $\phi$-field   is given in section 6. It depends on {\em Gaussian bounds} which hold in great generality, and on {\em
domination arguments} which are very specific of our model.

\bigskip

{\bf Notations.} Cluster expansions imply the use of many indices and letters. Let us summarize here some
of our most important conventions, in the hope that this will help the reader not to get lost.

\begin{enumerate}

\item Quite generally, $\psi=(\psi_1(x),\ldots,\psi_d(x))$ is a $d$-dimensional field living on $\R^D$,
$D\ge 1$, with scaling dimensions $\beta_1,\ldots,\beta_d$. Roughly speaking, for probabilists, $\beta$ is the H\"older continuity index, at least when $\beta>0$; physicists usually call scaling dimension $-\beta$. Most of the model-independent results here are valid for arbitrary $D$. The Fourier transform is denoted by $\cal F$.

\item If $\psi=\psi(x)$ is a Gaussian field, then its scale decomposition (see section 1) reads
$\psi=\sum_{j\ge 0} \psi^j$. The {\em low-momentum}, resp. {\em high-momentum field} with respect to scale $j$
is denoted by $\psi^{\to j}=\sum_{h\le j} \psi^h$, resp. $\psi^{j\to}=\sum_{k\ge j} \psi^k$. All
through the text, we observe the following convention: if $h,j,k$ are scale indices,
then $h\le j\le k$; any primed scale index (for instance $j',\rho',\ldots$) is less than the
original scale index $j,\rho,\ldots$  {\em Secondary fields} (i.e. low-momentum fields $\psi$  minus their average) are denoted by
$\del \psi$. We also introduce {\em restricted high-momentum fields}, see section 1, denoted by $Res \ \psi$.

\item (products of fields) If $\psi=(\psi_1(x),\ldots,\psi_d(x))$ is a $d$-dimensional
field, and $I=(i_1,\ldots,i_n)\in\{1,\ldots,d\}^n$ $(n\ge 1$) is a multi-index, we denote
by $\psi_I$ the product of fields
$\psi_I(x):=\psi_{i_1}(x)\ldots\psi_{i_n}(x).$ The interaction Lagrangian ${\cal L}_{int}$ is written in general
as $\sum_{q=1}^p K_q \lambda^{\kappa_q} \psi_{I_q}$ for some constants $K_q$ and exponents $\kappa_q$.
 Cluster expansions produce propagators and products of fields
which are linear combinations of {\em monomials} (also called: {\em G-monomials}), generically denoted by $G$.
\item (constants) $K$ is a constant depending only (possibly) on the details of the model, such as the
degree of the interaction, the scaling dimension of the fields... It may vary from line to line. $M>1$
is the basis of the scale decomposition; it is an absolute constant, whose value is unimportant. $\gamma$ is
some constant $>1$, or sometimes simply $\ge 1$. $N_{ext,max}$ is such that every Feynman diagram with
$\ge N_{ext,max}$ external legs is superficially convergent; this perturbative notion also plays a central
r\^ole in constructive field theory. Our model $(\phi,\partial\phi,\sigma)$ has $N_{ext,max}=4$.
\item (variables) The maximum scale index is $\rho$. Other scale indices are denoted by $h,j,k$, or
any of those with primes. Indices $i$ are summation indices, with finite range, used in various contexts,
for instance for the field components. The parameters of the horizontal, resp. Mayer, resp. vertical
(also called momentum-decoupling)  cluster expansion,  are denoted by $s$, resp. $S$, resp. $t$.
 The scale of an interval $\Del$ is denoted by $j(\Del)$. If $\Del$ is an interval, then $n(\Del)$ is the coordination
number of $\Del$ inside the tree defined by the horizontal cluster expansion,
 while $N(\Del)$, resp. $N_i(\Del)$ is the total number of fields, resp. the number of
fields $\psi_i$ located in the interval $\Del$. $\tau$ stands for a number of derivations, usually
with respect to the $t$-parameters, in which case it is assumed to be $\le N_{ext,max}+O(n(\Del))$.
\end{enumerate}


\section{Statement of the problem and heuristics}


The quantity we want to define in the case of fractional Brownian motion is the following.

\begin{Definition}[L\'evy area]

The L\'evy area of a two-dimensional path $\Gamma:\R\to\R^2$ between $s$ and $t$ is
the area between the straight line connecting $(\Gamma_1(s),\Gamma_2(s))$ to $(\Gamma_1(t),\Gamma_2(t))$ and the graph $\{(\Gamma_1(u),\Gamma_2(u)); s\le u\le t\}$.
It is given by the following antisymmetric quantity,
\BEQ {\cal LA}_{\Gamma}(s,t):=\int_s^t d\Gamma_1(t_1) \int_s^{t_1} d\Gamma_2(t_2) -
\int_s^t d\Gamma_2(t_2) \int_s^{t_2} d\Gamma_1(t_1).\EEQ

\end{Definition}

\subsection{A Fourier analysis of the L\'evy area}

In order to understand the analytic properties of the L\'evy area of fBm, we shall resort to
a Fourier transform. One obtains, using the harmonizable representation of fBm introduced in subsection 2.2,
\BEQ {\cal A}(s,t):=\int_s^t dB_1(t_1) \int_s^{t_1} dB_2(t_2)=\frac{1}{2\pi c_{\alpha}}\int \frac{dW_1(\xi_1) dW_2(\xi_2)}{|\xi_1|^{\alpha-1/2} |\xi_2|^{\alpha-1/2}} \int_s^t dt_1 \int_s^{t_1} dt_2 \ \cdot\ e^{\II (t_1\xi_1+t_2\xi_2)}.
\EEQ

The L\'evy area ${\cal LA}(s,t):={\cal LA}_B(s,t)$ is obtained from this twice iterated integral by
antisymmetrization. Note that ${\cal LA}(s,t)$ is homogeneous of degree $2\alpha$ in $|t-s|$ since $B(ct)-B(cs)$, $c>0$ has same law as $c^{\alpha}(B(t)-B(s))$ by
self-similarity.

Expanding the right-hand side yields an expression which is not homogeneous in $\xi$.
Hence it is preferable to define instead the {\it skeleton integral}, depending only
on {\em one} variable,
\BEA {\cal A}(t):=\int^t dB_1(t_1)\int^{t_1} dB_2(t_2) &=& \frac{1}{2\pi c_{\alpha}}
\int \frac{dW_1(\xi_1) dW_2(\xi_2)}{|\xi_1|^{\alpha-1/2} |\xi_2|^{\alpha-1/2}} \int^t dt_1 \int^{t_1} dt_2 \ \cdot\ e^{\II (t_1\xi_1+t_2\xi_2)} \nonumber\\
&=& \frac{1}{2\pi c_{\alpha}} \int \frac{dW_1(\xi_1) dW_2(\xi_2)}{|\xi_1|^{\alpha-1/2} |\xi_2|^{\alpha-1/2}}
\cdot\ \frac{e^{\II t(\xi_1+\xi_2)}}{[\II(\xi_1+\xi_2)][\II\xi_2]}, \EEA
where by definition $\int^t e^{\II u\xi} du=\frac{e^{\II t\xi}}{\II\xi}.$ From
${\cal A}(t)$ and the one-dimensional skeleton integral
\BEQ \phi_i(t)=(2\pi c_{\alpha})^{-\half} \int^t dB_i(u)=\int \frac{dW_i(\xi)}{|\xi|^{\alpha-1/2}} \ \cdot\ \frac{e^{\II t\xi}}{\II \xi},\EEQ
which is the {\em infra-red divergent stationary process} associated to $B$, see section 2,
one easily retrieves ${\cal A}(s,t)$ since
\BEA {\cal A}(s,t) &=& \int_s^t dB_1(t_1) \left( \int^{t_1} dB_2(t_2) -\int^s dB_2(t_2)\right) \nonumber\\
&=& {\cal A}(t)-{\cal A}(s)+{\cal A}_{\partial}(s,t),  \label{eq:dec1} \EEA
where $(2\pi c_{\alpha})^{\half} {\cal A}_{\partial}(s,t):=(B_1(t)-B_1(s))\phi_2(s)$ (called {\em boundary term})
is a {\em product of first-order integrals}.

One may easily estimate these quantities in each sector $|\xi_1|\gtrless|\xi_2|$. In
practice, it turns out that estimates are easiest to get {\em after} a permutation of the
integrals (applying Fubini's theorem) such that (for twice or multiple iterated integrals
equally well) {\em innermost (or rightmost) integrals bear highest Fourier frequencies};
this is the essence of {\em Fourier normal ordering} \cite{Unt-Holder,FoiUnt,Unt-resume}. This gives a somewhat different
decomposition with respect to (\ref{eq:dec1}) since $\int_s^t dB_1(t_1)\int_s^{t_1}dB_2(t_2)$ is rewritten as $-\int_s^t dB_2(t_2)\int_t^{t_2} dB_1(t_1)$
in the "negative" sector $|\xi_1|>|\xi_2|$. After some elementary computations, one gets
the following.

\begin{Lemma}
Let
\BEA {\cal A}^+(t)&:=&2\pi c_{\alpha} \int^t dt_1 \int^{t_1} dt_2 {\cal F}^{-1}\left(
(\xi_1,\xi_2)\mapsto {\bf 1}_{|\xi_1|<|\xi_2|} ({\cal F}B'_1)(\xi_1) ({\cal F}B'_2)(\xi_2)
\right) (t_1,t_2) \nonumber\\
&=&  \int_{|\xi_1|<|\xi_2|} \frac{dW_1(\xi_1)dW_2(\xi_2)}{|\xi_1|^{\alpha-1/2} |\xi_2|^{\alpha-1/2}} \ \cdot\ \frac{e^{\II t(\xi_1+\xi_2)}}{[\II(\xi_1+\xi_2)][\II\xi_2]}
\EEA
and
\BEA {\cal A}^-(t)&:=& 2\pi c_{\alpha} \int^t dt_2 \int^{t_2} dt_1 {\cal F}^{-1}\left(
(\xi_1,\xi_2)\mapsto {\bf 1}_{|\xi_2|<|\xi_1|} ({\cal F}B'_1)(\xi_1) ({\cal F}B'_2)(\xi_2)
\right) (t_1,t_2) \nonumber\\
&=& \int_{|\xi_2|<|\xi_1|} \frac{dW_1(\xi_1)dW_2(\xi_2)}{|\xi_1|^{\alpha-1/2} |\xi_2|^{\alpha-1/2}} \ \cdot\ \frac{e^{\II t(\xi_1+\xi_2)}}{[\II(\xi_1+\xi_2)][\II\xi_1]}.
\EEA

Then
\BEQ {\cal A}(s,t)=\frac{1}{2\pi c_{\alpha}} \left\{ ({\cal A}^+(t)-{\cal A}^+(s))-({\cal A}^-(t)-{\cal A}^-(s))+
({\cal A}^+_{\partial}(s,t)-{\cal A}^-_{\partial}(s,t)) \right\}, \label{eq:1.8} \EEQ
where
\BEA && {\cal A}^+_{\partial}(s,t)-{\cal A}^-_{\partial}(s,t)= \left\{ -
\int_{|\xi_1|<|\xi_2|} \frac{(e^{\II t\xi_1}-e^{\II s\xi_1})e^{\II s\xi_2}}{[\II \xi_1][\II \xi_2]} \right. \nonumber\\
&&\left. \qquad \qquad +\int_{|\xi_2|<|\xi_1|} \frac{(e^{\II t\xi_2}-e^{\II s\xi_2})e^{\II t\xi_1}}{[\II \xi_1][\II \xi_2]} \right\}\ \cdot\   \frac{dW_1(\xi_1)dW_2(\xi_2)}{|\xi_1|^{\alpha-1/2} |\xi_2|^{\alpha-1/2}}. \nonumber\\ \EEA

\end{Lemma}

Two lines of computations show immediately that the variance of ${\cal A}^{\pm}_{\partial}$
is always finite (essentially, considering e.g. ${\cal A}^+_{\partial}$,
 because $\int_1^{+\infty} \frac{d\xi_1}{\xi_1^{2\alpha+1}}
\int_{\xi_1}^{+\infty} \frac{d\xi_2}{\xi_2^{2\alpha+1}}<\infty$, while the artificial infrared divergence at $\xi_1=0$ disappears when Taylor expanding $e^{\II t\xi_1}-e^{\II s\xi_1}$). On the other hand, letting $\xi:=\xi_1+\xi_2$ and introducing an ultra-violet cut-off at $|\xi_2|=\Lambda\gg 1$, one may see
for instance ${\cal A}^+(t)$ as an inverse random Fourier transform of the integral
$\xi\mapsto\int_{|\xi-\xi_2|<|\xi_2|}^{\Lambda}
 \frac{dW_2(\xi_2)}{\xi_2} \frac{1}{|\xi-\xi_2|^{\alpha-1/2} |\xi_2|^{\alpha-1/2}}$,
whose variance diverges like $\int^{\Lambda} \frac{d\xi_2}{\xi_2^{4\alpha}}=O(\Lambda^{1-4\alpha})$
or $O(\ln\Lambda)$  in the limit $\Lambda\to\infty$
as soon as $\alpha\le 1/4$. Note that the ultraviolet divergence is in the region $|\xi_1|,|\xi_2|\gg
|\xi|$. Later on, we shall introduce the field $\sigma$, and this region will correspond
to a splitting of the vertex $\int \partial\phi_1(x)\phi_2(x)\sigma(x) dx=\int d\xi_1 d\xi_2 d\xi
\del(\xi_1+\xi_2+\xi=0) ({\cal F}\partial\phi_1)(\xi_1){\cal F}\phi_2(\xi_2){\cal F}\sigma(\xi)$
such that the momentum of the $\sigma$-field, $\xi$, is much lower than that of the two $\phi$-fields.  Note also that, had we not operated the Fourier normal ordering, the above integral would
have been infrared divergent near $\xi_2=0$.

\bigskip

It is apparent that the central r\^ole in this decomposition is played by the Fourier
projection operator $D({\bf 1}_{|\xi_1|<|\xi_2|})={\cal F}^{-1} \left( {\bf 1}_{|\xi_1|<|\xi_2|} \ \cdot\ {\cal F}(\ .\ ) \right).$ Since ${\cal A}_{\partial}^{\pm}$
are obtained by Fourier projecting $(B_1(t)-B_1(s))\phi_2(s)$, or $(B_2(t)-B_2(s))\phi_1(t)$, which are perfectly well-defined products of continuous
fields \footnote{apart from the spurious infra-red divergence (see above)}, it
was clear from the onset that these would be regular terms. Hence singularities come only
from the one-time quantity  ${\cal A}^{\pm}(t)$, which {\em does not} split into a
 product of first-order integrals.

\medskip

In order to go one step further towards our final formulation, let us consider instead
of the brute-force projection $D({\bf 1}_{|\xi_1|<|\xi_2|})$ the following projection
operators, which are more illuminating, and also analytically more robust, allowing for
the extension of our theory to more general multiscale Gaussian fields
as defined in section 2.

\begin{Definition}[Fourier projections] \label{def:1.3}
 Let $\psi_i=\psi_i(x)$, $i=1,2$ be Gaussian fields, and $\psi_i=\sum_{j=0}^{\infty} \psi_i^j$ their scale decomposition along the Fourier partition of unity $(\chi^j)$ as in Definition \ref{def:Cjphi}. Then one lets
  \BEQ D^+(\psi_1\otimes\psi_2):= \sum_{0\le j<k<\infty} \psi^j_1\otimes \psi_2^k
+\half\sum_{j=0}^{+\infty}\psi_1^j\otimes \psi^j_2 \EEQ
and similarly
\BEQ D^-(\psi_1\otimes\psi_2):=\sum_{0\le j<k<\infty} \psi^k_1\otimes\psi^j_2
+\half\sum_{j=0}^{+\infty}\psi^j_1\otimes\psi^j_2. \EEQ
Clearly $\left(D^+ + D^-\right)(\psi_1\otimes\psi_2)=\psi_1\otimes\psi_2$.
\end{Definition}

Thus one abandons -- to the great benefit of notations -- the above complicated expression
for ${\cal A}^{\pm}(t)$ in favour of the analytically equivalent
\BEQ {\cal A}^+(t):=D^+\left( \int^t \partial \phi_1(t_1)\phi_2(t_1)dt_1 \right)
=\int^t \left[ \sum_{0\le j<k<\infty} \partial \phi_1^j \phi_2^k+\half\sum_{j=0}^{+\infty}
\partial \phi_1^j \phi_2^j \right](t_1) dt_1, \label{eq:1.12} \EEQ
\BEQ {\cal A}^-(t):=D^-\left( \int^t \partial \phi_2(t_1)\phi_1(t_1)dt_1 \right)
=\int^t \left[ \sum_{0\le j<k<\infty} \partial \phi_2^j \phi_1^k+\half\sum_{j=0}^{+\infty}
\partial \phi_2^j \phi_1^j \right](t_1) dt_1.  \label{eq:1.13} \EEQ

\subsection{Definition of the interaction}

We recall that $\int_s^t dB_1(t_1)\int_s^{t_1}dB_2(t_2)$ represents the
area between the straight line connecting $(B_1(s),B_2(s))$ to $(B_1(t),B_2(t))$ and the graph. If the curve turns right, resp. left, then the
L\'evy area increases, resp. decreases. We have seen that  ${\cal A}^{\pm}$ represents in some
sense the {\em singular part} of the L\'evy area.  So ${\cal A}^+$,
resp. ${\cal A}^-$ plays the r\^ole of a {\em right-, resp. left-turning} (singular)
{\em field}.

\medskip

It is conceivable that $B_1,B_2$ or $\phi_1,\phi_2$ represent the idealized, strongly
self-correlated motion in $\R^2$ of a particle, which -- although rotation-invariant --
may not (probably as a consequence of a mechanical or electromagnetic rigidity
due to the macroscopic dimension of the particle, or any other similar phenomenon)
turn absolutely freely. A natural quantum field theoretic description of this rigidity
phenomenon is to add an interaction Lagrangian of the form ${\cal L}_{int}=(\partial
A^{\pm})^2$.  The fundamental intuition here is that the field $B$ is in some sense a {\em mesoscopic field},
while ${\cal A}^{\pm}$ depends on {\em microscopic details of the theory}. This is explained in great accuracy in
\cite{Lej-bis}, in a mathematical language.
 A. Lejay explains how the paths of $B$ may be modified by inserting microscopic bubbles along the paths of $B$,
resulting in the limit in  paths which are indistinguishable from the original ones, while the L\'evy area has been corrected
by an {\em arbitrary} amount.  Hence one must search for an interaction which cures the ultra-violet divergences
of the microscopic scale,  without modifying the theory at mesoscopic scale. Understanding the mesoscopic scale as
low-frequency (which is not really appropriate, since there are {\em two} reference scales here, instead of one),
a natural candidate would be  an interaction theory which is asymptotically free at large (mesoscopic scale) distances. The
best-known example of that is probably the infra-red $\phi^4$-theory in $4$ dimensions \cite{FMRS}; in that case, the
  coupling constant increases indefinitely at short distances, and forces an ultra-violet cut-off. In our case, the
coupling constant $\lambda$ will not flow, and the theory will be well-defined at {\em all} scales, which suggests
a  {\em just renormalizable}
theory (or, in other terms, an integrated  interaction which is homogeneous of degree $0$).
 Since $(\partial A^{\pm})^2$ is homogeneous of
degree $(4\alpha-2)$ in time, one shall use in fact a {\em non-local} interaction
lagrangian, $\half c'_{\alpha} \int\int |t_1-t_2|^{-4\alpha} {\cal L}_{int}(\phi_1,\phi_2)(t_1,t_2)$, where
\BEQ {\cal L}_{int}(\phi_1,\phi_2)(t_1,t_2)= \lambda^2\left\{ \partial {\cal A}^+(t_1)\partial {\cal A}^+(t_2)+ \partial {\cal A}^-(t_1)\partial {\cal A}^-(t_2)\right\},\EEQ
which is {\em positive} for $\alpha<1/4$ since the kernel $|t_1-t_2|^{-4\alpha}$ is
locally integrable  and positive definite. Thus the Gaussian measure is {\em penalized} by
the singular exponential weight $e^{-\frac{c'_{\alpha}}{2} \int\int {\cal L}_{int}(\phi_1,\phi_2)(t_1,t_2) |t_1-t_2|^{-4\alpha}
dt_1 dt_2}$. Equivalently, using the so-called Hubbard-Stratonovich transformation
\footnote{which is nothing else but an infinite-dimensional extension of the Fourier
transform $\esper[e^{\II\lambda X}]=e^{-\sigma^2\lambda^2/2}$ for  a random variable
$X\sim {\cal N}(0,\sigma^2)$}, we introduce two independent exchange particle fields
$\sigma_{\pm}=\sigma_{\pm}(t)$ with covariance kernel $\esper \sigma_{\pm}(s)
\sigma_{\pm}(t)=c'_{\alpha} |s-t|^{-4\alpha}$ and rewrite (letting $d\mu(\phi)$, resp.  $d\mu(\sigma)$ be the Gaussian
measure associated to $\phi$, resp.  $\sigma=(\sigma_+,\sigma_-)$)
the partition function $Z:=Z(\lambda)$,
\BEQ Z:=\int e^{-\frac{c'_{\alpha}}{2} \int\int_{\R^2} |t_1-t_2|^{-4\alpha} {\cal L}_{int}(\phi_1,\phi_2)(t_1,t_2)
dt_1 dt_2} d\mu(\phi)\EEQ as
\BEQ Z:=\int  e^{- \int_{\R} {\cal L}_{int}(\phi_1,\phi_2,\sigma)(t) dt}
d\mu(\phi) d\mu(\sigma),\EEQ where
 \BEQ {\cal L}_{int}(\phi_1,\phi_2,\sigma)(t)=\II\lambda\left(\partial A^{+}(t)\sigma_+(t)
-\partial {\cal A}^-(t)\sigma_-(t)\right). \label{eq:Lint} \EEQ

\medskip All of this is ill-defined mathematically since (1) $\sigma$ is a distribution-valued process and $\partial A^{\pm}$ is not defined {\em at all} when
$\alpha\le 1/4$; (2) one integrates over $\R$ a translation-invariant quantity (note
that $\phi_1,\phi_2,\sigma$ are all {\em stationary fields}).


\subsection{Towards a heuristic expression for the L\'evy area}


Assume the coupling parameter $\lambda$ is small enough. Perturbative quantum field theory
suggests then to expand formally the exponential of the Lagrangian and compute polynomial
moments, $\frac{1}{Z}
\esper\left[ \psi_1(x_1)\ldots\psi_n(x_n) e^{- \int {\cal L}_{int}(\phi_1,\phi_2,\sigma)(t) dt}
\right]$,
 also called {\em $n$-point functions} and
denoted by  $\langle \psi_1(x_1)\ldots\psi_n(x_n)\rangle_{\lambda}$, $\psi_i=\phi_1,\phi_2,
\sigma_+$ or $\sigma_-$, as $\frac{1}{Z}\sum_{n\ge 0} \frac{(-1)^n}{n!}
\esper\left[ \psi_1(x_1)\ldots\psi_n(x_n) \left( \int {\cal L}_{int}(\ .\ ;t) dt\right)^n \right]$.
Using Wick's formula (recalled in \S 6.1), one may represent
 $\langle \psi_1(x_1)\ldots\psi_n(x_n)\rangle_{\lambda}$ as a sum over Feynman diagrams, $\sum_{\Gamma} A(\Gamma)$, where $\Gamma$ ranges over a set of diagrams with $n$
external legs, and $A(\Gamma)\in\C$ is the evaluation of the corresponding diagram (see
section 4). The Gaussian integration by parts formula \footnote{an infinite-dimensional
extension of the well-known formula for Gaussian vectors, $\esper \left[
\partial_{X_i}F(X_1,\ldots,X_n)\right]=\sum_j C^{-1}(i,j) \esper \left[ X_j F(X_1,\ldots,X_n)\right]$ if
$C$ is the covariance matrix of $(X_1,\ldots,X_n)$.} yields a so-called {\em Schwinger-Dyson identity},
\BEA  \langle \partial A^{\pm}(x)\partial A^{\pm}(y)\rangle_{\lambda}
&=&-\frac{1}{\lambda^2 Z(\lambda)} \esper\left[ \frac{\del}{\del \sigma_+(y)} \frac{\del}{\del\sigma_+(x)}
e^{-\int {\cal L}_{int}(\phi_1,\phi_2,\sigma_+)(t)dt} \right]\nonumber\\
&=& -\frac{1}{\lambda^2 Z(\lambda)} \esper\left[ (C_{\sigma_+}^{-1}\sigma_+)(y) \frac{\del}{\del\sigma_+(x)}
e^{- \int {\cal L}_{int}(\phi_1,\phi_2,\sigma_+)(t)dt} \right] \nonumber\\
&=& -\frac{1}{\lambda^2} \left[ -C_{\sigma_+}^{-1}(x,y)+ \langle (C_{\sigma_+}^{-1}\sigma_+)(x)
(C_{\sigma_+}^{-1}\sigma_+)(y)\rangle_{\lambda}\right], \nonumber\\ \EEA
with Fourier transform
\BEQ  \langle |{\cal F}(\partial A^{\pm})(\xi)|^2\rangle_{\lambda}=\frac{1}{\lambda^2}
|\xi|^{1-4\alpha} \left[ 1-|\xi|^{1-4\alpha}  \langle |({\cal F}\sigma_+)(\xi)|^2\rangle_{\lambda} \right].
\EEQ

\begin{figure}[h]
  \centering
   \includegraphics[scale=0.3]{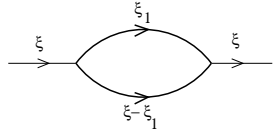}
   \caption{Bubble diagram. Bold lines are $\phi$-fields, plain lines are $\sigma$-fields.}
  \label{Fig-bubble}
\end{figure}

\begin{figure}[h]
  \centering
   \includegraphics[scale=0.4]{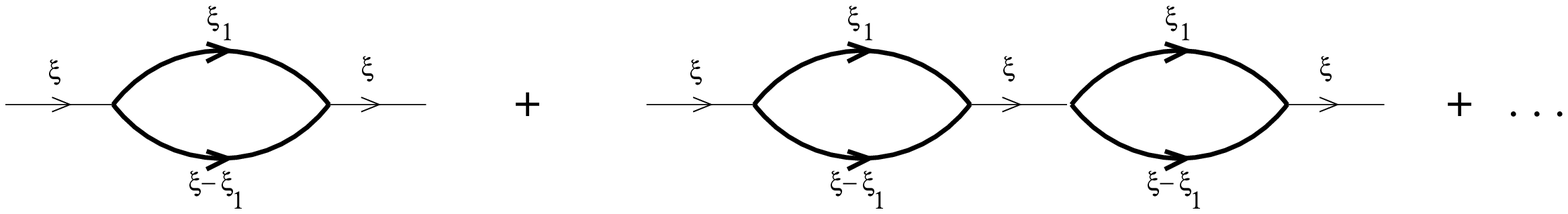}
   \caption{First two terms of the bubble series.}
  \label{Fig-bubble-series}
\end{figure}

Consider (using a brute-force ultraviolet cut-off at momentum $\Lambda$)
the term of lowest degree in $\lambda$ in the term between square brackets: it
is essentially equal, up to a sign, to the evaluation of the half-amputated bubble diagram, see Fig. \ref{Fig-bubble},
\BEA &&  - |\xi|^{1-4\alpha}\ \cdot\ (-\II\lambda)^2  \int^{\Lambda}_{|\xi_1|<|\xi-\xi_1|} d\xi_1
\left\{  \left( \esper[ |{\cal F}\sigma_+(\xi)|^2 ]  \right)^2 \
\esper[ |{\cal F}(\partial\phi_1)(\xi_1)|^2 ]\  \esper[ |{\cal F}\phi_2(\xi-\xi_1)|^2 ] \right\} \nonumber\\
&&=  \lambda^2 |\xi|^{4\alpha-1}  \int_{|\xi_1|<|\xi-\xi_1|}^{\Lambda} d\xi_1 |\xi_1|^{1-2\alpha}
|\xi-\xi_1|^{-1-2\alpha} \sim_{\Lambda\to\infty} K\lambda^2 (\Lambda/|\xi|)^{1-4\alpha},\EEA
hence a diverging {\em positive} quantity. (Using the Fourier truncation of Definition
\ref{def:1.3} only changes the constant $K$.)  However, resumming {\em formally} the bubble series as in
Fig. \ref{Fig-bubble-series}
yields -- taking into account the possible insertion of $\sigma_-$-propagators between $\sigma_+$-propagators --
\BEA  \frac{1}{\lambda^2}|\xi|^{1-4\alpha} \left[1-\frac{1}{1+K'\lambda^2 (\Lambda/|\xi|)^{1-4\alpha}}\right] &=&\frac{1}{\lambda^2}|\xi|^{1-4\alpha}  \ \cdot\
\frac{K'\lambda^2 (\Lambda/|\xi|)^{1-4\alpha}}{1+K'\lambda^2 (\Lambda/|\xi|)^{1-4\alpha}} \nonumber\\
&\to_{\Lambda\to\infty} & \frac{1}{\lambda^2} |\xi|^{1-4\alpha}. \label{eq:1.21} \EEA
Thus the bare $\sigma$-propagator $\frac{1}{|\xi|^{1-4\alpha}}$ has been replaced with the renormalized propagator $\frac{1}{|\xi|^{1-4\alpha}+K\lambda^2 \Lambda^{1-4\alpha}}$, which vanishes in the limit $\Lambda\to\infty$.
In physical terms, the interaction in $\frac{1}{|\xi|^{1-4\alpha}}$ has been totally screened by an infinite
mass term $K'\lambda^2 \Lambda^{1-4\alpha}.$ More complicated diagrams contributing to $\langle |({\cal F}\sigma_+)(\xi)|^2\rangle_{\lambda}$, and involving internal $\sigma$-links, also vanish when $\Lambda\to\infty$. There remains simply:

\BEQ \langle |{\cal F} {\cal A}^{\pm}(\xi)|^2|\rangle_{\lambda}=\frac{1}{\lambda^2}
|\xi|^{-1-4\alpha}. \EEQ

As for the mixed term $\langle \partial {\cal A}^{\pm}(x)\partial{\cal A}^{\mp}(y)\rangle_{\lambda}$, its
Fourier transform is given by $\frac{1}{\lambda^2}|\xi|^{1-4\alpha} \left[-\frac{1}{1+K''\lambda^2 (\Lambda/|\xi|)^{1-4\alpha}}\right]$, where $K''<K'$ due to the constraints on the scales for
bubbles of mixed type with one $\sigma_+$- and one $\sigma_-$-leg \footnote{Namely, coupling $(\phi_1,\phi_2)$ simultaneously to $\sigma_+$ and $\sigma_-$ leaves only the scale-diagonal projected
field $D^+ D^- (\phi_1\otimes\phi_2)=\frac{1}{4}\sum_j \phi_1^j\otimes\phi_2^j$.}, which vanishes
in the limit $\Lambda\to\infty$ (note the disappearance of the factor $1$ compared to eq. (\ref{eq:1.21}), due to the fact that $\esper \sigma^+(x)\sigma^-(y)=0$). Thus the covariance of the two-component
$\sigma$-field has been renormalized to $\frac{1}{|\xi|^{1-4\alpha}\Id+b^{\rho}}$, where $b^{\rho}$ is
a two-by-two positive matrix with eigenvalues $\thickapprox \lambda^2 M^{\rho(1-4\alpha)}$.

\medskip

Using eq. (\ref{eq:1.8}), one obtains:
\BEA (2\pi c_{\alpha})^2 \langle {\cal A}(s,t)^2 \rangle_{\lambda} &=& \langle \left| {\cal A}^+(t)-
{\cal A}^+(s)\right|^2 \rangle_{\lambda}+\langle \left| {\cal A}^-(t)-
{\cal A}^-(s)\right|^2 \rangle_{\lambda} \nonumber\\
&& \qquad \qquad \qquad +\esper \left|{\cal A}^+_{\partial}(s,t)-
{\cal A}_{\partial}^-(s,t)\right|^2 \nonumber\\
&=&   \frac{4}{\lambda^2} \int (1-\cos(t-s)\xi)|\xi|^{-1-4\alpha} d\xi+
 \esper \left|{\cal A}^+_{\partial}(s,t)-
{\cal A}_{\partial}^-(s,t)\right|^2 \nonumber\\
&=& (\frac{4}{\lambda^2}K_1+K_2)|t-s|^{4\alpha} \EEA
for some constants $K_1,K_2$.

\begin{figure}[h]
  \centering
   \includegraphics[scale=0.4]{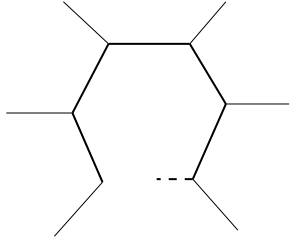}
   \caption{Higher connected moments of the L\'evy area.}
  \label{Fig-polygonal}
\end{figure}

A simple power-counting argument (see subsection 4.1) shows that the overall degree of divergence of a connected diagram with $2n$ external $\sigma$-legs is
$1-4n\alpha$. For $n\ge 2$, this is $\le 1-8\alpha<0$ since $\alpha>\frac{1}{8}$ by hypothesis, so such diagrams are convergent. By the above arguments, there remain only the connected diagrams in the limit $\Lambda\to\infty$,  see Fig. \ref{Fig-polygonal}, whose evaluation is independent of $\lambda$.
Then these may be shown by standard multi-scale arguments for single Feynman diagrams (see e.g. \cite{FVT}) to be of order $K^n$, which implies that the generating series for connected moments has
a non-zero radius of convergence.

On the other hand, the law of the field $\phi$ is left unchanged by the interaction. Namely,
 all non-trivial diagrams contributing e.g.  to $\langle \phi_1(x)\phi_2(x)\rangle_{\lambda}$ involve internal $\sigma$-links which (as previously "shown") vanish in the limit $\Lambda\to\infty$.

On the whole, this is the content of Theorem 0.1.

\bigskip

The art of constructive field theory is to make the previous speculations rigorous.


\section{Multiscale Gaussian fields}



\subsection{Scale decompositions}


We fix some constant $M>1$. The next definition is borrowed from \cite{Trie}.

\begin{Definition}[Fourier partition of unity]

 Let $\chi^0:\R\to[0,1]$ be an even, $C^{\infty}$ function such that $\chi^0\big|_{[-1,1]}\equiv 1$ and $\supp \chi^0\subset
[-M,M]$, and
\BEQ \chi^j(\xi):=\chi^0(M^{-j}\xi)-\chi^0(M^{1-j}\xi),\quad j\ge 1\EEQ
so that $\supp \chi^j\subset [M^{j-1},M^{j+1}]\cup[-M^{j+1},-M^{j-1}].$

The $(\chi_j)_{j\ge 0}$ define a $C^{\infty}$ partition of unity, namely,
$\sum_{j=0}^{+\infty} \chi_j\equiv 1$.

\end{Definition}

Note that $\supp \chi_j\cap\supp \chi_{j'}$ has empty interior if $|j-j'|=2$, and is empty if
$|j-j'|\ge 3$.

\bigskip

\begin{Definition}[$M$-adic intervals]

Let $\D^j:=\{ [kM^{-j},(k+1)M^{-j}), k\in\Z\}$, $j\ge 0$ be the set of $M$-adic intervals of scale $j$, and
$\D:=\uplus_{j\ge 0} \D^j$ the disjoint sum of these sets over all  scales.
The set $\D$ is a tree with links (called: {\em inclusion links}) connecting each interval $\Del\in\D^j$ to the unique interval $\Del'\in\D^{j-1}$ such that
$\Del\subset\Del'$ (see Definition \ref{def:polymer} below, and left part of Fig. \ref{Fig-polymer} in
subsection 3.3). An element of $\D^j$ is usually
denoted by $\Del^j$, or simply $\Del$ if no confusion may arise. The volume $|\Del^j|$ is simply $M^{-j}$. If $\Del\in \D^j$, then one
denotes by $j(\Del)=j$ the scale of $\Del$.

If $x\in \R$, then $x$ belongs to a single $M$-adic interval of scale $j$, denoted by $\Del^j_x$.

If $\Del^j\in\D^j$, then
the set of intervals $\Del\in\uplus_{h<j} \D^h$ such that $\Del$ lies below $\Del^j$, i.e. $\Del\supset\Del^j$,
is denoted by $(\Del^j)^{\downdownarrows}$.

We denote by $d^j(\Del,\Del')$, $\Del,\Del'\in \D^j$, the distance in terms of number of $M$-adic intervals of scale $j$
between $\Del$ and $\Del'$, namely,
\BEQ d^j([kM^{-j},(k+1)M^{-j}),[k'M^{-j},(k'+1)M^{-j}))=|k'-k|.\EEQ

By extension, one may also define the $d^j$-distance of two points or two
$M$-adic intervals of scale $j'>j$, namely, $d^j(x,y)=M^j|x-y|$ and
\BEQ d^j(\Del,\Del'):=M^{j-j'} d^{j'}(\Del,\Del'),\quad \Del,\Del'\in\D^{j'}.\EEQ

\end{Definition}

{\bf Remark.} It is preferable {\em not} to define $d^j(\Del,\Del')$ for $\Del,\Del'\in\D^{j'}$ with
$j'<j$, since $d^j(x,x')$, $x\in\Del,x'\in\Del'$ depends strongly on the choice of the points $x,x'$ then.

The definition extends in a natural way to a $D$-dimensional setting by decomposing
$\R^D$ into a disjoint union of hypercubes of size side $M^{-j}$.

\begin{Definition}[multiscale decomposition]  \label{def:multi-scale-decomposition}
\begin{itemize}

\item[(i)] Let $\psi=\psi(x)$ be some integrable function or distribution in ${\cal S}'$.
Then the multiscale decomposition of $\psi$ is $\psi=\sum_{j\ge 0} \psi^j$, where
\BEQ \psi^j(x):={\cal F}^{-1}\left( \xi\mapsto \chi^j(\xi)({\cal F}\psi)(\xi)\right)(x)=\frac{1}{\sqrt{2\pi}}
\int e^{\II x\xi} \chi^j(\xi) ({\cal F}\psi)(\xi) \ d\xi.\EEQ

\item[(ii)] (low-momentum field)
Let $\psi^{\to k}:=\sum_{j=0}^k \psi^j.$

\item[(iii)] (high-momentum field)
Let $\psi^{k\to}:=\sum_{j=k}^{+\infty} \psi^j.$
\end{itemize}
\end{Definition}

Now comes a general remark. Let $\hat{f}= {\cal F}{f}(\xi)$ be some function with
support in $|\xi|\le M^j$ such that $|{\cal F}(f')(\xi)|=|\xi \hat{f}(\xi)|$ is
bounded. Then
\BEQ |{\cal F}^{-1} \hat{f}(x)-{\cal F}^{-1}\hat{f}(y)|\le |x-y| \int_0^{M^j} |\xi|
|\hat{f}(\xi)| \ d\xi \le K\cdot M^j |x-y|,\EEQ

so $\psi^j(x)$ or $\psi^{\to j}(x)$ varies slowly inside intervals of scale $k$ if $k>j$. Hence it makes
sense in first approximation to consider $\psi^{j}(x)$ or $\psi^{\to j}(x)$
to be approximately equal to the averaged, locally constant function
\BEQ \psi^{j}_{av}(x):=\sum_{\Del^j\in \D^j} {\bf 1}_{x\in \D^j} \frac{1}{|\Del^j|} \int_{\Del^j} \psi^j(y)dy
=\frac{1}{|\Del^j_x|} \int_{\Del^j_x} \psi^{j}(y)dy, \EEQ
or similarly for the low-momentum field $\psi^{\to j}$
\BEQ \psi^{\to j}_{av}(x):=\sum_{\Del^j\in \D^j} {\bf 1}_{x\in \D^j} \frac{1}{|\Del^j|} \int_{\Del^j} \psi^{\to j}(y)dy=\frac{1}{|\Del^j_x|} \int_{\Del^j_x} \psi^{\to j}(y)dy. \EEQ
  Summing the $\psi^j_{av}$
over $j$ would give a new function $\sum_{j\ge 0} \psi^{j}_{av}(x)$ which is a sort of ``naive'' wavelet expansion of the original
function $\psi$; with a little extra care, one could arrange that the two functions be equal, but we shall not need to do
so.

Conversely, if $\psi$ is a 'reasonable' random Gaussian field, then the
covariance
$\langle \psi^j(x)\psi^j(y)\rangle$ -- or, more generally, $\langle\psi^{j\to}(x)\psi^{j\to}(y)\rangle$  --
is usually small if the corresponding $M$-adic intervals are far apart, i.e. if
$d^j(\Del^j_x,\Del^j_y)\gg 1$.

These two remarks may be made precise in the case when $\psi$ is a multiscale Gaussian
field, see Definition \ref{def:multiscale-Gaussian-field} below. Then, for any  scale $j$,

\begin{itemize}
\item the field $\psi^{\to j}$ may be decomposed into the sum of the {\em locally averaged field at scale $j$}, namely,
 $\psi^{\to j}_{av}(x)$,
and a {\em secondary field}, denoted by $\del^j\psi^{\to j}$, whose low momentum components of scale $h<j$ decrease like $M^{-\gamma(j-h)}$ for
some $\gamma>0$,
see Lemma \ref{lem:2.6} and Corollary \ref{cor:spring-factor}. For reasons explained below, it is customary
to use the nickname of {\em spring factor} for a decrease factor of the type $M^{-\gamma(j-h)}$, $\gamma>0$.

\item since the covariance  decreases with the distance in terms of number of $M$-adic intervals of scale $j$, it makes sense to
try and find some expansion of the functional  ${\cal L}(\psi)$ in which the field values over far enough $M$-adic intervals have
been made independent. This is called a {\em cluster expansion} (see section 2).
\end{itemize}

\begin{Definition}[multiscale Gaussian field]  \label{def:multiscale-Gaussian-field}

A {\em multiscale Gaussian field} with scaling dimension $\beta<1$ is a field $\psi=\psi(x)$ such that,
for every $r\ge 0$ and $\tau,\tau'=0,1,2,\ldots$,
\BEQ |\langle \partial^{\tau} \psi^j(x) \partial^{\tau'}\psi^j(y)\rangle|\le K_{\tau,\tau',r}
\frac{M^{(\tau+\tau'-2\beta)j}}{(1+M^j|x-y|)^r}\EEQ
with some constant $K_{\tau,\tau',r}$ depending only on $\tau,\tau'$ and $r$.
\end{Definition}

{\bf Remark.} A multiscale Gaussian field $\psi$ with scaling dimension $\beta\in(0,1)$
 has almost surely $\beta^-$-H\"older paths. Namely, $\esper
|\psi^j(x)-\psi^j(y)|^2=\int_x^y \int_x^y dz dz' \langle \partial \psi^j(z)
\partial \psi^j(z')\rangle$ is bounded by $K(x-y)^2 M^{(2-2\beta) j}$ if $|x-y|<M^{-j}$, and by
\BEQ 2\int_{-|x-y|}^{|x-y|} dz' |\langle \psi^j(0)\partial\psi^j(z')\rangle|\le K \int_0^{|x-y|}dz' \frac{M^{(1-2\beta)j}}{(1+M^jz')^2}\le K' M^{-2\beta j} \EEQ
otherwise. Summing over $j$ yields
\BEA  && \esper |\psi(x)-\psi(y)|^2\le K\left ((x-y)^2 \sum_{j\le-\log|x-y|}
M^{(2-2\beta)j}+\sum_{j\ge -\log|x-y|} M^{-2\beta j} \right) \nonumber\\
&& \qquad \qquad \qquad \le K'|x-y|^{2\beta},\EEA and one concludes by using the Kolmogorov-Centsov lemma.

\bigskip

As we shall see in the next paragraph, fractional Brownian motion with Hurst index $\alpha$ is the paramount example of multiscale
Gaussian field with scaling dimension $\beta=\alpha\in(0,1)$.

\begin{Definition}[averaged and secondary fields] \label{def:2.5}

 Choose some scale $k\ge 0$. Let $f:\R\to\R$ be a function, typically, $f(x)=\psi^{\to k}(x)$ or $\psi^j(x)$ for some
$j<k$. Then:
\begin{itemize}
\item[(i)] the {\em averaged field at scale $k$} is the locally constant field
\BEQ f_{av}(x):=\sum_{\Del^k\in \D^k} {\bf 1}_{x\in \D^k} \frac{1}{|\Del^k|}
 \int_{\Del^k} f(y)dy.\EEQ
It is more convenient  to consider $f_{av}$ as some function on $\D^k$, still denoted by $f$,
\BEQ f(\Del^k):=\frac{1}{|\Del^k|} \int _{\Del^k} f(y)dy, \quad \Del^k\in \D^k,\EEQ
so that $f_{av}(x)=f(\Del^k_x).$ This notation fixes unambiguously the scale of the average.

\item[(ii)] the {\em secondary field at scale $k$} is the difference between the original field and the averaged field at scale $k$,
namely,
\BEQ \del^k f(x):=f(x)-f(\Del^k_x).\EEQ
\end{itemize}
\end{Definition}

These definitions imply the following easy Lemma:

\begin{Lemma} \label{lem:2.6}

\BEQ \del^k f(x)=\frac{1}{|\Del^k|} \int_{\Del^k_x}  \left( \int_u^x f'(v)dv\right) du= \int_{\Del^k_x} dv f'(v) \del^k(x;v) \EEQ
where
\BEQ \del^k(x;v):=\frac{1}{|\Del^k|} \left\{ (v-\inf \Del_x^k){\bf 1}_{v<x}+(v-\sup \Del^k_x) {\bf 1}_{v>x}\right\} \in [-1,1] \EEQ
is a signed distance from $v$ to the boundary of the interval $\Del_x^k=[\inf\Del_x^k,\sup\Del_x^k)$ measured in terms
of the rescaled $d^k$-distance.
\end{Lemma}

{\bf Proof.} Straightforward. \hfill \eop

\begin{Corollary}  \label{cor:spring-factor}

Let $k,k'>j$. Assume $\psi=\psi(x)$ is a multiscale Gaussian field with scaling dimension $\beta<1$. Then
\BEA  && |\langle\del^k \psi^j(x)\del^{k'} \psi^j(y)\rangle| \le K_r \frac{M^{-2\beta j}}{(1+M^j|x-y|)^r}
M^{-(k-j)} M^{-(k'-j)} \label{eq:deltapsi1}\\
&& \qquad \qquad = K_r M^{-\beta (k+k')} \frac{M^{-(1-\beta)(k-j)} M^{-(1-\beta)(k'-j)} }{(1+M^j|x-y|)^r}. \label{eq:deltapsi2}\EEA

\end{Corollary}

{\bf Proof.}  Straightforward. \hfill \eop

Eq. (\ref{eq:deltapsi1}), resp. (\ref{eq:deltapsi2})
 emphasizes the ``spring-factor'' $M^{-(k-j)}$, resp. $M^{-(1-\beta)(k-j)}$ gained for each secondary field with respect to the
covariance of a multiscale Gaussian field with scaling dimension $\beta$ {\em at scale $j$}, resp.
{\em at scale $k$}. Had one considered directly $|\langle \psi^j(x)\psi^j(y)\rangle|$, the spring factor -- called {\em rescaling spring factor} -- in  (\ref{eq:deltapsi2}) would have been simply $M^{\beta(k-j)}$, see
introduction to \S 6.1.2.

{\bf Remarks.}
 \begin{enumerate}
 \item
 The same spring factors appear in $D$ dimensions. It is sometimes useful to
take a cleverer definition of the secondary field by replacing the simple average
$f(\Del_x^k)$ with a wavelet component of $f$, where the wavelets have vanishing first
moments up to order $\tau\ge 1$. This allows one to enhance the spring-factor from
$M^{-(k-j)}$ to $M^{-(\tau+1)(k-j)}$, resp. from $M^{-(1-\beta)(k-j)}$ to
$M^{-(\tau+1-\beta)(k-j)}$.

\item The separation of low-momentum fields into  a sum (field average)+(secondary field) is required only for fields with scaling dimension $\beta>-D/2$, and is {\em not}  performed otherwise
(see explanation after Definition \ref{def:3.10} and in subsection 6.1).
\end{enumerate}

\bigskip

Consider conversely a high-momentum field $\psi^j(x)$, $j>h$ with $x\in \D^h$.
Typically (see section 3), $\psi^j(x)=\psi^j(x_{\ell^h})$,
$\psi^j(x'_{\ell^h})$ or $\psi^j(x_{\Del^h,\tau})$, $x\in\Del^h$ coming from the scale $h$ horizontal
cluster expansion or a scale $h$ $t$-derivation in an interval $\Del^h\in\D^h$, and one must integrate over $x\in\Del^h$.
Then $\psi^j(x')$ and $\psi^j(x'')$ are almost decorrelated if $x',x''\in \Del^h$
but $d^j(x,x')\gg 1$. Hence it makes sense to restrict $\psi^j$ over each sub-interval $\Del^j\subset\Del^h$ of scale $j$:

\begin{Definition}[restriction of high-momentum fields] \label{def:restriction}

Let $\Del^h\in\D^h$ and $j>h$. Then the high-momentum field $\psi^j(x)$, $x\in\Del^h$, splits into
\BEQ \psi^j(x)=\sum_{\Del^j\in\D^j,\Del^j\subset\Del^h} Res^h_{\Del^j}\psi^j(x),\quad x\in\D^h\EEQ
where $Res^h_{\Del^j}\psi^j(x):={\bf 1}_{x\in\Del^j} \psi^j(x).$

\end{Definition}


\subsection{Multiscale Gaussian fields in one dimension}


We introduce here more specifically the infra-red divergent stationary field $\phi$ associated to fBm $B$, and the fields $\sigma=\sigma_{\pm}$ conjugate to the turning fields ${\cal A}^{\pm}$. In this paragraph $D=1$. All fields come implicitly with an ultra-violet cut-off at scale $\rho$, so that $\psi=\phi,\sigma$ should be understood as $\psi^{\to\rho}$, $\psi^{j\to}$ as
$\psi^{j\to\rho}$, and so on.

\begin{Definition}[Harmonizable representation of fBm]

Let $W(\xi),\xi\in\R$ be a complex Brownian motion such that $W(-\xi)=\overline{W(\xi)}$, and
\BEQ B_t:=(2\pi c_{\alpha})^{-\half} \int_{-\infty}^{+\infty} \frac{e^{\II t\xi}-1}{\II\xi} |\xi|^{\half-\alpha} dW(\xi),
\quad t\in\R.\EEQ
\end{Definition}

The field $B_t, t\in\R$ is called {\em fractional Brownian motion}
\footnote{The constant $c_{\alpha}$ is  conventionally chosen
 so that $\esper (B_t-B_s)^2=|t-s|^{2\alpha}$.} . Its paths are almost surely
$\alpha^-$ H\"older, i.e. $(\alpha-\eps)$-H\"older for every $\eps>0$. It has dependent
but identically distributed (or in other words, stationary) increments $B_t-B_s$.
In order to gain translation invariance, we shall rather use the closely related {\em stationary
process}
\BEQ \phi(t):= \int_{-\infty}^{+\infty} \frac{e^{\II t\xi}}{\II\xi} |\xi|^{\half-\alpha} dW(\xi),
\quad t\in\R\EEQ
-- with covariance
\BEQ \langle \phi(x)\phi(y)\rangle =  \int e^{\II \xi (x-y)} \frac{1}{|\xi|^{1+2\alpha}} d\xi
\label{eq:cov-phi} \EEQ
--
which is  infrared divergent. However, as already discussed in section 1, the increments $\phi(t)-\phi(s)$ are well-defined for any pair of variables $(s,t)$.

\begin{Definition} \label{def:Cjphi}
\begin{enumerate}
\item
Let
\BEQ C^j_{\phi}(x,y):=  \int e^{\II \xi (x-y)} \frac{\chi^j(\xi)}{|\xi|^{1+2\alpha}} d\xi,
\quad j\in\N.\EEQ
Then $ C_{\phi}:=\sum_{j\in\N} C^j_{\phi}$ is the covariance of the field $\phi$. We denote by
$\phi:=\sum_{j\in\N} \phi^j$ the corresponding multiscale decomposition of the field $\phi$ into independent
components $\phi^j$, $j\in\N$.
\item Let $\phi^{\to j}=\sum_{h=0}^j \phi^h$ and $\phi^{j\to}=\phi^{j\to\rho}=\sum_{k=j}^{\rho} \phi^k$. The covariance of $\phi$, resp. $\phi^{\to j}$, resp. $\phi^{j\to}$ is $D(\chi^j)C_{\phi}$, resp. $D(\chi^{\to j})C_{\phi}$, resp. $D(\chi^{j\to})C_{\phi}$, where $\chi^{\to j}:=\sum_{h=0}^j \chi^h$, $\chi^{j\to}=\sum_{k=j}^{\rho} \chi^k.$
    \end{enumerate}
\end{Definition}

{\bf Remark.} Note that this multiscale decomposition does not coincide exactly
with that of Definition \ref{def:multi-scale-decomposition}. With this slightly modified definition, the $\phi^j$ have been made independent.

\begin{Lemma} \label{lem:2.11}
The stationary field $\phi$ associated to fBm  is a multiscale Gaussian field with scaling dimension $\alpha$.
\end{Lemma}

{\bf Proof.} A simple scaling of the variable of integration yields
\BEQ  C_{\phi}^j(x,y)=M^{-2\alpha(j-1)} \int e^{\II M^{j-1}(x-y)\xi'}
 \frac{1}{|\xi'|^{1+2\alpha}} \chi^1(\xi') d\xi'\EEQ

with $\supp\chi^1\subset[1,M^2]\cup[-M^2,-1]$ bounded away from $0$, hence $|C_{\phi}^j(x,y)|\lesssim M^{-2\alpha j}$, and more precisely (by integrating by parts
$n$ times)
$ |C_{\phi}^j(x,y)|=O\left(M^{-2\alpha j} (M^j |x-y|)^{-n}\right)$ when $M^j |x-y|\to\infty$.
The bound for $\langle \partial^{\tau} \phi^j(x)\partial^{\tau'} \phi^j(y)\rangle$ is obtained in the same way
(simply multiply by $|\xi|^{\tau+\tau'}$ in Fourier coordinates).
\hfill \eop

\bigskip

We may now define the field $\sigma$. As indicated in the Introduction and in section 1, the Fourier transform of
 its {\em bare} covariance $\langle \sigma_i\sigma_{i'}\rangle$, $i,i'=\pm$ is $\frac{\del_{i,i'}}{|\xi|^{1-4\alpha}}$. On the other hand, the renormalized covariance (up to
 the overlap between the supports of the Fourier multipliers $\chi^j$)
is essentially $\frac{1}{|\xi|^{1-4\alpha}\Id+\sum_{j=0}^{\to\rho} b^j \chi^{\to(j-1)}(\xi)}$, where the scale $j$
mass  counterterm $b^j$ for diagrams of scale $j$ with two low-momentum external legs $\sigma^{\to(j-1)}$ -- a two-by-two, positive matrix -- is defined
 inductively (see subsection 3.5, section 4 and subsection 6.3) and shown to be of order
 $\lambda^2 M^{j(1-4\alpha)}$. For technical reasons one chooses to retain in the covariance of $\sigma$ only
a simplified version of the counterterm (essentially the term of highest scale $\rho$), which is of the
same order as the sum of all mass counterterms.

\begin{Definition}
\begin{itemize}
\item[(i)]
Let $\sigma$ be the {\em stationary} two-component massive Gaussian field with covariance kernel
\BEQ C^{\to\rho}_{\sigma}(x-y):=\int_{-\infty}^{+\infty} \frac{e^{\II \xi(x-y)}}{|\xi|^{1-4\alpha}\Id+b^{\rho}}
\ \chi^{\rho}(\xi) \ d\xi.\EEQ
\item[(ii)] Decompose $C^{\to\rho}_{\sigma}$ into $\sum_{j\ge 0} C_{\sigma}^j$ and $\sigma$ into a sum of
independent fields, $\sum_{j\ge 0} \sigma^j$ as in Definition \ref{def:Cjphi} by setting
\BEQ C^j_{\sigma}(x,y)=\int e^{\II\xi(x-y)} \frac{\chi^j(\xi)}{|\xi|^{1-4\alpha}\Id+b^{\rho}} \ d\xi.\EEQ
\end{itemize}
\end{Definition}

\begin{Lemma} \label{lem:2.13}
$\sigma$ is a multiscale Gaussian field with scaling dimension $-2\alpha$. More precisely,
\BEQ \left| \langle \partial^{\tau}\sigma^j(x)\partial^{\tau'}\sigma^j(y)\rangle\right|\le K_{\tau,\tau',r}
\frac{M^{(\tau+\tau'-4\alpha)j}}{(1+M^j |x-y|)^r} \ .\ \inf(1,\frac{M^{j(1-4\alpha)}}{b^{\rho}}).\EEQ
\end{Lemma}

For $\rho$ large enough, $\inf(1,\frac{M^j(1-4\alpha)}{b^{\rho}})\thickapprox \lambda^{-2} M^{-(\rho-j)(1-4\alpha)}$.
Thus $\sigma$ vanishes in the limit $\rho\to\infty$ because of the infinite mass counterterm.

{\bf Proof.} Same as for Lemma \ref{lem:2.11}. Note that the rescaled denominator
$\frac{\chi^1(\xi)}{|\xi|^{1-4\alpha}+b^{\rho} M^{-(j-1)(1-4\alpha)}}$ is bounded by
$\inf(1,\frac{M^j(1-4\alpha)}{b^{\rho}})$. \hfill\eop

\bigskip

Note the following to conclude:

-- multiscale Gaussian fields with {\em positive scaling dimension} (as $\phi$ for instance)
 are ultra-violet convergent and define continuous fields, but may be  infra-red divergent, necessitating
some infra-red cut-off (however, if $\beta<1$,  the increments $\phi_t-\phi_s$
 are well-defined);

-- multiscale Gaussian fields with {\em negative scaling dimension} (as in the case of the field $\sigma$,
or of most fields in quantum field theory) are ultra-violet divergent and define distribution-valued fields
whose covariance kernel is singular on the diagonal. On the other hand, they are infra-red convergent.


\section{Cluster expansions: an outline}


The general aim of horizontal and vertical cluster expansions has already been explained in the Introduction and
in \S 2.1. The
 (non restricted) {\em horizontal cluster expansion} has been given by D. Brydges and T. Kennedy a beautiful combinatorial
structure in terms of forests of intervals (see \S 3.1 and 3.2). The vertical or {\em momentum-decoupling expansion}, in
 terms of $t$-parameters, is somewhat looser, relying on a Taylor expansion to some order $\tau_{\Del}$ in each interval
 $\Del$. Putting together these two expansions, one obtains so-called {\em polymers} (see \S 3.3). Polymers with too few
 external legs must still be renormalized by resumming the local part of diverging graphs into an exponential making up an interaction counterterm  (see section 4 for detailed explanations); a Mayer expansion makes it possible  to get
rid of their non-overlap constraints with the other polymers, and finally  to divide out the
 polymers with no external legs (called: {\em vacuum polymers})
 when computing $n$-point functions (see \S 3.4).  All these provide series depending
on  the
small parameter $\lambda$ which will be shown to converge in section 6.

Summarizing, the general expansion scheme which one must have in mind is made up
of  a sequence of transformations for each scale, starting from highest scale $\rho$, namely

( {\tiny  Horizontal cluster expansion at scale $\rho$ $\longrightarrow$
Vertical cluster expansion at scale $\rho$ $\longrightarrow$ separation of local part of diverging graphs of lowest scale $\rho$ $\longrightarrow$ Mayer expansion at scale $\rho$
$\longrightarrow$ resummation of local parts of diverging graphs of lowest scale $\rho$, also called : renormalization stage of scale $\rho$ \Big)
$\longrightarrow \ldots\longrightarrow$
\Big(Horizontal cluster expansion at scale $j$ $\longrightarrow$
Vertical cluster expansion at scale $j$ $\longrightarrow$ separation of local part of diverging graphs of lowest scale $j$ $\longrightarrow$ Mayer expansion at scale $j$
$\longrightarrow$ resummation of local parts of diverging graphs of lowest scale $j$ \Big)$\longrightarrow\ldots$ }

down to lowest scale $j=0$.


\subsection{The general Brydges-Kennedy formula}


Let us start with the following general definition.

\begin{Definition}
Let:
\begin{itemize}
\item[(i)] $\cal O$ be an arbitrary set, whose elements are called {\em
objects};
\item[(ii)] $L({\cal O})$ be the set of links of the total graph
associated to $\cal O$,
or in other words, the set of pairs of objects, so that $\ell\in L({\cal
O})$ is
represented as a pair, $\ell\sim\{o_{\ell}, o'_{\ell}\}\subset {\cal
O}$,
$o_{\ell}\not= o'_{\ell}$;
\item[(iii)] $[0,1]^{L({\cal O})}:=\{z=(z_{\ell})_{\ell\in L({\cal O})},
0\le z_{\ell}\le 1 \} $ be the convex set of {\em link weakenings of $\cal O$};
\item[(iv)] ${\cal F}({\cal O})$ be the set of forests connecting (some, not
necessarily all) vertices of $\cal O$. A typical element of ${\cal
F}({\cal O})$ is
denoted by $\F$, and its set of links by $L(\F)\subset L({\cal O})$.
\end{itemize}

Assume  finally that link weakenings have been made to act in some smooth
way on $Z$ (a functional
depending on some extra parameters), thus defining a $C^{\infty}$
{\em link-weakened functional}
 $Z:L({\cal O})\to\R,
\vec{z}=(z_{\ell})_{\ell\in L({\cal O})}\mapsto Z(\vec{z})$ on the set of
pairs of objects, still denoted by $Z$ by a slight abuse of notation, such
that $Z$ (the
original functional) is equal to $Z(1,\ldots,1)=Z({\bf 1}_{L({\cal O})})$.

\end{Definition}

\begin{Definition}[one step of BK expansion] \label{def:one-step-BK}

The Brydges-Kennedy decoupling expansion consists in the following steps:

\begin{itemize}
\item[(i)] Taylor-expand $Z(1,\ldots,1)$ with respect to the parameters $(z_{\ell})$
simultaneously, namely,
\BEA  && Z(1,\ldots,1) = Z\left( (z_{\ell})_{\ell\in L({\cal O})}=
{\bf 0}\right)\\
&& +  \sum_{\ell_{1}\in L({\cal O})} \int_0^1 dw_1 \partial_{z_{\ell_{1}}} Z\left(
(z_{\ell})_{\ell\in L({\cal O})} =w_1 \right). \nonumber \EEA

\item[(ii)] Choose some link $\ell_{1}$ in the above sum.  Draw a link of strength $w_1$ between $o_{\ell_1}$ and $o'_{\ell_1}$ and
consider the new set ${\cal O}_1$ made up of the simple objects $o\in {\cal O}\setminus\{o_{\ell_1},o'_{\ell_1}\}$ and of the composite object
$\{(o_{\ell_1},o'_{\ell_1})\}$, with set of links $L({\cal O}_1)=L({\cal O})\setminus \{\ell_1\}$.
 Then   $\partial_{z_{\ell_{1}}} Z\left(z_{\ell_{1}}=w_1 , (z_{\ell})_{\ell\in  L({\cal O}_{1}) }\right)$
 must now be considered as  a functional
  of $(z_{\ell})_{\ell\in L({\cal O}_{1})} $.

\end{itemize}
\end{Definition}

 When iterating this procedure, composite objects grow up to be trees. This leads to the following result:

\begin{Proposition}[Brydges-Kennedy  or BK1 cluster formula] \label{prop:BK1}

$Z(1,\ldots,1)$ may be computed as an integral over weakening
parameters $w_{\ell}\in[0,1]$,
where $\ell$ does not range over the links of the total graph, but, more
restrictively,
over the  links of a forest on $\cal O$:

\BEQ Z(1,\ldots,1)=\sum_{\F\in {\cal F}({\cal O})} \left( \prod_{\ell \in
L(\F)}
\int_0^1 dw_{\ell} \right) \left( \left(\prod_{\ell\in L(\F)}
\frac{\partial}{\partial
z_{\ell} } \right) Z\right)(\vec{z}(\vec{w})),\EEQ
where $z_{\ell}(\vec{w}), \ell\in L({\cal O})$ is the infimum of the
$w_{\ell'}$ for
$\ell'$ running over the unique path from $o_{\ell}$ to $o'_{\ell}$ if
$o_{\ell}$ and
$o'_{\ell}$ are connected by $\F$, and $z_{\ell}(\vec{w})=0$ otherwise.

\end{Proposition}

If  it is not desirable to make a cluster expansion with respect to the links
between certain objects (of type $2$ in Proposition \ref{prop:BK2} below), then it
is sufficient to consider all these objects as belonging to the same composite object. This yields:

 \begin{Proposition}[restricted 2-type cluster or BK2 formula] \label{prop:BK2}

Assume ${\cal O}={\cal O}_1\amalg {\cal O}_2$. Choose as initial object an object
$o_1\in {\cal O}_1$ of type 1, and stop the Brydges-Kennedy expansion as soon as a
link to an object of type 2 has appeared. Then choose a new object of type $1$, and so on. This  leads to a  restricted expansion, for which {\em only} the  link variables $z_{\ell}$, with $\ell\not\in {\cal O}_2\times {\cal O}_2$, have been weakened. The following closed formula holds. Let ${\cal F}_{res}({\cal O})$ be the set of  forests $\F$ on $\cal O$, each component of which is (i) either a tree of
objects of type $1$, called {\em unrooted tree}; (ii)or a {\em rooted tree} such that only the root is of type $2$. Then
\BEQ Z(1,\ldots,1)=\sum_{\F\in {\cal F}{res}({\cal O})}
\left( \prod_{\ell \in L(\F)}
\int_0^1 dw_{\ell} \right) \left( \left(\prod_{\ell\in L(\F)}
\frac{\partial}{\partial
S_{\ell} } \right) Z(S_{\ell}(\vec{w})) \right),
 \EEQ

where $S_{\ell}(\vec{w})$ is either $0$ or the minimum of the $w$-variables running along the
unique path in $\bar{\F}$ from $o_{\ell}$ to $o'_{\ell}$, and $\bar{\F}$ is the
forest obtained from $\F$ by merging all roots of $\F$ into a single vertex.

\end{Proposition}

This restricted cluster expansion will  be useful for the Mayer expansion (see
section 3.4).


\subsection{Single scale cluster expansion}

\begin{Definition} \label{def:cluster}

\begin{itemize}

\item[(i)]

A {\em horizontal cluster forest} of level $\rho'\le \rho$, associated to a $d$-dimensional
vector Gaussian field $\psi=(\psi_1(x),\ldots,\psi_d(x))$, is a finite number of $M$-adic
intervals in $\D^{\rho'}$, seen as {\em vertices}, connected by {\em links}, without loops.  Any link
 $\ell$ connects $\Del_{\ell}$ to $\Del'_{\ell}$
$(\Del_{\ell},\Del'_{\ell}\in \D^{\rho'})$,
bears a double index, $(i_{\ell},i'_{\ell})\in\{1,\ldots,d\}\times\{1,\ldots,d\}$, and
may be represented as two half-propagators or simply half-segments put end to end, one
starting from $\Del_{\ell}$ with index $i_{\ell}$, and the other starting from $\Del'_{\ell}$ with  index $i'_{\ell}$.

 The
set of all horizontal cluster forests of level $\rho'$ is denoted by ${\cal F}^{\rho'}$. If $\F^{\rho'}\in {\cal F}^{\rho'}$,
then $L(\F^{\rho'})$ is its set of links.

\item[(ii)] A {\em horizontal cluster tree} is a {\em connected} horizontal cluster forest. Any
horizontal cluster forest decomposes into a product of cluster trees which are its {\em connected components}.

\item[(iii)]  If there exists a link between $\Del$ and $\Del'$ $(\Del,\Del'\in D^{\rho'})$ then we shall write $\Del\sim_{\F^{\rho'}} \Del'$, or
simply (if no ambiguity may arise) $\Del\sim\Del'$.

\end{itemize}
\end{Definition}

The following result is an easy consequence of  Proposition \ref{prop:BK1}, see
\cite{AbdRiv}.

\begin{Proposition}[single-scale cluster expansion] \label{prop:ssce}

Let
\BEQ Z_V(\lambda):=\int e^{- \int_V {\cal L}(\psi)(x) dx} d\mu(\psi),\EEQ
where $d\mu$ is the Gaussian measure associated to a  Gaussian field with $d\ge 1$ components,
$\psi(x)=(\psi_1(x),\ldots,\psi_d(x))$, with covariance matrix $C=C(i,x;j,y)$, and ${\cal L}$ is a local
functional.

Choose $\rho'\le \rho$.

Let $d\mu_{\vec{s}}(\psi)$, $\vec{s}=(s_{\Del,\Del'})_{(\Del,\Del')\in \D^{\rho'}\times\D^{\rho'}}
\in [0,1]$ such that $s_{\Del,\Del'}=s_{\Del',\Del}$ and $s_{\Del,\Del}=1$,  be the Gaussian measure with covariance kernel (if positive definite)
\BEQ C_{\vec{s}}(i,x;i',x')=s_{\Del_x,\Del_{x'}} C(i,x;i',x')\EEQ
if $\Del_x\ni x$, resp. $\Del_{x'}\ni x'$ $(\Del_x,\Del_{x'}\in \D^{\rho'})$ are the intervals of
size $M^{-\rho'}$ containing $x$, resp. $y$.

 Then
\BEA &&  Z_V(\lambda):=\sum_{\F^{\rho'}\in{\cal F}^{\rho'}} \left[ \prod_{\ell\in L(\F^{\rho'})}
\int_0^1 dw_{\ell} \int_{\Del_{\ell}} dx_{\ell} \int_{\Del'_{\ell}} dx'_{\ell}
C_{\vec{s}(\vec{w})}(i_{\ell},x_{\ell};i'_{\ell},x'_{\ell})  \right]  \nonumber\\
&& \int d\mu_{\vec{s}(\vec{w})}(\psi)  {\mathrm{Hor}}^{\rho'}
(e^{-\int_V {\cal L}(\psi)(x)dx})  \label{eq:3.6}  \EEA
where:
 \BEQ  {\mathrm{Hor}}^{\rho'}={\mathrm{Hor}}^{\rho'}(\F^{\rho'}; (x_{\ell})_{\ell\in L(\F^{\rho'})}, (x'_{\ell})_{\ell\in L(\F^{\rho'})}):= \prod_{\ell\in L(\F^{\rho'})}
\left(\frac{\del}{\del \psi^{\rho'}_{i_{\ell}}(x_{\ell})}  \frac{\del}{\del \psi^{\rho'}_{i'_{\ell}}(x'_{\ell})} \right); \EEQ
    $\vec{s}(\vec{w})=(\vec{s}_{\Del,\Del'}(\vec{w}))_{\Del,\Del'\in \D^{\rho'}}$, $s_{\Del,\Del'}(\vec{w})$,
$\Del\not=\Del'$  being the infimum of
the $w_{\ell}$ for $\ell$ running over the unique path from $\Del$ to $\Del'$ in $\F^{\rho'}$
if $\Del\sim_{\F^{\rho'}}\Del'$, and  $s_{\Del,\Del'}(\vec{w})=0$ else. \end{Proposition}

Note that $C_{\vec{s}(\vec{w})}$ is positive-definite with this definition
\cite{AbdRiv}, as a convex sum of evidently positive-definite kernels.


\subsection{Multi-scale cluster expansion}


As explained in the Introduction, cluster expansion has two main objectives. The first one is to express the
partition function as a sum over quantities depending essentially on a {\em finite} number of degrees of
freedom. For a given scale $j$, the horizontal cluster expansion perfectly meets this aim. The second one is to get
rid of the invariance by translation, which necessarily produces divergent quantities.

Multi-scale cluster expansion fulfills this program by building inductively, starting from the highest scale
$\rho$ and going down to scale $0$, {\em connected multi-scale clusters} called {\em polymers}. These are
{\em finite, connected} subsets of $\D$, extended over several scales, and containing at least one fixed interval
$\Del^j$ at the bottom which breaks invariance by translation. Like the horizontal cluster expansion, they are obtained
by a Taylor expansion with respect to some {\em $t$-parameters} $(t_{\Del})_{\Del}$ depending on the interval.

\begin{figure}[h]
  \centering
   \includegraphics[scale=0.5]{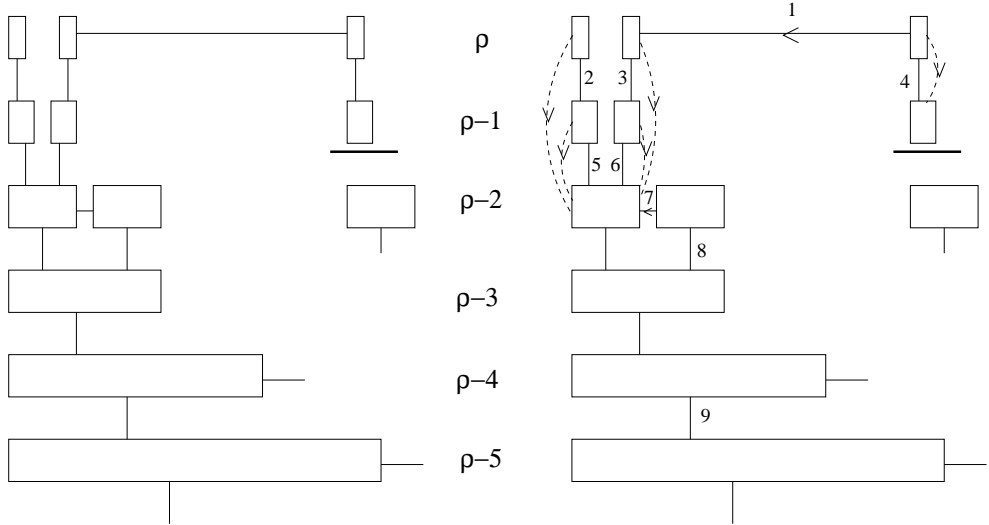}
   \caption{Polymer. A horizontal bold line below an interval means that it is {\em not} connected to the
intervals below it.}
  \label{Fig-polymer}
\end{figure}

Let us give a precise definition of
such an object (see left part of Fig. \ref{Fig-polymer}).

\begin{Definition} \label{def:polymer}
\begin{itemize}
\item[(i)] A {\em polymer of scale $j$} (typically denoted by $\P^j$) is a
tree
of $M$-adic intervals of scale $j$, i.e. a connected element of ${\cal
F}^j$. The set
of all polymers of scale $j$ is denoted by ${\cal P}^j$.
\item[(ii)] A {\em polymer down to scale $j$} (typically denoted by
$\P^{j\to}$) is
a connected graph connecting $M$-adic intervals of scale
$k=j,j+1,\ldots,\rho$, whose
links are of two types:

-- {\em horizontal links} connecting $M$-adic intervals of the same scale;
the restriction
of $\P^{j\to}$ to any given scale, $\P^{j\to}\cap \D^k$, $k\ge j$, is
required to
be a disjoint union of polymers of scale $k$, or simply in other words, an
element
of ${\cal F}^k$; furthermore, $\P^{j\to}\cap \D^j$ is assumed to be
non-empty;

-- {\em vertical links} or more explicitly {\em inclusion links} connecting an $M$-adic interval $\Del\in
\D^k$, $k>j$ to the unique interval $\Del'\in \D^{k-1}$ {\em below} $\Del$, i.e.
such that $\Del\subset\Del'$. These may be multiple links, with {\em degree of multiplicity $\tau_{\Del}$}.

A polymer $\P^{j\to}$ is also allowed to have an {\em external structure}, characterized by a subset ${\bf \Del}^j\subset \P^{j\to}\cap\D^j$ of intervals of scale $j$, such that each $\Del^j\in{\bf\Del}^j$ is connected below by an {\em external} inclusion link to $\D^{j-1}$.

\item[(iii)] The {\em horizontal skeleton} of a polynomial $\P^{j\to}$ is the disjoint union of the single-scale cluster
forests $\P^{j\to}\cap\D^k$, $k\ge j$ (all vertical links removed).

\end{itemize}

One denotes by ${\cal P}^{j\to}$ the set of polymers down to scale $j$, and
${\cal P}^{j\to}_{N_{ext}}\subset {\cal P}^{j\to}$ the subset of polymers with $N_{ext}:=\sum_{\Del^j\in
\P\cap\Del^j} \tau_{\Del}$  external links.

In particular, polymers in ${\cal P}^{j\to}_0$, without external links, are called {\em vacuum polymers}.

\end{Definition}

The easiest example is that of so-called {\em full-inclusion polymers}
 $\P^{j\to}\in {\cal P}^{j\to}$ containing {\em
all}
inclusion links. Such polymers may be obtained from arbitrary multiscale horizontal
cluster forests $\F^{j\to}=(\F^j,\F^{j+1},\ldots,\F^{\rho})$ by linking all pairs
$(\Del,\Del'), \Del,\Del'\in\F^{j\to}$ such that $\Del'$ is the unique interval lying
below $\Del$.

In general, if $\Del^k\subset\Del^{k-1}$, $\Del^k\in {\cal P}^{j\to} $, then {\em there is an
 inclusion link from $\Del^k$ to $\Del^{k-1}$ inside ${\cal P}^{j\to}$ (by definition) if and only if $\tau_{\Del^k}\not=0$. }

The integer $\tau_{\Del}$ corresponds in the momentum-decoupling expansion to
{\em the number of derivatives}  $\partial/\partial t_{\Del}$.
\bigskip

Whereas horizontal links are built in by independent horizontal cluster expansions at each scale,
the construction of vertical links depends on a procedure called {\em momentum decoupling} (or sometimes
{\em vertical cluster}) expansion, which we now set about to describe.

\medskip

Let us first give some definitions.

Assume positive numbers $t^{j}_x$, $0\le j\le \rho$, depending on $x$ but locally constant in
each $M$-adic interval $\Del^j\in\D^j$, have been defined.
Let us introduce  operators $T^{\to k}_x$ $(0\le k\le \rho,x\in\R)$ acting on the {\em low-momentum components}
of the fields at the point $x$, i.e. on $\phi^j(x)$ and $\sigma^j(x)$, $j<k$, by
\BEQ (T^{\to k}\psi)^j(x)=T^{\to k}\psi^j(x)=t^k_x t^{k-1}_x\ldots t^{j+1}_x \psi^j(x),\quad
k>j,\EEQ
\BEQ (T^{\to j}\psi)^j(x)=T^{\to j}\psi^j(x)=\psi^j(x),\EEQ
where $\psi=\phi$, $\partial\phi$ or $\sigma$.

The shorthand $(T^{\to k}\psi)^j$ emphasizes the idea that $T^{\to k}$ is not simply a
multiplication by some product of $t$-variables, but an operator acting diagonally
on the whole field $\psi$,
or rather on $\psi^{\to k}$ (seen as a vector).

In other words, writing $t_x^j=t_{\Del^j}$ if $x\in \Del^j$, so that $t$ is a real-valued function
on $\D^{\to\rho}=\uplus_{0\le j\le \rho}\D^j$, one has for $k\ge j$:
$ (T^{\to k}\psi)^j(x)=\prod_{k'=j+1}^k t_{\Del^{k'}_x}   \ \cdot \ \psi^j(x).$

We shall also use the following notation:
\BEQ (T\psi)^{\to k}:=\sum_{j\le k} (T^{\to k}\psi)^j. \label{eq:Tpsik} \EEQ

This weakened field, depending on the reference scale $k$,  will be called the {\em dressed low-momentum field at scale $k$}. Separating the $\psi^k$-component and the $t^k$-variables from the others yields equivalently:
\BEQ (T\psi)^{\to k}(x)=\psi^k(x)+t^k_x \sum_{j<k} (T^{\to(k-1)}\psi)^j(x). \label{eq:Tpsik-bis} \EEQ

Dressed interactions may be built quite generally out of dressed low-momentum fields. Let us
give a general definition.

\begin{Definition}[dressed interaction] \label{def:dressed-interaction}

\begin{itemize}

\item[(i)] Let ${\cal L}_q^{\to\rho}({\bf \psi})(x):=\lambda^{\kappa_q} \psi_{I_q}(x)=\lambda^{\kappa_q}
 \prod_{i\in I_q} \psi_i(x)$ be some arbitrary local functional
built out of a product of fields. Then the {\em momentum decoupling of ${\cal L}_q$},  or simply {\em dressed interaction}, is the following quantity,
\BEQ
{\cal L}_q^{\to\rho}(. ;\vec{t})(x):= \lambda^{\kappa_q} \left\{\prod_{i\in I_q}   (T\psi_i)^{\to\rho}(x)\ +\  \sum_{\rho'\le \rho}
(1-(t_x^{\rho'})^{|I_q|}) \prod_{i\in I_q}    (T\psi_i)^{\to(\rho'-1)}(x).  \right\} \label{eq:dressed-interaction}
\EEQ

\item[(ii)] (generalization to scale-dependent case)

  Let ${\cal L}_q^{\to\rho}(\psi)(x):=(\lambda^{\rho})^{\kappa_q}
\sum_{(j_i)_{i\in I_q}} K_q^{(j_i)} \prod_{i\in I_q}\psi_i^{j_i}(x)$ be a local functional with coupling constant $\lambda^{\rho}$ and coefficients $K_q^{(j_i)}$ depending on the scales. Then
one defines
\BEA && {\cal L}_q^{\to\rho}(\ .\ ;\vec{t})(x):=(\lambda^{\rho})^{\kappa_q}
\sum_{0\le (j_i)_{i\in I_q}\le\rho}   K_q^{(j_i)}\prod_{i\in I_q} (T^{\to\rho}\psi_i)^{j_i}(x)  \nonumber\\
  &&+ \sum_{\rho'\le\rho} (\lambda^{\rho'-1})^{\kappa_q} (1-(t_x^{\rho'})^{|I_q|})
   \sum_{0\le (j_i)_{i\in I_q}\le\rho'-1}   K_q^{(j_i)}\prod_{i\in I_q}
   (T^{\to (\rho'-1) }\psi_i)^{j_i}(x).  \nonumber\\
\EEA
where $\lambda^{\rho'}$, $\rho'\le\rho$ are the {\em renormalized} coupling constants.
\end{itemize}

Summing up in general the contribution of the different interaction terms, ${\cal L}_{int}:=\sum_{q=1}^p K_q {\cal L}_q$,  this defines a {\em dressed interaction}
${\cal  L}_{int}^{\to\rho}(\ .\ ;\vec{t})$ and a {\em dressed partition function}
\BEQ Z_V^{\to\rho}(\lambda;\vec{t})=\int e^{-\int_V {\cal L}^{\to\rho}_{int}(\ .\ ;\vec{t})(x)dx} d\mu(\psi).\EEQ

\end{Definition}

\bigskip

For the moment, vertical links have not been constructed, and nothing prevents a priori the different
horizontal clusters to move freely in space one with respect to the other, nor on the contrary to generate
via Wick pairings infinite multi-scale clusters. As we shall now see, the momentum-decoupling expansion
provides the mechanism responsible both for the translation-invariance breaking and the separation
of polymers.

\medskip Since this is an inductive procedure, let us first consider the result of horizontal cluster
expansion at highest scale $\rho$. It may be expressed as a {\em sum of monomials split over each connected
component} $\T_1^{\rho},\ldots,\T_{c^{\rho}}^{\rho}$ of $\F^{\rho}$, {\em with generic term (called $G$-monomial or simply monomial)}
\BEA && G_{c}^{\rho}:=\left[ \prod_{\ell\in L(\T_c^{\rho})} \left( \psi^{\rho}_{I_{\ell}}(x_{\ell})
\psi^{\rho}_{I'_{\ell}}(x'_{\ell}) \right) \right] \ \cdot \nonumber\\
&& \qquad \qquad \cdot \
\left[ \prod_{\ell\in L(\T_c^{\rho})} \left( t_{x_{\ell}}^{\rho} (T\psi_{J_{\ell}})^{\to(\rho-1)}(x_{\ell})\right)
\ \cdot\ \left( t_{x'_{\ell}}^{\rho} (T\psi_{J'_{\ell}})^{\to(\rho-1)}(x'_{\ell})\right) \right]  \nonumber\\
\label{eq:Gcrho} \EEA
for some (possibly empty) index subsets $I_{\ell},I'_{\ell},J_{\ell},J'_{\ell}\subset\{1,\ldots,d\}$, {\em multiplied by a product of propagators} as in eq. (\ref{eq:3.6}); the
$t^{\rho}$-variables dressing low-momentum components have been written explicitly as in eq.
(\ref{eq:Tpsik-bis}).

Let us draw an {\em oriented, downward dashed line}  from $\Del^{\rho}\in\T_c^{\rho}$ to some $\Del^{\rho'}\in(\Del^{\rho})^{\downdownarrows}$, $\rho'<\rho$ below $\Del^{\rho}$ if $G_{c}^{\rho}$ contains some {\em low-momentum}  field $\psi_i^{\rho'}(x_{\ell})$
with $x_{\ell}\in\Del^{\rho}$. Either $\Del^{\rho'}\not\in\F^{\rho'}$, or $\Del^{\rho'}$ belongs to
some connected component $\T^{\rho'}_{c'}$ of $\F^{\rho'}$. In the latter case, one has attached $\T^{\rho}_c$
to some cluster tree below, $\T^{\rho'}_{c'}$, by the {\em inclusion constraint} $\Del^{\rho}\subset\Del^{\rho'}$,
which prevents $\T_c^{\rho}$ from moving freely in space with respect to $\T^{\rho'}_{c'}$. It may also happen
that $G_{c}^{\rho}$ contains no low momentum field component $\psi^{\rho'}_i$, $\rho'<\rho$, so that $\T_c^{\rho}$
looks isolated; unfortunately, nothing prevents some horizontal cluster expansion at a lower scale $\rho'$
from generating some {\em high-momentum} field component $\psi_i^{\rho}(x)$, $x\in\Del^{\rho'}\in(\Del^{\rho})^{\downdownarrows}$, which may be represented as a reversed {\em upward} dashed line from $\Del^{\rho'}$ to $\Del^{\rho}$.

\medskip

In order to have an effective mechanism of separation of scales, we shall make a Taylor expansion to order 1 with respect
to the $(t_{\Del}^{\rho})_{\Del\in\T_c^{\rho}}$-variables of the product (G-monomial)$\times$
($(t^{\rho}_{\Del})_{\Del\in\T_c^{\rho}}$-dependent part of the dressed interaction),
$\tilde{G}_{c}^{\rho}:=G_c^{\rho}\ e^{-\int_{\T_c^{\rho}}{\cal L}_{int}(x)dx}$, namely (splitting $\tilde{G}_c^{\rho}$
into a product $\prod_{\Del} \tilde{G}^{\rho}_{\Del}=\prod_{\Del} G^{\rho}_{\Del}
e^{-\int_{\Del} {\cal L}_{int}(x) dx}$ over the fields located in each interval $\Del$ of scale $\rho$ in $\T^{\rho}_c$)

\BEQ \tilde{G}_{\Del}^{\rho}(t_{\Del}^{\rho}=1)=
\tilde{G}_{\Del}^{\rho}((t_{\Del}^{\rho})_{\Del\in\T_c^{\rho}}=0)+
\int_0^1 dt_{\Del}^{\rho} \partial_{t_{\Del}^{\rho}} \tilde{G}_{\Del}^{\rho}(t_{\Del}^{\rho}),\EEQ
thus producing a new set of monomials multiplied by the interaction.

Setting all $(t_{\Del}^{\rho})_{\Del\in\T_c^{\rho}}$  to zero has the effect, see eq.
(\ref{eq:dressed-interaction}) and (\ref{eq:Gcrho}),  of killing in the
interaction -- and hence in $G_c^{\rho}$ which is a derivative of the interaction -- {\em all mixed
terms} containing both $\psi^{\rho}_i(x)$ and $\psi^{\rho'}_{i'}(x')$, with $x\in\Del^{\rho}\subset\D^{\rho}$,
$x'\in\Del^{\rho'}\subset\D^{\rho'}$, $\Del^{\rho'}\supset\Del^{\rho}$ as above. Hence, in all the terms
$\tilde{G}_{c}^{\rho}((t_{\Del}^{\rho})_{\Del\in\T_c^{\rho}}=0)$, $\T_c^{\rho}$ has been effectively isolated
from all other intervals in $\D$; it constitutes a (single-scale) isolated polymer. Thus, by letting $t_{\Del}=0$,
 one {\em cuts} all  dashed lines crossing from $\Del$ (or above $\Del$ in general) to
$\Del^{\downdownarrows}$, and sets up a wall between $\Del$ and $\Del^{\downdownarrows}$.

On the contrary, derivating with respect to $t_{\Del}^{\rho}$, $\Del\in\T_c^{\rho}$, produces necessarily
low-momentum components, and hence vertical links materialized by dashed lines as above. In this case,
one links $\Del^{\rho}$ to $\Del^{\rho-1}$ by a {\em full line}, signifying that all corresponding monomials
will contain some low-momentum field; this implies in turn the existence of a {\em downward dashed line} or {\em wire}
connecting some $\Del^{\rho}\in\T_c^{\rho}$ to some $\Del^{\rho'}\in(\Del^{\rho})^{\downdownarrows}$ as above
(note however that there isn't necessarily a  dashed line from $\Del^{\rho}$ to $\Del^{\rho-1}$).

Due to the necessity of {\em renormalization}
 (see section 4) and to the {\em domination problem} (see \S 6.2), one may need to Taylor expand
to higher order, up to order $N_{ext,max}+O(n(\Del))$. One obtains thus
a polymer with a certain number of external legs per interval. Choosing the Taylor integral rest term
 for some $\Del$ leads to a polymer with $\ge N_{ext,max}$ external legs, which does not need to be renormalized.
On the other hand, a polymer whose {\em total} number of external legs is $<N_{ext,max}$ requires renormalization (see section 4).

Let us emphasize at this point that high-momentum fields have {\em two scales} attached to them: $j$ and
$k$ for a field $\psi_i^k(x)$, $x\in\Del^j$ ($k>j$), produced by the horizontal/vertical cluster
expansion at scale $j$; but low-momentum fields $\psi_i^h$ have {\em three scales}. There are in fact
two cases:

\begin{itemize}

\item[(i)] If $\beta_i<-D/2$ then $\psi_i$ is not separated into a sum (low-momentum field average)+(secondary
    field).  The genesis of $\psi_i^h$ (contrary to the high-momentum case) is actually a process which
    may not be understood apart from the multi-scale cluster expansion. At their {\em production scale} $k$, low-momentum fields are of the form $(T\psi_i)^{\to(k-1)}(x)=\sum_{j=0}^{k-1}
    \left[\prod_{k'=j+1}^{k-1} t_{\Del_x^{k'}}\right] \psi^j(x)$. Successive $t$-derivations of
    scale $k-1$, $k-2,\ldots$ "push" $(T\psi_i)^{\to(k-1)}(x)$ downward like a down-going elevator,
    in the sense that $\partial_{t_{\Del_x^{k-1}}} (T\psi_i)^{\to(k-1)}(x)$ is by construction
    of scale $\le k-2$, $\partial_{t_{\Del_x^{k-2}}} (T\psi_i)^{\to(k-1)}(x)$
    of scale $\le k-3$ and so on. But of course, $t$-derivations may act on other fields instead.
    The last $t$-derivation acting on $(T\psi_i)^{\to(k-1)}(x)$ {\em drops} $(T\psi_i)^{\to(k-1)}(x)$
    at  a certain scale $j<k$. Then $(T\psi_i)^{\to(j-1)}(x)$ leaves the elevator and is torn apart
    into its scale components $\left( (T^{\to(j-1)}\psi_i)^h\right)_{h<j}$ through {\em free falling}. Thus $\psi_i^h$ has a {\em production scale} $k$ and a {\em dropping scale} $j$, while $h$ itself may be called its {\em free falling scale} or simply its scale.

    \item[(ii)] If $\beta_i\ge -D/2$, then, at the {\em dropping scale}, $(T\psi_i)^{\to(j-1)}(x)$ is
    separated from its average $(T\psi_i)^{\to(j-1)}(\Del_x^j)$ which must be dominated apart,
    while $(T\psi_i)^{\to(j-1)}(x)-(T\psi_i)^{\to(j-1)}(\Del_x^j)$ splits into its scale components
    $\left( (T^{\to(j-1)}\del^j\psi_i)^h\right)_{h<j}$ through {\em free falling} as in case (i).

    \end{itemize}

\bigskip

 The extension of the above procedure to lower scales is straightforward and leads
to the following result.

\begin{Definition}[multi-scale cluster expansion] \label{def:3.10}

\begin{enumerate}

\item
Fix a multi-scale horizontal cluster expansion $\F^{j\to}=(\F^j,\ldots,\F^{\rho})$ and consider a polymer
down to scale $j$, $\P^{j\to}$, with horizontal skeleton $\F^{j\to}$. To such a polymer is associated
a sum of products ($G$-monomial)$\times$(dressed interaction ${\cal L}_{int}(\ .\ ;\vec{t})$), where
all $t_{\Del}$-variables such that $\Del\in\F^{j\to}$ and $\tau_{\Del}< N_{ext,max} +O(n(\Del))$ have been set to $0$,
and $G$ is one of the monomials obtained by expanding $\left[\prod_{k\ge j} {\mathrm{Vert}}^k . {\mathrm{Hor}}^k \right]
 \cdot\
 e^{-\int_V {\cal L}_{int}(\ .\ ;\vec{t})(x) dx}$, where
 ${\mathrm{Hor}}^k$ is as in Proposition \ref{prop:ssce}, and ${\mathrm{ Vert}}^k$ is the following operator, with $N'_{ext,max}=N_{ext,max}+O(n(\Del))$,
\BEQ {\mathrm{ Vert}}^k=  \prod_{\Del\in\F^k} \left( \sum_{\tau_{\Del}=0}^{N'_{ext,max}-1} \partial^{\tau_{\Del}}_{t_{\Del}}\big|_{t_{\Del}=0} +  \int_0^1  dt_{\Del}
\frac{(1-t_{\Del})^{N'_{ext,max}-1}}{
(N'_{ext,max}-1)!} \partial_{t_{\Del}}^{N'_{ext,max}} \right).\EEQ

\item (definition of $n(\Del)$)
Fix $\F^{j\to}$ and let $\Del\in \D^{j\to}$. Then $n(\Del)$ is the number of intervals of scale $j(\Del)$ connected to $\Del$ by the forest $\F^{j(\Del)}$.

\item (definition of $N(\Del)$)
Fix $\F^{j\to}$ and some $G$-monomial, and let $\Del\in \D^{j\to}$. Then
$N(\Del)$ is the number of fields $\psi^{j(\Del)}(x)$, $x\in\Del$ of scale $j(\Del)$ lying in the interval $\Del$.

\end{enumerate}
\end{Definition}

Classically, one sets apart polymers made up of one
 isolated interval $\Del^j$ where {\em no vertex} has been produced; this means that $t_{\Del^{j+1}}=0$ for all intervals $\Del^{j+1}\subset\Del^j$; $t_{\Del^j}=0$; and
 $s_{\Del^j,\Del'}=0$ if $\Del'\in\D^j\setminus\{\Del^j\}$.   Write as usual
 ${\cal L}_{int}:=\sum_{q=1}^p K_q\lambda^{\kappa_q} {\cal L}_q$. Their contribution to the partition functions reads simply
\BEA && \int d\mu(\psi) e^{-\int_{\Del} {\cal L}_{int}^{\to\rho}(\ .\ ;\vec{t}=0)(x) dx} \nonumber\\
&& \qquad \qquad =1-\sum_{q=1}^p \kappa_q K_q\lambda^{\kappa_q}\int d\mu(\psi)\ \cdot\ \int_0^1 dv
\int {\cal L}_q^{\to\rho}(\ .\ ;\vec{t}=0)(x)dx \ \cdot\
e^{-v \int_{\Del} {\cal L}_{int}^{\to\rho}(\ .\ ;\vec{t}=0)(x) dx } \nonumber\\
&& \qquad \qquad =: 1-A.
\label{eq:isolated-interval}
\EEA
by using a Taylor expansion to order $1$. In the above expression, all fields have same scale $j$ since all $t$-coefficients connecting $\Del^j$ have been set to zero.  In all cases where the
exponentiated interaction is bounded by $1$ -- which is the case in our $(\phi,\partial\phi,\sigma)$-model -- the Cauchy-Schwarz inequality implies
that $|A|^2\le  K\lambda^{2\kappa}$, where $\kappa=\min\{\kappa_q; q=1,\ldots,p\}$ since the  integrated vertex is homogeneous of degree $0$.  Now, one eliminates the factor $1$ by releasing
the constraint that the disjoint union of all polymers must span the whole volume $V$, and ends
up with a polymer whose evaluation is of order $O(\lambda^{\kappa})$. This makes it possible to
treat such polymers on an equal footing with the other polymers where some vertex has been produced (see introduction to section 6.1 for the general principles of the bounds for cluster
expansions).

{\bf Remark.}
Note that $N(\Del)$ is at most of order $O(n(\Del))$ for a {\em single-scale} cluster expansion. This is not true for a {\em multi-scale} cluster expansion.  Namely, fix $\Del^h\in\D^h$. Assume a low-momentum
field $\psi_i^h(x)$ ($\beta_i< -D/2$) or $\del^j\psi_i^h(x)$ ($\beta_i\ge -D/2$) has been produced at scale $k$ and dropped 
inside $\Del^j$ at scale $j$, with $k>j>h$. Although the number of fields produced at scale $k$ in an interval $\Del^k\subset\Del^j$
 increases exponentially 
with $k-j$, there are $\le N_{ext,max}+O(n(\Del^j))$ low-momentum fields dropped inside $\Del^j$  for a {\em given} $G$-monomial, since $\partial_{t_{\Del^j}}$ occurs
at a power $\le N_{ext,max}+O(n(\Del^j))$. On the other hand,  the number of low-momentum fields $\psi^h(x)$, $x\in\Del^j$
 with $\Del^j\subset\Del^h$ originating from  vertices of scale $j$ may be of order $\# \{\Del^j\in\D^j; \Del^j\subset\Del^h\}=M^{D(j-h)}$.
 This is a well-known phenomenon, called {\em accumulation of low-momentum fields}. This "negative" spring-factor must be combined with
 the rescaling spring factor $M^{2\beta(j-h)}$, see Corollary \ref{cor:spring-factor}, resulting in $M^{(D+2\beta)(j-h)}$, a {\em positive}
 spring-factor if $\beta<-D/2$. This accounts for the (already mentioned) fact that secondary fields need not be produced for fields with 
scaling dimension $<-D/2$, see subsection 6.1 for details.

At a given scale $j$, composing the horizontal cluster and momentum-decoupling expansions  at all scales $\ge j$
yields the following result, easy to show by induction:

\begin{Lemma}[result of the expansion above scale $j$]  \label{lem:3.10}

\begin{enumerate}

\item

\BEQ  Z_V^{\rho}(\lambda;\vec{t}^{\to (j-1)})=  \int
 d\mu(\psi^{\to(j-1)}) Z_V^{j\to\rho} (\lambda; (\psi^h)_{h\le j-1}; \vec{t}^{\to(j-1)}),\EEQ with
{\tiny \BEA &&  Z_V^{j\to\rho} (\ .\ )=\left( \sum_{\F^{j}\in{\cal F}^{j}}
 \left[ \prod_{\ell^{j}\in L(\F^{j})}
\int_0^1 dw_l^{j} \int_{\Del_{\ell^{j}}} dx_{\ell^{j}} \int_{\Del'_{\ell^{j}}}
 dx'_{\ell^{j}}
C_{\vec{s}^{j}(\vec{w^{j}})}(x_{\ell^{j}},x'_{\ell^{j}})  \right] \right)
 \nonumber\\
&& \ldots \left( \sum_{\F^{\rho}\in{\cal F}^{\rho}} \left[ \prod_{\ell^{\rho}\in L(\F^{\rho})}
\int_0^1 dw_l^{\rho} \int_{\Del_{\ell^{\rho}}} dx_{\ell^{\rho}} \int_{\Del'_{\ell^{\rho}}} dx'_{\ell^{\rho}}
C_{\vec{s}^{\rho}(\vec{w^{\rho}})}(x_{\ell^{\rho}},x'_{\ell^{\rho}})  \right] \right) \nonumber\\
&&  \int d\mu^{j\to\rho}_{\vec{s}(\vec{w})}(\psi^{j\to\rho})
  \left[\prod_{k\ge j} {\mathrm{Vert}}^k\ {\mathrm{Hor}}^k \right] e^{-
\int_V {\cal L}(.;\vec{t})(x)dx}  \nonumber\\  \label{eq:ZVrho'}
\EEA }
where $ d\mu^{j\to\rho}_{\vec{s}(\vec{w})}(\psi^{j\to\rho})$
is a short-hand for
 $\prod_{k=j}^{\rho} d\mu_{\vec{s}^k(\vec{w}^k)}(\psi^k)
$.

\item The right-hand side (\ref{eq:ZVrho'}) depends on the low-momentum fields $(\psi_i^{h})_{h<j}$
only through the dressed fields $(T\psi_i)^{\to(j-1)}$, since the $t$-variables of scale $\le j-1$ have not been touched.
Hence it makes sense to consider the quantity $Z_V^{j\to \rho}(\lambda):=Z_V^{j\to \rho}(\lambda; (\psi^h)_{h\le j-1}=0; \vec{t}^{\to(j-1)}=
{\bf 1})$.

\item The partition function $Z_V^{j\to \rho}(\lambda)$ writes
\BEA  Z_V^{j\to \rho}(\lambda) &=& \sum_{N=1}^{\infty} \frac{1}{N!}
\sum_{{\mathrm{non-overlapping}}\
\P_1,\ldots,\P_N\in {\cal P}_0^{j\to}} \prod_{n=1}^N F_{HV}(\P_n) \nonumber\\
&=& \sum_{N=1}^{\infty} \frac{1}{N!} \sum_{
\P_1,\ldots,\P_N\in {\cal P}_0^{j\to} } \prod_{n=1}^N F_{HV}(\P_n) \ \cdot \
\prod_{\ell=\{\P,\P'\}} {\bf 1}_{\P,\P'\ {\mathrm{non-overlapping}}} \nonumber\\  \label{eq:3.19} \EEA
where $F_{HV}$, called {\em polymer functional associated to horizontal (H) and vertical (V) cluster expansion},
is the contribution of each polymer $\P^{j\to}$ to the right-hand side of (\ref{eq:ZVrho'}).

\end{enumerate}
\end{Lemma}

Later on, the polymer functional $F_{HV}$ will be replaced with the {\em renormalized (R) polymer functional} $F_{HVR}$,
or simply $F$, hence eq. (\ref{eq:3.19}) shall not be used in this form . The general discussion in the next subsection, which does not depend on the precise form of $F$, shows how to get rid
of the non-overlapping conditions.



\subsection{Mayer expansion}


Recall the polynomial evaluation function $F_{HV}$ introduced in the expansion of $Z_V^{j\to\rho}$, see eq. (\ref{eq:3.19}). It is unsatisfactory in this form because of the non-overlapping conditions which make it impossible to compute directly a finite quantity out of it, see Introduction. Were this the only source of trouble, it would suffice to make a global Mayer expansion for $Z_V^{0\to\rho}(\lambda)$ after the multi-scale cluster expansion has been completed. However, it is also unsatisfactory when renormalization is required; local parts of diverging graphs (for reasons accounted for in section 4) must be discarded and resummed into an exponential, thus leading to a counterterm in the interaction. This forces upon us a sequence of three moves at {\em each scale} $j$,
starting from highest scale $\rho$:  (1) a separation of the local part of diverging graphs; (2) the Mayer expansion proper, at scale $j$; (3) the construction of the interaction counterterm (also called renormalization phase).

\bigskip

Let us formalize this into the following:

\medskip

{\bf Induction hypothesis at scale $j$}. {\em After completing all expansions of scale $\ge j+1$ and the horizontal/vertical cluster expansion at scale $j$, $Z_V^{\to\rho}(\lambda;\vec{t})$ has been rewritten as \BEQ Z_V^{\to\rho}(\lambda;\vec{t}^{\to(j-1)})=\prod_{k=j+1}^{\rho} e^{|V| M^k f^{k\to\rho}(\lambda)} \ \cdot\ \int d\mu(\psi^{\to(j-1)}) Z_V^{j\to\rho}(\lambda;\vec{t}^{\to(j-1)}; \psi^{\to(j-1)}),\EEQ
where $f^{k\to\rho}(\lambda)$ may be reinterpreted as a scale $k$ free energy density per degree of freedom,
with
\BEA && Z_V^{j\to\rho}(\lambda;\vec{t}^{\to(j-1)};\psi^{\to(j-1)})= \nonumber\\
&& \qquad \sum_{N=1}^{\infty} \frac{1}{N!} \sum_{{\mathrm{non}} -j- {\mathrm{overlapping}} \ \P_1,\ldots,\P_N\in {\cal P}^{j\to}}
 \int d\vec{w}^{j\to} \int d\mu_{\vec{s}(\vec{w})}(\psi^{j\to})
 \prod_{n=1}^N F_{HV}^j(\P_n;\psi), \nonumber\\  \EEA
 see Lemma \ref{lem:3.10} for notations,
where "non$-j-$overlapping" means that $\P_1\cap\D^j,\ldots,\P_N\cap\D^j$ are
non-overlapping single-scale polymers, and $F^j_{HV}(\P_n;\psi)$ depends only
on the values of $\psi$ on the support of $\P_n$. }

For $j=\rho$, this formula is equivalent to the single-scale expansion of $Z_V^{\rho\to\rho}(\lambda)$ of  eq. (\ref{eq:3.19}).

According to the general expansion scheme (see introduction to section 3), we must now perform three tasks.

\begin{enumerate}

\item {\em Separation of local part of diverging graphs}

Consider the polymer evaluation function $F_{HV}^j(\P_n;\psi)$. If
$\P_n$ has $N_{ext}=\sum_{\Del^j\in \P_n\cap\D^j} \tau_{\Del^j}<N_{ext,max}$
external legs, situated in (possibly coinciding) intervals $\Del_i^j$, it may be {\em superficially divergent}, in which case one separates its local part according to the rule explained in section 4: letting
$(\psi_i^{\to(j-1)}(x_i))_{i\in I_q}$, with $x_i\in\Del_i^j$, be its external structure, $F_{HV}^j(\P_n;\psi)$ writes
\BEA && \prod_{i\in I_q} \int_{\Del_i^j} dx_i F^j_{HV,{\mathrm{amputated}}}(\P_n;\psi;(x_i)) \prod_{i\in I_q} \psi_i^{\to(j-1)}(x_i)
  \nonumber\\
&& =F^j_{HV,{\mathrm{local}}}(\P_n;\psi)+
\del F^j_{HV}(\del^{N_{ext,max}-N_{ext}} \P_n;\psi), \nonumber\\ \EEA
where
\BEQ F^j_{HV,{\mathrm{local}}}(\P_n;\psi)=\prod_{i\in I_q} \int_{\Del_i^j} dx_i  F^j_{HV,{\mathrm{amputated}}}(\P_n;\psi;(x_i)) \ \cdot\ \left(\frac{1}{|I_q|} \sum_{i\in I_q} \psi_{I_q}^{\to(j-1)}(x_i) \right)
 \EEQ
and $F_{HV}^j(\P_n;\psi)$ minus its local part is now thought as if it had $N_{ext,max}-N_{ext}$ supplementary external legs shared in an arbitrary way among the intervals $\Del_i^j$, which produces a polymer denoted by $\del^{N_{ext,max}-N_{ext}} \P_n$ belonging to ${\cal P}^{j\to}_{N_{ext,max}}$.

Equivalently (see subsection 4.1), considering only local parts, according
to our new convention,
\BEQ F^j_{HV}(\P_n;\psi)=\prod_{i\in I_q} \int_{\Del_i^j} dx_i  F^j_{HV,{\mathrm{amputated}}}(\P_n;\psi;(x_i)) \ \cdot\ \prod_{i\in I_q}
\psi_i^{\to(j-1)}(x_i),\EEQ
where now
\BEQ F^j_{HV,{\mathrm{amputated}}}(\P_n;\psi;(x_i))=\frac{1}{|I_q|}
\sum_{i\in I_q} \int \prod_{i'\in I_q,i'\not=i} dx_{i'}
F^j_{HV}(\P_n;\psi;(x_i)).\EEQ

\item {\em Mayer expansion at scale $j$}

We shall now   apply
the restricted cluster expansion, see Proposition \ref{prop:BK2}, to
the functional $Z_V^{j\to\rho}(\lambda; \vec{t}; \psi^{\to(j-1)})$.  The "objects" are
multi-scale polymers $\P$ in    ${\cal O}=\{\P_1,\ldots,\P_N\}$ (see induction hypothesis above); a link
$\ell\in L({\cal O})$ is a pair of polymers $\{\P_n,\P_{n'}\}$, $n\not=n'$.
Objects of type 2 are  polymers with $\ge N_{ext,max}$ external legs whose non-overlap conditions may not be removed  properly at this stage, because they are in any case  dependent of each other in an as yet unforeseeable way through the  low-momentum fields.
All other objects are of type 1, they belong to
 ${\cal P}^{j\to}_{<N_{ext,max}}:=\uplus_{N=0}^{N_{ext,max}-1} {\cal P}^{j\to}_N$; they
are either {\em vacuum polymers}, i.e. polymers with no external legs, or superficially divergent polymers whose contribution to the interaction counterterm one would like to compute.

 Link weakenings $\vec{S}=(S_{ \{\P,\P'\}})_{ \{\P,\P'\}\in L({\cal
O}) } \in [0,1]^{L({\cal O})} $ act
on $Z_V^{j\to\rho}(\lambda; \vec{t}^{\to(j-1)}; \psi^{\to(j-1)})$ by replacing the non-overlapping condition
 \BEA {\mathrm{NonOverlap}}(\P_1,\ldots,\P_N) &:=& \prod_{(\P_n,\P_{n'}) } {\bf 1}_{\P_n,\P_{n'}\ {\mathrm{non}}-j-{\mathrm{overlapping}}} \nonumber\\
 &=&
\prod_{(\P_n,\P_{n'}) } \prod_{\Del\in\P_n\cap\D^j,\Del'\in\P_{n'}\cap\D^j}
 \left( 1 + \left(
 {\bf 1}_{\Del\not=\Del'}-1 \right)\right) \nonumber\\ \EEA

with a weakened non-overlapping condition

\BEA && \prod_{(\P_n,\P_{n'}) } \prod_{\Del\in{\bf\Del}_{ext}(\P_n),\Del'\in
{\bf\Del}_{ext}(\P_{n'})} {\bf 1}_{\Del\not=\Del'} \ \cdot\nonumber\\
 &&\qquad \left( 1 +  S_{\{\P_n,\P_{n'}\}} \left(
\prod_{\Del\in\P_n\cap\D^j,\Del'\in\P_{n'}\cap\D^j,(\Del,\Del')\not\in
{\bf\Del}_{ext}(\P_n)\times
{\bf\Del}_{ext}(\P_{n'})}  {\bf 1}_{\Del\not=\Del'}-1 \right)\right),  \nonumber\\ \label{eq:3.23} \EEA

where ${\bf\Del}_{ext}(\P)\subset\P\cap\D^j$ is the subset of intervals $\Del^j$ with
external legs, i.e. such that $\tau_{\Del^j}\not=0$.

Note that each factor in (\ref{eq:3.23}) ranges in $[0,1]$.
We ask the reader to accept this definition as it is and wait till the remark
after Proposition \ref{prop:Mayer} for explanations.

Let us now give some necessary precisions. The Mayer expansion is really applied to the non-overlap function NonOverlap and {\em not} to
$Z_V^{j\to\rho}(\lambda;\vec{t}^{\to(j-1)};\psi^{\to(j-1)})$. Hence one must
still extend the function $\int d\vec{w}^{j\to} \int d\mu^{j\to}_{\vec{s}(\vec{w})}(\psi^{j\to}) \prod_{i=1}^N F^j_{HV}(\P_n;\psi)$  to the case when the $\P_n$, $n=1,\ldots,N$ have some overlap. The natural way to do this is to assume that
 the random variables $(\psi^j\big|_{\P_n})_{n=1,\ldots,N}$ remain independent even when they overlap. This may be understood in the following way. Choose a different color for each polymer $\P_n=\P_1,\ldots,\P_N$, and paint with that color
 {\em all} intervals $\Del^j\in \P_n\cap\D^j$. If $\Del^j\in{\bf \Del}_{ext}(\P_n)$, then its
 external links to the interval $\Del^{j-1}$ below it are left in black. The rules (\ref{eq:3.23}) imply that intervals with different colors may superpose; on the other hand, {\em external inclusion links} may {\em not}, so that: (i)
 {\em low-momentum fields} $\psi^{\to(j-1)}(x)$, $x\in\Del^j$ with $\Del^j\in{\bf\Del}_{ext}(\P_n)$,
do not superpose and may be left in black (till the next expansion stage at scale $j-1$ at least); (ii)
 the color of {\em high-momentum fields} of scale $j$ created at a later stage may be determined without ambiguity.

Hence one must see $\psi^j$ as living on a two-dimensional set, $\D^j\times\{ {\mathrm{colors}}\}$, so that copies of $\psi^j$ with different colors are independent of each other. This defines a new, {\em extended} Gaussian measure
$d\tilde{\mu}_{s^j(w^j)}(\tilde{\psi}^j)$ associated to an {\em extended} field $\tilde{\psi}^j:\R^D\times \{{\mathrm{colors}}\}\to\R$, and {\em Mayer-extended polymers}.
 By abuse of notation, we shall skip the tilde in the sequel, and always implicitly extend the fields and the measures of each scale.
 Mayer-extended polymers shall be considered as (colored)  polymers in section 6.

\bigskip

This gives the following expansion for $Z_V^{j\to\rho}(\lambda; \vec{t}^{\to(j-1)}; \psi^{\to(j-1)})$.

\begin{Proposition}[Mayer expansion] \label{prop:Mayer}

Let ${\cal F}({\cal P}^{j\to})$ be the set of all forests of polymers whose each component $\T$ is (i) either a tree of polymers of type 1 (called: {\em unrooted tree}); (ii) or a {\em rooted tree of polymers} such that {\em only} the root is of type 2. Then

\BEQ Z_V^{j\to\rho}(\lambda; \vec{t}^{\to(j-1)}; \psi^{\to(j-1)})
= \sum_{\F \in {\cal F}({\cal P}^{j\to})}
{\mathrm{Mayer}} (Z_V^{j\to\rho}(\lambda; \vec{t}^{\to(j-1)}; \psi^{\to(j-1)}); \F), \EEQ
with
\BEA &&  {\mathrm{Mayer}} (Z_V^{j\to\rho}(\lambda; \vec{t}^{\to(j-1)}; \psi^{\to(j-1)}); \F) \nonumber\\
&& \qquad =
\left( \prod_{\ell \in L(\F)}
\int_0^1 dW_{\ell} \right) \left( \left(\prod_{\ell\in L(\F)}
\frac{\partial}{\partial
S_{\ell} } \right) Z_V^{j\to\rho}(\lambda; \vec{t}^{\to(j-1)}; \psi^{\to(j-1)})\right)(\vec{S}(\vec{W})) \nonumber\\ \EEA

 where $S_{\ell}(\vec{W})$ is either $0$ or the minimum of the $W$-variables running along the
unique path in $\bar{\F}$ from $o_{\ell}$ to $o'_{\ell}$, and $\bar{\F}$ is the
forest obtained from $\F$ by merging all roots in ${\cal P}^{j\to}_{\ge N_{ext,max}}$ into a single vertex.

\end{Proposition}

As a result, (i) polymers $\P_{\ell},\P'_{\ell}$ linked by a Mayer link are $j$-overlapping (otherwise the derivative $\partial_{S_{\ell}}$ would produce a zero factor); (ii) pairs of {\em vacuum} polymers $\P,\P'$ belonging to different Mayer trees come with the factor 1: they have lost
their non-overlap conditions and may superpose each other freely (in other words, they have become transparent to each other).  Hence  $Z_V^{j\to\rho}(\lambda; \vec{t}^{\to(j-1)}; \psi)$ factorizes as a product
\BEA &&  Z_V^{j\to\rho}(\lambda; \vec{t}^{\to(j-1)}; \psi)=
e^{F^{j\to\rho}_{HVM}}\int d\mu(\psi^{\to(j-1)})  \ \cdot\  \nonumber\\
&&\qquad \cdot\ \sum_{N=1}^{\infty} \frac{1}{N!} \sum_{{\mathrm{non}}-j-{\mathrm{overlapping}}\ \P'_1,\ldots,\P'_n\in {\cal P}^{j\to}} \prod_{n=1}^N F^j_{HVM}(\P'_n;\psi) \nonumber\\
 \EEA
where $F^{j\to\rho}_{HVM}$ is the contribution of all (unrooted) trees of {\em vacuum polymers}. Denote by $f^{j\to\rho}(\lambda)$ the quantity obtained by fixing {\em one} interval $\Del^j$ of scale $j$ belonging to one of the polymers of the tree.  Summing over all $\Del^j$, one obtains an overall factor $e^{|V| M^j f^{j\to\rho}(\lambda)}$.

\item {\em Renormalization phase}

For simplicity we shall assume that only 2-point functions $\langle \psi_i\psi_{i'}\rangle$, $1\le i,i'\le d$ need to be renormalized (which is the case for the $(\phi,\partial\phi,\sigma)$-model). Denote by $-\half\int_{\Del^j} dx b_{i,i'}^j \psi_i^{\to(j-1)}(x)
\psi_{i'}^{\to(j-1)}(x)$ the contribution to $F^{j\to\rho}_{HVM}$ of all unrooted trees of polymers containing exactly {\em one} polymer with 2 external legs in $\Del^j$, plus a {\em cloud} of vacuum polymers, attached
to it, directly or indirectly, by Mayer links.  The intervals without external legs of the different unrooted trees of polymers have become transparent to each other, and to the rooted trees too, hence $b^j_{i,i'}$ may be computed by considering {\em one} unrooted tree of polymers with 2 external legs, irrespectively of the position of the other trees of polymers.   By translation invariance, $b_{i,i'}^j$ is a constant, which is fixed by induction on $j$ by demanding that
\BEQ 0\equiv \int dy \frac{\partial^2}{\partial \psi_i^{\to(j-1)}(x) \partial \psi_{i'}^{\to(j-1)}(y)} Z^{j\to\rho}(\lambda;\psi^{\to(j-1)})\big|_{\psi^{\to(j-1)}=0}, \EEQ
or better (so as to eliminate higher scale polymers which depend only on $b_{i,i}^{(j+1)\to}$)
\BEQ 0\equiv \int dy \frac{\partial^2}{\partial \psi_i^{\to(j-1)}(x) \partial \psi_{i'}^{\to(j-1)}(y)} \left(Z^{j\to\rho}-Z^{(j+1)\to\rho}\right)
(\lambda;\psi^{\to(j-1)})\big|_{\psi^{\to(j-1)}=0}. \label{eq:3.32} \EEQ

Hence eq. (\ref{eq:3.32}) yields $b_{i,i'}^j$ as a sum over polymers (containing at least one interval of scale $j$) of an  expression depending
itself on $b_{i,i'}^j$' -- an {\em implicit} equation.

Local parts of 2-point functions are generated by the
exponential over the whole volume, $e^{-\half\int_V dx b_{i,i'}^j \psi_i^{\to(j-1)}(x)
\psi_{i'}^{\to(j-1)}(x)}$, which may be seen as a counterterm in the interaction. This counterterm disappears
  by renormalizing \footnote{Note that $K$ -- contrary to the original bare kernel -- is only almost-diagonal,
 namely, $\int\psi_i^j(x)\psi_{i'}^{j'}(x)dx=0$ by momentum conservation if $|\supp(\chi^j)\cap\supp(\chi^{j'})|=0$, i.e.
 if $|j-j'|\ge 2$.} the Gaussian covariance kernel $C_{\psi}=K_0^{-1}$  -- previously renormalized down to scale $j+1$ --  to $K^{-1}$,  where
\BEQ K(\psi,\psi)=K_0(\psi,\psi)+\del K^j(\psi,\psi)=K_0(\psi,\psi)+
\int_V dx b_{i,i'}^j \psi_i^{\to(j-1)}(x)
\psi_{i'}^{\to(j-1)}(x). \EEQ
Note that the  renormalization of $C_{\sigma}$ in the case of the $(\phi,\partial\phi,\sigma)$ shall only be performed
at scale $\rho$, while the other counterterms shall be left out as a counterterm in the interaction (see section 5).

After applying a horizontal/vertical expansion at scale $(j-1)$ to $Z^{(j-1)\to\rho}$, we are finally back to the induction hypothesis, at scale $(j-1)$ this time.

\end{enumerate}

Let us now show the following {\em single scale} Mayer bounds.

\begin{Proposition}[Mayer bounds] \label{prop:Mayer-bound}

\begin{enumerate}
\item ({\em vacuum polymers}) Let ${\cal P}_0^{j\to}(\Del^j)$ be the set of vacuum polymers down to scale $j$ containing some fixed interval $\Del^j\in\D^j$. Assume
\BEQ \sum_{\P\in {\cal P}_0^{j\to}(\Del^j)} e^{|\P|-1} |F_{HV}^j(\P)|\le K'\lambda, \label{eq:3.35} \EEQ
where \footnote{Since $\P$ is a vacuum polymer, no high-momentum fields of scale $j$ may be produced at a later stage, hence one may integrate out the field components $\psi^{j\to}$.}  $F^j_{HV}(\P)=\int d\mu(\psi^{j\to}) F^j_{HV}(\P;\psi^{j\to}).$
Then \BEQ |f^{j\to\rho}(\lambda)|\le K'\lambda(1+O(\lambda)) \label{eq:3.36} \EEQ
with the same constant $K'$.

\item (counterterm) Let ${\cal P}_{i,i'}^{j\to}(\Del^j)$ be the set of  polymers down to scale $j$ with exactly two external legs, $\psi_i^{\to(j-1)}$ and $\psi_{i'}^{\to(j-1)}$, displaced into the same point in a fixed interval $\Del^j$. Rewrite $F_{HV,amputated}^j(\P):=\int
    d\mu(\psi^{j\to}) F^j_{HV,amputated}(\P;\psi^{j\to})$, $\P\in {\cal P}_{i,i'}^{j\to}(\Del^j)$ as
\BEQ M^{j(1+\beta_i+\beta_{i'})} F_{HV,amputated}^{j,{\mathrm{rescaled}}}
(\lambda,\P;\psi^{\to(j-1)}),\EEQ
with:
 \begin{itemize}
 \item  $\sum_{\P\in {\cal P}_{i,i'}^{j\to}(\Del^j),\#\{{\mathrm{vertices\ of\ }} \P\}=2}
     F^{j,{\mathrm{rescaled}}}_{HV,amputated}(\lambda;\P)=-K\lambda^2(1+O(\lambda))$;
 \item
\BEA &&  \sum_{\P\in {\cal P}_{i,i'}^{j\to}(\Del^j), \#\{{\mathrm{vertices\ of\ }} \P\} \ge 3} e^{|\P|} F_{HV,amputated}^{j,{\mathrm{rescaled}}}(\lambda,\P) \nonumber\\
&& \qquad =-K'\lambda^3
(1+O(\lambda)+O(M^{-j(1+\beta_i+\beta_{i'})}b_{i,i'}^j)). \nonumber\\ \EEA
    \end{itemize}
Then \BEQ b_{i,i'}^j=K\lambda^2 M^{j(1+\beta_i+\beta_{i'})} (1+O(\lambda)). \EEQ

\end{enumerate}

\end{Proposition}

We shall see in section 6.1 how to obtain estimates which are valid for a succession of Mayer expansions at each scale.

{\bf Proof.}

\begin{enumerate}
\item Fix some interval $\Del^j\in\D^j$ and compute $f^{j\to\rho}(\lambda)$ using Proposition \ref{prop:Mayer} as
\BEQ \sum_{N\ge 1} \sum_{\P_1\in {\cal P}_0^{j\to}(\Del^j),\P_2,\ldots,\P_N
\in {\cal P}_0^{j\to}} \sum_{{\mathrm{trees}} \ \T \ {\mathrm{ over\ }} \ \P_1,\ldots,\P_N}
{\mathrm{Mayer}}(Z_V^{j\to\rho}(\lambda);\T), \EEQ
where
\BEQ |{\mathrm{Mayer}}(Z_V^{j\to\rho}(\lambda);\T)|\le \frac{1}{N!}
\prod_{n=1}^N |F_{HV}^j(\P_n).\EEQ
The $\frac{1}{N!}$ factor is matched by Cayley's theorem, which states that the number of trees over $\P_1,\ldots,\P_N$ with fixed coordination number $(n(\P_i))_{i=1,\ldots,N}$ equals $\frac{N!}{\prod_i (n(\P_i)-1)!}$.
Recall that $\P_{\ell}$ and $\P'_{\ell}$ are necessarily overlapping if
 $\ell\in L(\T)$. Start from
the leaves and  go down the branches of the tree inductively. Let $\P_1,\ldots,\P_{n(\P')-1}$ be the leaves attached onto
one and the same vertex $\P'$ of $\T$. Choose $n(\P')-1$  (possibly
non-distinct) {\em vertices} of $\P'\cap\D^j$ (there are $|\P'\cap\D^j|^{n(\P')-1}$ possibilities), fix their spatial location, $\Del^j_1,\ldots,\Del^j_{n(\P')-1}$, and assume that $\Del_i^j\in\P_i\cap\D^j$.  For each choice of polymer $\P'$, this gives a supplementary factor $\le \left( K'\lambda |\P'\cap\Del^j|\right)^{n(\P')-1}$, to be multiplied by $\frac{1}{(n(\P')-1)!}$ coming from Cayley's theorem. Summing over $n(\P')=2,3,\ldots$, yields $e^{K'\lambda |\P'\cap\Del^j|}-1$,
which is $\le 2^{|\P'|} (K'\lambda)$ for $\lambda$ small enough. By induction, one gets $|f^{j\to\rho}(\lambda)|\le \sum_{h\ge 1} (K'\lambda)^h=K'\lambda(1+O(\lambda))$, where $h$ is the height of the tree.

\item The definition of $b_{i,i'}^j$ and arguments analogous to those used in (1) yield (letting $\tilde{b}_{i,i'}^j:=\lambda^{-2} M^{-j(1+\beta_i+\beta_{i'})} b_{i,i'}^j$)
\BEQ \tilde{b}_{i,i'}^j=-K(1+O(\lambda)+O(\lambda^2 \tilde{b}_{i,i'}^j)),\EEQ
hence the result by the implicit function theorem.

\end{enumerate}

\hfill\eop

{\bf Remark.} Consider the case when the Mayer expansion has surrounded  $\P\in{\cal P}^{j\to}$, a polymer with $\ge N_{ext,max}$ external legs, by a cloud of {\em vacuum} polymers. Then the contribution of the cloud to $F_{HV}^j(\P;\psi)$ may be computed as in Proposition \ref{prop:Mayer-bound} (1), see eq.
(\ref{eq:3.36}), leading at most to an overall mutliplicative coefficient $\le 1+K'\lambda(1+O(\lambda))+\ldots+\frac{(K'\lambda(1+O(\lambda)))^n}{n!}+\ldots=1+O(\lambda)$
per interval $\Del\in\P$ instead of $e$ as in eq. (\ref{eq:3.35}). This is important if $\Del\in
{\bf\Del}_{ext}(\P)$ is crossed by downward dashed lines as in \S 3.3 but contains
no field, for an arbitrary product of factors $e$ down to the scale of the low-momentum fields, possibly produced by successive Mayer expansions, may lead to a disaster. On the other hand (see end of \S 6.1.2) a factor $1+O(\lambda)$ per interval of $\P$ is not a problem.


\section{Power-counting and renormalization}


This section is divided into two parts. The first one gives a quick overview of how divergences in quantum field theory
are discarded by extracting the {\em local part of diverging diagrams}; these ideas come from classical
{\em power-counting} arguments which are recalled here. The whole idea of renormalization is then
to   transfer the sum of these local parts  to the interaction as a counterterm (see \S 3.4).
The second one is an informal discussion of the {\em domination problem} of low-momentum field averages, in particular in the
case of the   $(\phi,\partial\phi,\sigma)$-model). The full treatment of this problem is postponed to \S 6.2.


\subsection{Power-counting and diverging graphs}


\begin{Definition}[Feynman diagrams]  \label{def:Feynman}

Let $\psi=(\psi_1(x),\ldots,\psi_d(x))$ be a multiscale Gaussian field with covariance kernel $C_{\psi}$, and
${\cal L}_{int}=\sum_{q=1}^p K_q \lambda^{\kappa_q} \psi_{I_q}$ be an interaction. Then:

\begin{enumerate}

\item a {\em Feynman diagram}
for this theory is a connected  graph $\Gamma$ (whose lines, resp. vertices are generically denoted by $\ell$, resp. $z$)  with
\begin{itemize}
\item[(i)] external vertices $\vec{y}=(y_1,\ldots,y_{n})$ of type $1,\ldots,p$;
\item[(ii)] internal  vertices $x$ of type $1,\ldots,p$;
\item[(iii)]  internal lines $\ell$ connecting
$z_{\ell}$ to $z'_{\ell}$ with double index $(i_{\ell},i'_{\ell})$.  Since the lines do not have a
prefered orientation, and in order to avoid confusion, one
writes $\ell\simeq (i_{\ell},z_{\ell};i'_{\ell},z'_{\ell})$ or indifferently
 $(i'_{\ell},z'_{\ell};i_{\ell},z_{\ell})$.
For every vertex $z$ of type $q$, one may order the lines leaving or ending in $z$ as $\ell_{z,1}\simeq
(i_1,z;i'_1,z'_1),\ldots
\ell_{z,n(z)}\simeq(i_{n(z)},z;i'_{n(z)},z'_{n(z)})$ so that $n(z)=q$ and $(i_1,\ldots,i_q)=I_q$ if $z$ is an internal vertex, and $n(z)<q$,
$(i_1,\ldots,i_q)\subsetneq I_q$ if $z$ is external.

\item[(iv)] external lines $\ell\simeq (i_{\ell},y_{\ell};i'_{\ell},y'_{\ell})$, where $y'_{\ell}\in\R^D$ is an external point not belonging to $\Gamma$;  $|I_q|-n(y)$ of them per external
vertex $y$ of type $q$. The total number of external lines is $N_{ext}(\Gamma):=
\sum_q\sum_{y\ {\mathrm{of\ type\ }} q} (|I_q|-n(y)).$
\end{itemize}

\smallskip

The evaluation $A(\Gamma)$ of an {\em amputated}  Feynman diagram is  given by
\BEQ A(\Gamma)(\vec{y})= \int \prod_{x}
dx
\left( \prod_{\ell\in L_{int}(\Gamma)} C_{\psi}(i_{\ell},z_{\ell};i'_{\ell},z'_{\ell}) \right), \EEQ
 where $x$ ranges over all internal vertices, and $L_{int}(\Gamma)$
is the set of  all  internal lines.

The evaluation $A(\Gamma)$ of a {\em full}  Feynman diagram is  given by
\BEQ \bar{A}(\Gamma)(\vec{y}')= \int \prod_{y}
dy
\left( \prod_{\ell\in L_{ext}(\Gamma)} C_{\psi}(i_{\ell},y_{\ell};i'_{\ell},y'_{\ell}) \right)
A(\Gamma)(\vec{y}),\EEQ
where $L_{ext}(\Gamma)$
is the set of  all  external lines.

\item A {\em multi-scale  Feynman diagram} $(\Gamma;(j(\ell)) )$
is obtained from $\Gamma$ by choosing a scale $j({\ell})$ for each (internal or external) line
and splitting each vertex into these different scales, with the following constraint,
\BEQ {\mathrm{height}}(\Gamma):=\min\{j(\ell);\ell\in L_{int}(\Gamma)\}-\max\{j(\ell);
\ell\in L_{ext}(\Gamma)\}\ge 0. \EEQ
Note that {\em height}$(\Gamma)$ measures the height of internal lines of $\Gamma$ with respect to the external
lines it is attached to.

 The evaluation $A(\Gamma)$ of a  multi-scale
  Feynman diagram is  given by
\BEA  &&  A(\Gamma;(j(\ell))_{\ell\in L_{int}(\Gamma)})(\vec{y})=  \nonumber\\
&& \qquad \qquad \int \prod_{x}
dx
\left( \prod_{\ell\in L_{int}(\Gamma)} C_{\psi}^{j(\ell)}(i_{\ell},z_{\ell};i'_{\ell},z'_{\ell}) \right), \nonumber\\ \EEA
while

\BEQ \bar{A}(\Gamma;(j(\ell))_{\ell\in L(\Gamma)})(\vec{y}')= \int \prod_{y}
dy
\left( \prod_{\ell\in L_{ext}(\Gamma)} C^{j(\ell)}_{\psi}(i_{\ell},y_{\ell};i'_{\ell},y'_{\ell}) \right)
A(\Gamma;(j(\ell)))(\vec{y}).\EEQ

\end{enumerate}

\end{Definition}

Cluster expansions produce {\em single-scale (horizontal) Feynman} trees, {\em multiplied by some
$G$-polynomial and by the exponential of the interaction},
up to the modification of the measure by the weakening coefficients $s$. In
 the final bounds, see section 6, one bounds the exponential by a constant and
 applies the Cauchy-Schwarz inequality to the rest, yielding a multi-scale
 Feynman {\em diagram}. The scale $j(z)$ of a vertex $z$
is then the scale at which the vertex has been produced.

\medskip

Let us from now on assume that the interaction is {\em just renormalizable}. This means that, for every $q=1,\ldots,p$,
$\sum_{i\in I_q} \beta_i=-D$, or in other words, $\int \psi_{I_q}(x) dx$ is  homogeneous of degree $0$. In that case,
the following very simple power-counting rules hold.

\begin{Proposition}[power-counting] \label{prop:power-counting}

\begin{enumerate}
\item The power-counting of an {\em amputated} Feynman diagram $\Gamma$ is the product $\Lambda^{\omega(\Gamma)}:=
\Lambda^D \prod_z \Lambda^{-D} \prod_{\ell\in L_{int}(\Gamma)}
\Lambda^{-(\beta_{i_{\ell}}+\beta_{i'_{\ell}})}$, see Definition \ref{def:Feynman} for notations, where $\Lambda$
is some large, indefinite constant representing an ultra-violet cut-off.

Then the {\em degree of divergence} $\omega(\Gamma)$ of $\Gamma$ is equal to $D+\sum_{\ell\in L_{ext}(\Gamma)}
\beta_{i_{\ell}}.$

\item  The power-counting of a {\em full} multi-scale Feynman diagram $\Gamma$  is the product
\BEQ M^{\omega_{m.s.}
(\Gamma)}:=
M^{D\cdot {\mathrm{height}}(\Gamma)} \prod_z M^{-Dj(z)} \prod_{\ell\in L(\Gamma)}
M^{-j(\ell)(\beta_{i_{\ell}}+\beta_{i'_{\ell}} )}.\EEQ

The {\em multi-scale degree of divergence}  $\omega_{m.s.}(\Gamma)$ is equal to
\BEQ \omega_{m.s.}(\Gamma)=D \cdot {\mathrm{height}}(\Gamma)+
\sum_z \sum_{\ell\in L_{z}}\beta_{i_{\ell}}(j(z)-j(\ell)),\EEQ
where $L_z$ is the set of {\em internal or external lines} leaving or ending in $z$.

\end{enumerate}

\end{Proposition}

Let us give brief explanations. The power-counting in 1. may be obtained from Definition
\ref{def:multiscale-Gaussian-field} by assuming all internal lines of $\Gamma$ to be of scale
$\frac{\ln \Lambda}{\ln M}$ and summing over all  vertices {\em except one}, due to overall (approximate or
exact) translation invariance, see \cite{FVT}. The multi-scale power counting is a refined version of the previous one,
 which takes into account the scale of the external legs; the definition of the height
of $\Gamma$ is somewhat {\em ad-hoc} and measures the horizontal freedom of movement of $\Gamma$ with respect
to its external legs if these are supposed to be fixed. More precise computations in the spirit of
cluster expansion appear in the final bounds in section 6.

\bigskip

The principle of the  power-counting rule explained in section \S 6.1.2 is to rescale all fields produced at scale $j$
of a multiscale  Feynman diagram
{\em as if  they were
of   scale $j$}. The production scale of a given vertex $z$ being $j(z)$, this leads to a rescaled degree of divergence,
\BEQ \omega_{m.s.}^{{\mathrm{rescaled}}}(\Gamma):=\omega_{m.s.}(\Gamma)-\sum_z
\sum_{\ell\in L_{z,int}}\beta_{i_{\ell}}(j(z)-j(\ell)),\EEQ
where the sum on $\ell$ ranges over all {\em internal lines} leaving or ending in $x$. Hence
\BEQ \omega_{m.s.}^{{\mathrm{rescaled}}}(\Gamma)=D\cdot {\mathrm{height}}(\Gamma) +
\sum_y \sum_{\ell\in L_{y,ext}} \beta_{i_{\ell}}(j(y)-j(\ell)),
\label{eq:omega-rescaled} \EEQ
where the sum on $\ell$ ranges over all {\em external lines} leaving or ending in $y$.

If height$(\Gamma)$ and
 all $j(y)-j(\ell)$ in (\ref{eq:omega-rescaled}) are equal to the same positive constant,  then
$\omega_{m.s.}^{{\mathrm{rescaled}}}(\Gamma)$ is proportional to the naive degree of divergence defined in Proposition
\ref{prop:power-counting} (1).

\begin{Definition}[renormalization]

If the degree of divergence  $\omega(\Gamma)$  of a   Feynman
diagram, resp. multiscale Feynman diagram,  is $\ge 0$, then it needs to be renormalized.

The {\em local part}  of an amputated   Feynman diagram or multiscale Feynman diagram
 is obtained by integrating over all external vertices except one (in order to take into account the global invariance by translation),
\BEQ A(\Gamma;(j(\ell)))(\vec{y}=(y_1,\ldots,y_n))\rightsquigarrow
\int dy_2\ldots dy_n A(\Gamma;(j(\ell)))(\vec{y})=:
{\mathrm{Local}}(A(\Gamma;(j(\ell)))).\EEQ
By invariance by translation again, it does not depend on $y_1$.

Integrating over external points, one obtains the local part of a {\em full}
Feynman diagram,
\BEA  && \bar{A}(\Gamma;(j(\ell)))(\vec{y}')\rightsquigarrow
{\mathrm{Local}}(\bar{A}(\Gamma;(j(\ell)))(\vec{y}'):=  \nonumber\\
&& \qquad {\mathrm{Local}}(A(\Gamma;(j(\ell)))) \int \prod_y dy
\left( \prod_{\ell\in L_{ext}(\Gamma)} C^{j(\ell)}_{\psi}(i_{\ell},y_{\ell};i'_{\ell},y'_{\ell}) \right),\nonumber\\
\EEA
which may be rewritten as
\BEQ {\mathrm{Local}}(\bar{A}(\Gamma;(j(\ell)))(\vec{y}')=\int \prod_y dy
\left( \prod_{\ell\in L_{ext}(\Gamma)} C^{j(\ell)}_{\psi}(i_{\ell},y_1;i'_{\ell},y'_{\ell}) \right) A(\Gamma;(j(\ell)))(\vec{y}).\EEQ
In other terms, all external legs of $\Gamma$ have been displaced at the same arbitrary external vertex, here $y_1$.

The {\em renormalized} amplitude ${\cal R}A$ or ${\cal R}\bar{A}$ of a Feynman diagram or multi-scale Feynman diagram is the difference between
its amplitude and its local part.

\end{Definition}

{\bf Remark.} Using a Fourier transform, the local part is equivalent to the
classical
evaluation at zero external momenta.

A simple Taylor expansion of order 1 yields
\BEA  && {\cal R}\bar{A}(\Gamma;(j(\ell)))(\vec{y}')=\int \prod_y dy \sum_{\ell\in L_{ext}(\Gamma)}
\left( \prod_{\ell'\in L_{ext}(\Gamma),\ell'\not=\ell} C^{j(\ell')}_{\psi}(i_{\ell'},y_1;i'_{\ell'},y'_{\ell'}) \right)\nonumber\\
&& \qquad
\int_0^1 dt (y_{\ell}-y_1) \esper\left[ \partial \psi_{i_{\ell}}^{j(\ell)}(y_1+t(y_{\ell}-y_1))
 \psi_{i'_{\ell}}^{j(\ell)}(y'_{\ell}) \right]\ \cdot\  A(\Gamma;(j(\ell))(\vec{y}). \nonumber\\
 \label{eq:Te1R} \EEA

In principle, this formula shows that the renormalized diagram amplitude ${\cal R}\bar{A}(\Gamma;(j(\ell)))$
comes with a supplementary spring factor
$M^{-{\mathrm{height}}(\Gamma)}$, leading to an equivalent degree of divergence $\omega^*(\Gamma):=
\omega(\Gamma)-1$. Namely, $y_{\ell}-y_1$ should be
 at most of order $M^{-\min\{j(\ell);\ell\in L_{int}(\Gamma)\}}$, while by definition (see section 2)
 $\esper\left[ \partial \psi_{i_{\ell}}^{j(\ell)}(y_1+t(y_{\ell}-y_1))
 \psi_{i'_{\ell}}^{j(\ell)}(y'_{\ell}) \right] $ is of order $M^{j(\ell)}. M^{-(\beta_{i_{\ell}}+\beta_{i'_{\ell}})} $.
 If $\omega(\Gamma)<1-D$, then $\omega^*(\Gamma)<-D$ and the diagram has become
convergent \footnote{Otherwise one should renormalize to a higher order by removing the beginning of the
Taylor expansion of the diagram in the external momenta. We shall not describe this straightforward extension
of the procedure since we shall not use it in our context.}.

Let us give a formal proof of this general statement. This requires a little care because the
covariance is not exponentially decreasing at large distances in our setting, so using systematically
 the Taylor expansion (\ref{eq:Te1R}) does not work. Consider a given multi-scale diagram $(\Gamma;(j_{\ell}))$. We
consider only divergent two-point subdiagrams for the sake of notations.  Start from the
 divergent diagram of lowest scale, $j_{min}$, say, with external legs $y,y'$, and use the Taylor expansion (\ref{eq:Te1R}).
 Choose
 a tree of  propagators of scale $\ge j_{min}$
  connecting $y$ to $y'$ through intermediary points $x_2,\ldots,x_n$, and set $x_1=y,x_{n+1}=y'$, so that
$|y-y'|\le \sum_{m=1}^n |x_m-x_{m+1}|$.
 The divergent diagram with external legs $y,y'$ has been renormalized as
explained above. Letting $j_m$ be the scale of the link $(x_m,x_{m+1})$, one has lost a rescaled factor
$\sum_{m=1}^n d^{j_{min}}(x_m,x_{m+1})=\sum_{m=1}^n M^{-(j_m-j_{min})} d^{j_m}(x_m,x_{m+1})$.
Each term in this sum has obtained a spring factor $M^{-(j_m-j_{min})}$, which shows that the corresponding diagram
with external legs $x_m,x_{m+1}$ is already {\em de facto} superficially renormalized, although subdiagrams of higher
scale may still need renormalization. In the process, it is clear that every propagator of the diagram may appear only
in {\em one} tree of  propagators at most. In the bounds of section 6, this yields an overall factor
per polymer $\P$ which is bounded by $\prod_{\ell\in L(\P)} (1+d^{j(\ell)}(x_{\ell},x'_{\ell}))$, which is easily controlled
by the polynomial decrease of the covariance at large distances.

\medskip

Let us  summarize our brief discussion.

\begin{Definition}[diverging graphs]\label{def:diverging-graphs}

A Feynman graph or multi-scale Feynman graph $\Gamma$ is {\em divergent} if and only if
\BEQ \omega(\Gamma):=D+\sum_{\ell\in L_{ext}(\Gamma)} \beta_{i_{\ell}}\ge 0. \EEQ

We shall call $N_{ext,max}$ the minimum value of $N_{ext}$ such that
{\em every diagram with $\ge N_{ext,max}$ external legs is convergent}

Renormalizing the evaluation of a  Feynman graph or multi-scale Feynman graph $\Gamma$  yields a quantity  ${\cal R}A(\Gamma)$
for which the displaced external legs have obtained a supplementary spring factor, which is equivalent to replacing
one of the scaling dimensions $\beta_{i_{\ell}}$, $\ell \in L_{ext}(\Gamma)$ with $\beta^*_{i_{\ell}}:=\beta_{i_{\ell}}-1$, or globally
$\omega(\Gamma)$ by $\omega^*(\Gamma):=\omega(\Gamma)-1$.

Let $\omega^*_{max}<0$
the maximal value of the set $\{\omega^*(\Gamma)\}$, where $\Gamma$ ranges over the set of all
 Feynman graphs.

\end{Definition}

{\bf Example} ($(\phi,\partial\phi,\sigma)$-model)

In our case $\beta_{\phi}=\alpha$, $\beta_{\partial\phi}=\alpha-1$, $\beta_{\sigma}=-2\alpha$.
Note however that the interaction (see section 5) is such that split vertices are either of type
 $\partial\phi^j\phi^k\sigma^{k'}$ or $\sigma^j \partial\phi^k\phi^{k'}$, with $k'=k\pm 1$, hence
 (assuming $k'=k$ which does not change anything) $\phi$ may not be an external leg of
 a multi-scale Feynman diagram (in particular, the vertex $\phi\partial\phi\sigma$ is {\em not}
 renormalized). Let $N_{\partial\phi}$, resp. $N_{\sigma}$ be the number of external
  $\partial\phi$- or $\sigma$-legs of such a diagram $\Gamma$. The power-counting
  rule yields $\sum_{\ell\in L_{ext}(\Gamma)}
  \beta_{i_{\ell}}=(\alpha-1)N_{\partial\phi}-2\alpha N_{\sigma}<-1$ as soon
  as $N_{\partial\phi}\ge 2$ or $N_{\sigma}\ge 4$, so $N_{ext,max}=4$. However, by symmetry arguments,
  local parts of diagrams with $N_{\partial\phi}$ or $N_{\sigma}$ odd vanish, so
   there remains only the $\sigma$-propagator
   $\langle \sigma_{\pm}(x)\sigma_{\pm}(y)\rangle$ to renormalize.

\bigskip

Let us write down precisely these power-counting rules in the framework of multi-scale cluster expansion, where Feynman diagrams are replaced by polymers.
Recall from \S 3.3 that a low momentum field $(T^{\to(j-1)}\psi_i)^h(x)$ -- or $(T^{\to(j-1)}\del^j\psi_i)^h(x)$ -- has  three scales attached to it. {\em It is the difference
between the two highest ones} -- {\em the production scale $k$ and the dropping scale $j$}, with
$k>j$ -- that counts here. Namely, considering $h<j'<j$, there are two cases:

\begin{itemize}
\item[(i)] either $\tau_{\Del^{j'}_x}$ (the number of $t$-derivatives produced inside $\Del_x^{j'}$)
is $\ge N_{ext,max}$. Then  rescaling the fields to which the $t$-derivatives are applied is
enough to ensure the convergence of the polymer $\P^{j'\to}$ containing $\Del_x^{j'}$;
\item[(ii)] or $\tau_{\Del^{j'}_x}<N_{ext,max}$. Then (by definition) $t_{\Del_x^{j'}}$ is set to $0$, which kills $(T^{\to(j-1)}\psi_i)^h(x)$ or $(T^{\to(j-1)}\del^j\psi_i)^h(x)$. Hence one must count only on the derived fields for the power-counting.
As explained above, this may give a non-negative degree of divergence, in which case one must renormalize by
removing the local part of the polymer.
\end{itemize}

This means that {\em the rescaling spring factor $M^{\beta_i(k-h)}$ of these low-momentum fields must be
split into $M^{\beta_i(k-j)}$ -- ensuring the horizontal fixing of the polymer, and counting in the
right-hand side of eq. (\ref{eq:omega-rescaled}) -- and $M^{\beta_i(j-h)}$ -- which helps control
the accumulation of low-momentum fields as explained briefly in \S 3.3}.


\subsection{Domination problem and boundary term in the interaction}


Assume $\psi_i$ is a multi-scale Gaussian field with scaling dimension
$\beta_i\in
(-D/2,0)$. Then $2\beta_i>-D$ so, by Definition \ref{def:diverging-graphs}, the
associated two-point functions must be renormalized.
If $\psi_i$
occurs in some interaction term $\lambda^{\kappa_q}\psi_{I_q}$, i.e. $i\in
I_q$, then
renormalization produces a scale-dependent counterterm of the form
\BEQ \del {\cal L}(\psi_i;x)=\sum_{j\ge 0} \lambda_j^{2\kappa_q} C_j
M^{(D+2\beta_i)j}
\left((T\psi_i)^{\to j}\right)^2(x), \label{eq:4.15} \EEQ
where $\lambda_j$ is the running coupling constant, and $C_j$ is a
scale-dependent constant.
In this paper we generally assume that $\lambda_j=\lambda$ is not
renormalized. (In the case of our $(\phi,\partial\phi,\sigma)$-model, $\del{\cal L}$ is of the form
$\sum_{j\ge 0}  b^j ((T\sigma)^{\to j})^2$, with $b^j\thickapprox \lambda^2 M^{(1-4\alpha)j}$ and $C_j\thickapprox 1$. See
the counterterm $\del {\cal L}_4$ in section 5 below for details.)
 Separating
the low-momentum field $\psi_i^{\to(j-1)}(x),x\in\Del^j$ into
$\psi_i^{\to(j-1)}(\Del^j)+\del^j
\psi_i^{\to(j-1)}(x)$ as in Definition \ref{def:2.5} produces a secondary field
$\del^j\psi_i^{\to(j-1)}$
whose contributions to the partition function are easily bounded thanks to
the "spring
factor" (see subsection 6.1). Unfortunately, the averaged fields
$\psi_i^{\to(j-1)}(\Del^j)$
do not come with such a spring factor and (unless $\beta_i<-D/2$, see
subsection 6.1 again)
must be {\em dominated apart}, using the positivity of the interaction.
{\em We assume that $C_j>0$.} Then
Lemma \ref{lem:domination} (1)
shows that $|\lambda^{\kappa_q} (T\psi_i)^{\to(j-1)}(\Del^j)|^{n(\Del)}
e^{-\int \del {\cal L}(\psi_i;x) dx}$ is bounded by $O\left(Kn(\Del)^{\half}
C_j^{-\half} M^{-j\beta_i}\right)^{n(\Del)}$, which agrees -- up to
unessential {\em local
factorials} $n(\Del)^{n(\Del)/2}$, see \S 6.1 -- with the correct power-counting, but the
"petit facteur" $\lambda_j^{\kappa_q}$ has been entirely used, in
contradiction with
the general guideline of cluster expansions. So -- unless $C_j^{-\half}$
is small --
one must use a different strategy.

\bigskip

Several strategies have been used, depending on the model. Here things are particularly simple because the mass counterterm of highest
scale, $b^{\rho}\thickapprox \lambda^2 M^{\rho(1-4\alpha)}$, couples with all scales of the field $\sigma$ and is much better than
the term of order $\lambda^2 M^{j(1-4\alpha)}$ appearing in the right-hand side of \ref{eq:4.15}. This is the reason why this counterterm
has been set apart from the counterterm and put into the covariance of $\sigma$ right from the beginning. The supplementary spring factor
$M^{-\half(1-4\alpha)(\rho-j)}$ per  field plays the r\^ole of a ``petit facteur'', except in a small scale
 interval $\rho-q,\ldots,\rho$,
where $q\thickapprox \ln(1/\lambda)$ is $\rho$-independent. Add now e.g.   to the interaction a term  of the
form $M^{-(4n\alpha-1)\rho} \lambda^{\kappa} \sum_{\rho'\le\rho} ||(T\sigma)^{\to\rho'}(x)||^{2n}$
{\em with $2n\ge 4$}, which is homogeneous
of degree $0$, and totally negligible away from the highest scales because of the evanescent coupling coefficient
$M^{-(4n\alpha-1)\rho}\le M^{-(8\alpha-1)\rho}<1$. If $\kappa<2n$, then each  $\sigma$-field in the interaction is coupled to $\lambda^{\kappa/2n}\gg
\lambda$, which makes it possible to dominate the low-momentum fields for the highest scales $\rho-q+1,\ldots,\rho$.
In this sense this term is a (Fourier) boundary term. The term        $\del L_{12}$ in section 5 below plays this r\^ole. The choice of $4n=12$ is arbitrary.


\section{Definition of the model}


{\em As a general rule, we shall denote by $\mu_K$ the Gaussian measure with covariance $K^{-1}$ if $K$
is a positive-definite kernel.}

Recall from the heuristics in section 1.2 that one wants to define as interaction
$\int {\cal L}_{int}(\phi,\sigma)(x)dx=i\lambda \int (\partial A^+(x)\sigma_+(x)-\partial A^-(x)
\sigma_-(x)) dx.$ Integrating over $\R$ and splitting the fields $\phi_1,\phi_2,\sigma_{\pm}$
into their different momentum scales yields a momentum conservation rule; as a result,
products $\partial\phi_1^{j_1}(x)\phi_2^{j_2}(x)\sigma_+^k(x)$ $(j_1\le j_2)$ or $\partial \phi_2^{j_2}(x)
\phi_1^{j_1}(x) \sigma_-^k(x)$ $(j_2\le j_1)$  contribute to the action integral only if the two highest
scales differ by at most $1$ (provided $M\ge 2$, say). The notation $\sum^{adm}_{j_1,j_2,k}$
means that the sum is restricted to this {\em admissible subset}; explicitly,
$\sum^{adm}_{(j_1,j_2,k)\in I} (\cdots)=\sum_{(j_1,j_2,k)\in I\cap I^{adm}}$, with
$I^{adm}=\{(j_1,j_2,k);\ j_1\simeq j_2\ge k\ {\mathrm{or}}\ j_2\simeq k\ge j_1\
{\mathrm{or}}\  k\simeq j_1\ge j_2\},$ where $j\simeq k$ means $j=k$ or $k\pm 1$.
This scale restriction is added explicitly by hand because the dressing procedure (which is not translation-invariant)
spoils momentum conservation.

\begin{Definition}[propagators]

\begin{itemize}
\item[(i)] Let $d\mu^{\to\rho} (\phi)=d\mu_{|\xi|^{1+2\alpha}}(\phi^{\to\rho})$ be the Gaussian measure associated to the
field $\phi^{\to\rho}$ defined in section 2.2.
\item[(ii)] Let $d\mu^{\to\rho}(\sigma_{\pm})=d\mu_{|\xi|^{1-4\alpha}\Id+b^{\rho}}(\sigma^{\to\rho}_{\pm})$ be the  Gaussian measure
 associated to the
fields $\sigma^{\to\rho}_{\pm}$ defined in section 2.2.
\end{itemize}

\end{Definition}

The two-by-two matrix coefficient $b^{\rho}$ is called the {\em renormalized mass coefficient}
of the $\sigma$-field at scale $\rho$. It is equal to the local part at scale $\rho$
 of the two-point function of the $\sigma$-field (see precise definition below). Note that
$d\mu_{|\xi|^{1-4\alpha}\Id+b^{\rho}}(\sigma^{\to\rho})=\frac{1}{Z'} e^{-\half b^{\rho} \int (\sigma^{\to\rho})^2(x)
dx} d\mu_{|\xi|^{1-4\alpha}},$ where $Z'$ is a normalization constant. This quadratic term appears with the {\em opposite}
sign in the interaction Lagrangian ${\cal L}_{int}^{\to\rho}(\phi,\sigma)(x)$ below, so the interacting measure has really
been constructed initially by using a $\sigma$-field with {\em bare} covariance $d\mu_{|\xi|^{1-4\alpha}}$, as explained
in section 1.

\medskip

Note that all summations over scales in the interaction start from the scale $j=1$. Namely,
 it is really meant to cure ultra-violet divergence and absolutely useless in the infra-red region.
However, the $\sigma$-field is not coupled to the infra-red part of the L\'evy area, which yields
a regularized area which is not homogeneous any more. This can easily be corrected by extending the
Fourier partition $(\chi^j)_{j\ge 0}$ to
the infra-red region  $|\xi|\lesssim 1$, with  $|\supp\chi^j|\subset [M^{j-1},M^{j+1}]$ for all
$j\in\Z$. One may then sum over
all scales from $j=\rho$ to $-\rho$, without bothering to renormalize at scales $j<0$,
and let $\rho\to\infty$. We leave
details to the reader and concentrate on the region $|\xi|\gtrsim 1$.

\medskip

In the sequel, expressions such as $b \sigma^2$ are to be understood as a scalar product $(b\sigma,\sigma)=b_{+,+}(\sigma_+)^2+
2b_{+,-}\sigma_+\sigma_-  +b_{-,-}(\sigma_-)^2.$

\bigskip

\begin{Definition}[interaction]

\begin{itemize}

\item[(i)]
Let \BEQ {\cal L}_{int}^{\to\rho}
(\phi,\sigma)(x):={\cal L}_{4}^{\to\rho}(\phi,\sigma)(x)-\frac{b^{\rho}}{2} (\sigma^{\to\rho}(x))^2, \EEQ
 where
{\tiny \BEA && {\cal L}_4^{\to\rho}(\phi,\sigma)(x):=\nonumber\\
&& \II
 \lambda \left\{ D_+ \left(  \sum^{adm}_{1\le j_1, j_2,k\le \rho}
\partial \phi_1^{j_1}(x) \phi_2^{j_2}(x)\sigma_+^k(x) \right) -
D_- \left(  \sum^{adm}_{1\le j_1, j_2,k\le \rho}
\partial \phi_2^{j_2}(x) \phi_1^{j_1}(x)\sigma_-^k(x) \right) \right\}
 \nonumber\\
&& =\II\lambda\left\{  \sum^{adm}_{1\le j_1< j_2\le \rho, 0\le k\le \rho} \! \!
\partial \phi_1^{j_1}(x) \phi_2^{j_2}(x)\sigma_+^k(x) -  \! \! \sum^{adm}_{1\le j_2< j_1\le \rho, 1\le k\le \rho}
\partial \phi_2^{j_2}(x) \phi_1^{j_1}(x)\sigma_-^k(x) \right\}  \label{eq:etoile2} \\
&&+\frac{\II\lambda}{2} \left\{    \sum^{adm}_{1\le j \le \rho, 1\le k\le \rho}
\partial \phi_1^{j}(x) \phi_2^{j}(x)\sigma_+^k(x) -  \sum^{adm}_{1\le j\le \rho, 1\le k\le \rho}
\partial \phi_2^{j}(x) \phi_1^{j}(x)\sigma_-^k(x) \right\}  \label{eq:etoile3}  ;\EEA }

\item[(ii)] (dressed interaction)

Let \BEA &&  {\cal L}_{int}^{\to\rho}(\phi,\sigma;\vec{t})(x):= \nonumber\\
&& \quad
{\cal L}_4^{\to\rho}(\ .\ ;\vec{t})(x)-\frac{b^{\rho}}{2} (t^{\rho}_x)^2 ||(T\sigma)^{\to\rho}(x)||^2
+\del {\cal L}_4^{\to\rho}(\ .\ ;\vec{t})(x)+\del {\cal L}_{12}^{\to\rho}(\ .\ ;
\vec{t})(x), \nonumber\\ \EEA

where:

\BEQ {\cal L}_4^{\to\rho}(\ .\ ;\vec{t})(x)={\cal L}_{4,+}^{\to\rho}(\ .\ ;\vec{t})(x)-
{\cal L}_{4,-}^{\to\rho}(\ .\ ;\vec{t})(x),\EEQ
{\tiny \BEA &&  {\cal L}^{\to\rho}_{4,+}:=\II \lambda D^{+} \left(\sum^{adm}_{1\le j_1,j_2,k\le \rho} \partial (T^{\to\rho}\phi_1)^{j_1}(x)
(T^{\to\rho}\phi_2)^{j_2} (x) (T^{\to\rho}\sigma_+)^k(x) \right. \nonumber\\
&& \left.
+\sum_{2\le\rho'\le \rho} (1-(t_x^{\rho'})^3) \sum^{adm}_{1\le j_1,j_2,k\le \rho'-1}
\partial (T^{\to(\rho'-1)}\phi_1)^{j_1}(x)
(T^{\to(\rho'-1)}\phi_2)^{j_2} (x) (T^{\to(\rho'-1)}\sigma_+)^k(x)  \right) \nonumber\\ \EEA }

and similarly for ${\cal L}_{4,-}^{\to\rho}$
(the Fourier projections $D_{\pm}$ have been defined in section 1);

\BEQ \del{\cal L}_4^{\to\rho}( \ .\ ;\vec{t})(x):=\half  \sum_{2\le\rho'\le \rho}  b^{\rho'-1} \sum_{2\le\rho''\le \rho'}
 (1-(t_x^{\rho''})^2) \left((T\sigma)^{\to(\rho''-1)}(x)\right)^2; \EEQ

\BEA  &&  \del{\cal L}_{12}^{\to\rho}(\ .\ ;\vec{t})(x):= \nonumber\\
&& \quad   M^{-(12\alpha-1)\rho} \lambda^3\left\{ ||(T\sigma)^{\to\rho}(x)||^6 +
\sum_{2\le \rho'\le \rho} (1-(t^{\rho'}_x)^6)
\int_{\Del} ||(T\sigma)^{\to(\rho'-1)}(x)||^6  dx\right\} \nonumber\\
\EEA
for the Euclidean norm $||\sigma(x)||=\sqrt{|\sigma_+(x)|^2+|\sigma_-(x)|^2}$.

By construction ${\cal L}_{int}^{\to\rho}(\phi,\sigma)={\cal L}_{int}^{\to\rho}(\phi,\sigma;\vec{t}=1).$

\item[(iii)]
Let
\BEQ Z_V^{\to\rho}(\lambda):=\int e^{- \int_V {\cal L}_{int}^{\to\rho}(\phi,\sigma;\vec{t}=1)(x)dx }
d\mu^{\to\rho}(\phi)d\mu^{\to\rho}(\sigma),\EEQ
and, more generally,
\BEQ Z_V^{\rho'\to\rho}(\lambda;(\phi^h)_{h<\rho'},(\sigma^h)_{h<\rho'})=\int
e^{-\int_V {\cal L}_{int}^{\to\rho}(\phi,\sigma;\vec{t}=1)(x) dx} d\mu^{\rho'\to\rho}(\phi)d\mu^{\rho'\to\rho}(\sigma)\EEQ
which is a function of the low-momentum components of the fields, considered as external fields.
\item[(iv)] (definition of the renormalized mass coefficient $b^{\rho'}$)
    Fix $b^{\rho'}$ by requiring that
    \BEQ 0=\int dy \frac{\partial^2}{\partial \sigma^{\to(\rho'-1)}(x)
    \partial \sigma^{\to(\rho'-1)}(y) } Z_V^{\rho'\to\rho}(\lambda;(\phi^h)_{h<\rho'},(\sigma^h)_{h<\rho'}) \big|_{\sigma^{\to(\rho'-1)}=\phi^{\to(\rho'-1)}=0}.  \label{eq:def-b} \EEQ

\end{itemize}

\end{Definition}

Note that $\del {\cal L}_4^{\to\rho} (\ .\ ;\vec{t})(x)$   may be seen as a
 dressed interaction as in Definition \ref{def:dressed-interaction}
if one sets $b^{\rho}:=0$ (the counterterm of scale $\rho$ has been treated separately).
The interaction {\em before} dressing therefore
vanishes, which is coherent with the fact that it has not been put into the model from the beginning, but built inductively to compensate local parts of diverging
graphs.

\medskip

Although the actual form of ${\cal L}^{\to\rho}_4$ looks complicated, it is very close to the interaction term
 $\II\lambda\left(\partial {\cal A}^+(x)\sigma_+(x)-
\partial {\cal A}^-(x)\sigma_-(x) \right)$ of section 1, see eq. (\ref{eq:1.12}), (\ref{eq:1.13}), (\ref{eq:Lint}).


\section{Bounds}


%


\subsection{Gaussian bounds}


This paragraph is the backbone of the section, since it provides a means (i) to bound the sum  {\em all} possible Wick
pairings of all possible $G$-monomial associated to a given polymer; (ii) to bound the sum over {\em all} possible
polymers containing some fixed interval $\Del$ at its lowest scale. The general idea (as explained in the introduction)
is that a polymer is connected
either by horizontal cluster links which are polynomially decreasing at large distances, or by vertical inclusion links
which create {\em spring factors}. The computations below for a polymer $\P$ (see \S 6.1.2 and 6.1.3)
give in the end a bound which is of order $\prod_v \lambda^{\kappa}
\prod_{\Del} M^{-\eps}$ for some positive exponents $\kappa,\eps$,
  where the product ranges over all vertices $v$ and intervals $\Del$ of $\P$. This is the general principle of the
bounds for cluster expansions. Both terms $\lambda^{\kappa}$ and $M^{-\eps}$ are called a {\em ``petit facteur par carr\'e''}
(small factor per interval, in French).  It does not include combinatorial factors (see \S 6.1.4), dominated averaged low-momentum
fields (see \S 6.2), and possibly some other terms (see first step in the Proof of Theorem 6.2 in \S 6.3), but all
these are proved to give for each choice of cluster expansion and $G$-polynomial a supplementary factor of the type
$\prod_{\Del} (1+O(\lambda^{\kappa'}))$ for some exponent $\kappa'>0$, which is eaten up by the factor $\prod_{\Del}
M^{-\eps}$.

\bigskip

Since all computations are {\em Gaussian} in this paragraph, we shall take the liberty to
 write $\langle \  (\cdots) \ \rangle$ instead
of $\esper[\ (\cdots)\ ]$, without any risk of confusion.


\subsubsection{Wick's formula and applications}


We first recall the classical Wick formula.

\begin{Proposition}[Wick's formula]

Let $(X_1,\ldots,X_{2N})$ be a (centered) Gaussian vector. Pair the indices $1,\ldots,2N$;
the result may be represented as  a  graph $\F$ with $n$ connected components, linking the vertices
$1,\ldots,2N$ two by two. As in Definition \ref{def:cluster}, we use the pair notation
$\ell=\{i_{\ell},i'_{\ell}\}$ for links. Then
\BEQ \langle X_1\ldots X_{2N}\rangle=\sum_{\F\  {\mathrm{pairing\ of}} \{1,\ldots,2N\} }
\prod_{\ell\in\F} \langle X_{i_{\ell}} X_{i'_{\ell}}\rangle. \EEQ

\end{Proposition}

{\bf Proof.} see e.g. \cite{LeBellac}, \S 5.1.2.

\hfill \eop

\begin{Corollary}[simple Wick bound] \label{cor:Wick-bound}
Let $(X_1,\ldots,X_{2N})$ be a Gaussian vector. Then, for every $K>0$,
\BEQ |\langle X_1\ldots X_{2N}\rangle|\le K^{-N} \prod_{i=1}^{2N-1} \left[ 1+K\sum_{j>i} |\langle
X_i X_j\rangle|
\right].  \label{eq:Wick-bound-bis}\EEQ
In particular,
\BEQ  |\langle X_1\ldots X_{2N}\rangle|\le \prod_{i=1}^{2N} \left[ 1+\sum_{j\not =i} |\langle
X_i X_j\rangle|
\right].  \label{eq:Wick-bound}\EEQ
\end{Corollary}

{\bf Proof.} Expand the right-hand side and use Wick's formula to get eq.
(\ref{eq:Wick-bound-bis}) with $K=1$. The bound with $K\not=1$ may be obtained from the previous one
by a simple rescaling $X_i\to \sqrt{K}X_i$, $i=1,\ldots,2N$. \hfill \eop

The above bound  eq. (\ref{eq:Wick-bound-bis}) depends
 on the ordering of the variables $X_1,\ldots,X_{2N}$, although
$\langle X_1\ldots X_{2N}\rangle$ doesn't, of course. The idea conveyed by this bound
is that it may be important to choose the right order. Similarly, eq. (\ref{eq:Wick-bound-bis})
is clearly optimal when the factors $K\sum_{j>i} |\langle X_i X_j\rangle|$, $1\le i\le 2N-1$, are
of order $1$.

  However this bound is too simple
to apply in most cases, and we shall need refined versions of it using the spatial structure
of the Gaussian  variables. The following lemmas, for the reasons we have just explained, are
to be used after a suitable rescaling.

\begin{Corollary}[Wick bound with spatial structure] \label{cor:Wick-bound-spatial}

\begin{enumerate}

\item (single-scale bound)
Let $(X_{(\Del,n)})_{\Del\in\D^j,1\le n\le N_{max}(\Del)}$ be a Gaussian vector
indexed by $M$-adic intervals of scale $j$. Denote by $I=\{(\Del,n);\Del\in\D^j,1\le n\le N(\Del)\}$ the total
set of indices.  Call {\em connecting pairing} a partial pairing
$\F$ of the indices $(\Del,n)$ such that its {\em spatial projection} $\bar{\F}$ with
vertices $\{\Del\in \D^j;\exists n\le N_{max}(\Del) \ |\  (\Del,n)\in\F\}$ and links
$\{ \{\Del,\Del'\}\in \D^j\times\D^j; \exists n\le N_{max}(\Del), n'\le N_{max}(\Del') \ |\
(\Del,n)\sim_{\F} (\Del',n')\}$ is connected. Fix
some $M$-adic interval $\Del_1\in \D^j$. Then
\BEA  && \sum_{\F\  {\mathrm{connecting\ pairing\ of\ }} I, |\F|=2N,\Del_1\in\bar{\F} }
\prod_{\ell\in \F} |\langle X_{(\Del_{\ell},n_{\ell})}X_{(\Del'_{\ell},n'_{\ell})}\rangle| \nonumber\\
&& \le
\left( 1+ \sup_{\Del\in\D^j} \left( \sum_{n=1}^{N_{max}(\Del)} \sum_{(\Del',n')\in I,(\Del',n')\not=(\Del,n)}
 |\langle X_{(\Del,n)} X_{(\Del',n')}\rangle|\right)
\right)^{3N}. \nonumber\\  \label{eq:6.14} \EEA

\item (single-scale bound, improved version)

More generally, assume $N_{max}(\Del)=\infty$ and $K:(\D^j\times\{1,\ldots,d\})\times(\D^j\times
\{1,\ldots,d\})\to\R_+$ is some kernel, copied an infinite number of times, so that
$K((\Del,pd+i),(\Del',p'd+i'))=K((\Del,i),(\Del',i'))$, with $p,p'=1,2,\ldots$ (see Remark below). Let
$I=\D^j\times\{1,2,\ldots\}$. Conecting pairings of $I$ must be understood modulo the pair identifications
$((\Del,pd+i),(\Del',p'd+i'))\sim((\Del,i),(\Del',i'))$.

Then, for any $\gamma\ge 1$, letting $N(\Del):=\#\{(\Del,n); (\Del,n)\in\F\}$ -- to be interpreted as the number of fields lying in a given interval $\Del$ --,
\BEA &&  \sum_{\F\  {\mathrm{connecting\ pairing\ of\ }} I, |\F|=2N,\Del_1\in\bar{\F} }
\prod_{\ell\in \F} \left( N(\Del_{\ell})N(\Del'_{\ell})\right)^{-\gamma}
 K((\Del_{\ell},i_{\ell}),(\Del'_{\ell},i'_{\ell}))  \nonumber\\
&& \le
\left( 1+ \sup_{\Del\in\D^j}\sum_{\Del'\in\D^j}\sum_{i,i'=1}^d
K((\Del,i),(\Del',i'))
\right)^{3N}. \nonumber\\  \EEA

\item (multi-scale bound)

Let  $K:(\D^{j\to}\times\{1,\ldots,d\})\times(\D^{j\to}\times
\{1,\ldots,d\})\to\R_+ $ be some kernel
indexed by $M$-adic intervals of scale $\ge j$, copied an infinite number of times as in 2.
 We denote once again by $I$ the total set of indices.
 Let $I^k:=\{(\Del,n)\in I; \ \Del\in
\D^k\}$, $k\ge j$  and $I^{\to k'}:=\uplus_{j'=j}^{k'} I^{j'}$, $k'\ge j$. Fix a certain number of (non necessarily
distinct) $M$-adic intervals for each scale $k=j,j+1,\ldots,\rho$, say, $\Del_1^k,\ldots,
\Del_{c_k}^k$ $(k\ge j)$, with  $c_k\ge 0$ $(k>j)$ and $c_j=1$; write for short $\vec{\Del}^{j\to}
=\{(\Del_c^k)_{k\ge j,1\le c\le c_k} \}\subset \D^{j\to}$. Let ${\cal F}^{j\to}(\vec{\Del}^{j\to})$ be the set of multi-scale cluster forests
$\F^{j\to}$ (called:
{\em $\vec{\Del}^{j\to}$-connected multiscale cluster forests}) such that, for
each $j'\ge j$, each vertex of $\F^{j'\to}$ is connected by horizontal cluster links or
inclusion links to some (possibly many) of the selected intervals $(\Del_c^{k'})_{k'\ge j',1\le c\le c_{k'}}$, and the intervals $\Del_1^{j'},\ldots,\Del_{c_{j'}}^{j'}$ are not
connected within $\F^{j'\to}$ by horizontal cluster nor inclusion links. (In other words,
each selected interval $\Del_1^{j'},\ldots,\Del_{c_{j'}}^{j'}$ lies within a different
horizontal cluster and inclusion connected component of $\F^{j'\to}$, and these
$c_{j'}$ connected components exhaust the set of horizontal cluster and inclusion connected components of $\F^{j'\to}$ which contain at least one $M$-adic interval of scale $j'$).
Call {\em $\vec{\Del}^{j\to}$-connecting pairing} a partial pairing
$\F$ of the indices $(\Del,n)$ such that its {\em spatial projection} $\bar{\F}$ has the
same set of vertices and links as some $\vec{\Del}^{j\to}$-connected forest, plus
possibly some supplementary links, possibly reducing the number of connected components.
Then :
\BEA && \sum_{\F \vec{\Del}^{j\to} {\mathrm{-connecting\ pairing\ of\ }} I, |\F|=2N }
 \prod_{\ell\in \F} \left( N(\Del_{\ell})N(\Del'_{\ell})\right)^{-\gamma}
 K((\Del_{\ell},i_{\ell}),(\Del'_{\ell},i'_{\ell})) \nonumber\\
&& \le
\left(1+\max_{k\ge j}  \sup_{\Del^k\in\D^k}   \left[ \sum_{\Del'\in\D^{j\to k}} \sum_{i,i'=1}^d K((\Del^k,i),(\Del',i'))\right. \right. \nonumber\\
&& \left. \left.
 + \sum_{k'=j}^{k-1} \sum_{\Del''\in\D^{j\to k'}}
\sum_{i',i''=1}^d K((\Del^{k'},i'),(\Del'',i'')) \right]\right)^{3N}, \nonumber\\ \label{eq:multiscale-Wick-bound} \EEA
where $\Del^{k'}$ is the unique interval of scale $k'$ such that $\Del^{k'}\supset\Del^k$.

\end{enumerate}

\end{Corollary}

{\bf Remark.} When $N_{max}(\Del)=\infty$, which is due to the fact that the total number of fields in a given
$M$-adic interval, $N(\Del)$, may be  of order $n(\Del)$,
hence unbounded, the bound in 1. is infinite. In practice,
a cluster expansion generates -- thanks to the polynomial decrease in the distance of the covariance of
multi-scale Gaussian fields  -- extremely small factors per interval when $n(\Del)$ is large. The idea
is then to bound $|\langle \psi_i^j(x)\psi_{i'}^j(x')\rangle|$, $x\in\Del,x'\in\Del'$
 by $\frac{1}{(1+d^j(\Del,\Del'))^r}K((\Del,i),(\Del',i'))$, where some of the polynomial decrease
in the distance has been retained in the kernel $K$ (see \S 6.1.2).

{\bf Proof.}

\begin{enumerate}
\item

Consider first the left-hand side of (\ref{eq:6.14}). Consider a connecting pairing $\F$ such that $|\F|=2n$ and
 containing the $M$-adic
interval $\Del_1$, and a spanning tree $\bar{\T}$ of $\bar{\F}$ containing $\Del_1$.
 Associate to $\F$ the following sequence of links and of factors $1$:

-- consider all the pairings of the indices $(\Del_1,n)$, $1\le n\le N_{max}(\Del_1)$ among themselves and
with indices $(\Del',n')$, $\Del'\not=\Del_1$; say (in some arbitrary order), $(\Del_1,n_1)$ pairs with
$(\Del'_1,n'_1)$, $\ldots$, $(\Del_1,n_{N_1-1})$ pairs with $(\Del'_{N_1-1},n'_{N_1-1})$. Insert after
these $N_1-1$ links a factor $1$, signifying that all paired Gaussian variables lying in the interval $\Del_1$
have been exhausted;

-- continue to explore new vertices of $\bar{\F}$ by going along the branches of $\bar{\T}$.
Always insert a factor $1$ after all the pairings of the Gaussian variables lying in a given interval
have been exhausted.

Since
$\bar{\T}$ is connected, all $M$-adic intervals in $\bar{\F}$ and all indices in $\F$ will eventually have
been explored. The number of factors is $\ell(\F)+|\bar{\F}|\le N+|\F|=3N$, to the completed by
the required number of factors $1$ so that there are exactly $3N$ factors.

\medskip

Consider now the right-hand side. Let
\BEQ K_{\Del}:=\sum_{1\le n\le N_{max}(\Del)} \sum_{(\Del',n')\not=(\Del,n)}
|\langle X_{(\Del,n)}X_{(\Del',n')}\rangle| \label{eq:6.12} \EEQ
 and $K_{\emptyset}=1$, and expand
$K^{3N}:=(K_{\emptyset}+\sup_{\Del\in\D^j} K_{\Del})^{3N}$. One gets
\BEA &&  K^{3N}=\sum_p \sum_{N_1<\ldots<N_p} (\sup_{\Del} K_{\Del})^{N_1-1} \ \cdot 1 \ \cdot \nonumber\\ && \qquad \qquad
(\sup_{\Del} K_{\Del})^{N_2-N_1-1} \ \cdot 1  \cdots 1\ \cdot\ (\sup_{\Del}K_{\Del})^{N_p-N_{p-1}-1}\ \cdot \ 1\cdots 1. \nonumber\\ \EEA
Replace the first sequence of $N_1-1$ factors by $K_{\Del_1}^{N_1-1}\le (\sup_{\Del} K_{\Del})^{N_1-1}$ and
expand them, which encodes in particular all possible pairings of the Gaussian variables lying in $\Del_1$.
Consider now an interval $\Del_2\not=\Del_1$ linked by $\bar{\T}$ to $\Del_1$, and replace the
second sequence of factors by $K_{\Del_2}^{N_2-N_1-1}$, and so on. Thus one has encoded all possible
connecting pairings containing the interval $\Del_1$.

\item Considering as in the previous case all the pairings of the indices $(\Del_1,n)$, $n=1,2,\ldots$
for a given pairing $\F$, one gets the factor
\BEA && K_{\Del_1}:=\sum
 \left( (N(\Del_1)N(\Del'_1))^{-\gamma} K((\Del_1,i_1),(\Del'_1,i'_1))\right)\ldots \nonumber\\
&& \qquad \qquad
\left( (N(\Del_1)N(\Del'_{N_1-1}))^{-\gamma} K((\Del_1,i_{N_1-1}),(\Del'_{N_1-1},i'_{N_1-1}))\right),
\nonumber\\ \EEA
where the sum ranges over all possible $(N_1-1)$-uple couplings $((\Del_1,i),(\Del',i'))$ with multiplicities
(note that $\frac{N(\Del_1)}{2}\le N_1-1\le N(\Del_1)$, depending on the number of couplings of fields inside $\Del_1$). Then

\BEA K_{\Del_1}&\le & \left( \sum_{\Del'\in\F^j} \sum_{1\le i,i'\le d} K((\Del_1,i),(\Del',i'))
\sum_{n=1}^{N(\Del_1)}\sum_{n'=1}^{N(\Del')} (N(\Del_1)N(\Del'))^{-\gamma} \right)^{N_1-1} \nonumber\\
&\le & \left( \sum_{\Del'\in\D^j} \sum_{i,i'=1}^d K((\Del_1,i),(\Del',i')) \right)^{N_1-1}.\EEA
Apart from this slight difference, the exploration procedure is the same.

\item Choose a spanning tree of $\bar{\F}\cap\D^{\rho}$, complete it into a spanning tree
of $\bar{\F}\cap\D^{(\rho-1)\to}$, and so on. As in 1., explore the horizontal cluster connected components of
scale $\rho$ starting from the selected intervals $\Del^{\rho}_c$, $1\le c\le c_{\rho}$,
then the connected components of scale $\rho-1$, and so on, down to scale $j$. The only
difference is that two different {\em horizontal cluster} connected components of $\bar{\F}$ of the same
scale $j'$ may be connected from above by inclusion links and horizontal cluster links
of higher scale; in this case, this procedure may not explore all vertices of $\F$.
Fortunately, the bound in eq. (\ref{eq:multiscale-Wick-bound}) gives the possibility, starting from some interval $\Del^k\in\D^k$, to go on to explore all the Gaussian variables
located below $\Del^k$, i.e. in some interval $\Del^{k'}\supset\Del^k$ with $k'<k$.
\end{enumerate} \hfill \eop


\subsubsection{Gaussian bounds for cluster expansions}


We assume here that ${\cal L}_{int}$ is just renormalizable, so that (assuming just for simplicity of notations that its coefficients are scale-independent)
${\cal L}_{int}=K_1\lambda^{\kappa_1}
\psi_{I_1}
+\ldots+K_p\lambda^{\kappa_p}\psi_{I_p}$, where $\kappa_1,\ldots,\kappa_p>0$ and
$\sum_{i\in I_1} \beta_i=\ldots=\sum_{i\in I_p}\beta_i=-D$. Each term
$\psi_{I_1},\ldots,\psi_{I_p}$ is called a {\em vertex} by reference to the Feynman diagram representation (see section 4). Let us recall briefly
that the $G$-monomials are produced:

-- either by horizontal cluster expansions; if $i_{\ell}\in I_q$, $\ell$ being a link at scale $j$,
then $\frac{\del}{\del \psi^j_{i_{\ell}}(x_{\ell})} e^{-\lambda^{\kappa_q}\int \psi_{I_q}(x) dx}$ produces
 $\lambda^{\kappa_q}
\psi_{I_q\setminus\{i_{\ell}\}}(x_{\ell}).$ On the other hand, $\frac{\del}{\del \psi^j_{i_{\ell}}(x_{\ell})}$ may
derivate the low-momentum components of monomials produced at scales $\ge j+1$, which lowers the degree of $G$;

-- or by $t$-derivations acting on $e^{-\int {\cal L}_{int}(x)dx}$,
 yielding (up to $t$-coefficients) some (scale components) of the
$\lambda^{\kappa_q} \psi_{I_q}(x_{\Del})$, $1\le q\le p$, integrated over some $M$-adic interval $\Del$. Again,
$t$-derivations may derivate the monomials produced at scales $\ge j$, which does not change the degree of $G$.

The above products of fields must now be split according to their scale decomposition. Thus one obtains
a certain number $v$ of {\em vertices} split into different scales.

\medskip

We suggest the following notations in order to avoid the proliferation of indices. Make a list $(\psi_1,
\ldots,\psi_d)$, $d=|I_1|+\ldots+|I_p|$,
 of {\em all} the fields involved in the interaction, {\em possibly with repetitions}. Thus the cluster expansion
at scale $j$ generates at the same scale $\lambda^{\kappa_q} \psi_{I_{\ell}}^j(x_{\ell})$, where $I_{\ell}
\subset I_q
\setminus\{i_{\ell}\} \subset \{1,\ldots,d\}\setminus\{i_{\ell}\}$ for some $q\le p$; on the other hand, each
$t$-derivation in an interval $\Del\in\D^j$ generates (up to $t$-coefficients) some $\lambda^{\kappa_q}
\psi^j_{I_{\Del}}(x_{\Del})$, $I_{\Del}\subset I_q$, or (with an extra index $\tau_{\Del}$ for the order
of derivation) $(\lambda^{\kappa_q}\psi^j_{I_{\Del,\tau}}(x_{\Del,\tau}))_{\tau=1,\ldots,\tau_{\Del}}$.

But other field components of scale $j$, lying in some fixed interval $\Del^j\in\D^j$, are produced, either
at an earlier stage $k>j$, in the form of a low-momentum field,
 $\psi^j(x)$ or $\del^k\psi^j(x)$ (secondary field) with
$x\in\Del^k$, $\Del^k\in\D^k$, $\Del^k\subset\Del^j$, or at a later stage $h<j$,
in the form of a high-momentum field,
$Res^h_{\Del^j}\psi^j(x)$, $x\in\Del^h$, where $\Del^h\in \D^h$, $\Del^h\supset\Del^j$.

\bigskip

The general principle of bounds for cluster expansions in quantum field theory (as explained at the beginning
of \S 6.1)  is to (1) use the polynomial decrease in the
distance of the covariance of the field components; (2) find out a ``petit facteur par carr\'e'' (small factor per
cube, or rather per interval in one dimension). This means essentially the following: chose some possibly derivated
interaction term $\lambda^{\kappa_q} \psi_{I_q\setminus\{i_{\ell}\}}(x_{\ell})$, $x_{\ell}\in\Del^j$ or $\lambda^{\kappa_q}
\psi_{I_q}(x_{\Del^j})$ coming from a vertex at scale $j$; the fields $\psi_i^j$, $i\in I_q$ scale like $M^{-\beta_i j}$,
and the integration over the interval $\Del^j\in \D^j$ produces a factor $M^{-j}$ (or $M^{-Dj}$ in general). As for
the cluster expansion at scale $j$, it has produced a factor $C^j_{\psi}(i_{\ell},x_{\ell};i'_{\ell},x'_{\ell})$
which  scales like $M^{-\beta_{i_{\ell}}j}$, times the same quantity with a prime. Supposing
one chooses a splitting of the vertex such that all fields are of scale $j$, then the product of these factors is
$\lambda^{\kappa_q} M^{-Dj} M^{-\sum_{i\in I_q} \beta_i}=\lambda^{\kappa_q}\ll 1$, which is the ``petit facteur''.
Unfortunately the splittings of the vertex produce much more complicated situations; however, the guideline is to
{\em compare the scalings} of the high-momentum fields (rewritten as a sum of restricted fields) and of the low-momentum fields
(possibly rewritten as secondary fields, modulo averaged fields) {\em with the scaling they would produce if they were
of scale $j$}. The high-, resp. low-momentum rescaled fields,  write
\BEQ \psi_i^k(x)=:M^{- \beta_i(k-j)}\psi_i^{k,{\mathrm{rescaled}}}(x)\quad (k>j), \ {\mathrm{resp.}} \  \psi_i^h(x)=: M^{\beta_i(j-h)}\psi_i^{h,{\mathrm{rescaled}}}(x) \quad (h<j),
\EEQ
 yielding positive or negative rescaling spring-factors (depending also on the sign of $\beta$), see after Corollary
 \ref{cor:spring-factor}. The factor $\lambda^{\kappa_q}$ is split into the different scales of the
 vertex, so that each field $\psi_i(x)$, $i=i_1,\ldots,i_q$ is accompanied by a small factor
 $\le \lambda^{\kappa_0/|I_q|}$ for some $\kappa_0>0$.

 Averaged fields must be treated apart and account for the so-called {\em domination problem}; bounding them may
require part of the small factor $\lambda^{\kappa_q}$, so that, generally speaking, the ``petit facteur'' is of order
$\lambda^{\kappa}$, for some $\kappa>0$ but small.

\medskip

Let us first consider the following single scale situation, throwing away all low- or high-momentum
fields for the time being.

\begin{Lemma}  \label{lem:single-scale-bound}

Let $\psi=(\psi_1(x),\ldots,
\psi_d(x))$ be a Gaussian field with $d$ components  such that
\BEQ |C^j_{\psi}(i,x;i',x')|=|\langle \psi^j_i(x) \psi^j_{i'}(x')\rangle| \le K_r
\frac{M^{-j(\beta_i+\beta_{i'})}}{(1+M^j|x-x'|)^r} \label{eq:Cjpsi} \EEQ
for every $r\ge 0$, with some constant $K_r$ depending only on $r$; these bounds hold in particular
if
$\psi_1,\ldots,\psi_d\subset\{(\partial^n\tilde{\psi}_1)_{n\ge 0},\ldots,
(\partial^n\tilde{\psi}_{\tilde{d}})_{n\ge 0}\}$ are derivatives of some independent multiscale Gaussian fields $\tilde{\psi}_1,\ldots,\tilde{\psi}_{\tilde{d}}$.

Consider a horizontal cluster forest $\F^j\in {\cal F}^j$ of scale $j$, and associated
cluster points $x_{\ell},x'_{\ell}$, $\ell\in L(\F^j)$, $x_{\ell}\in\Del_{\ell}$,
$x'_{\ell}\in\Del'_{\ell}$. Choose:

-- for each link $\ell\in L(\F^j)$, a subset $I_{\ell}$ of $\{1,\ldots,d\}\setminus
\{i_{\ell}\}$ and a subset $I'_{\ell}$ of $\{1,\ldots,d\}\setminus
\{i'_{\ell}\}$;

-- for each $M$-adic interval $\Del\in\F^j$, $\tau_{\Del}\le N_{ext,max}+O(n(\Del))$  subsets
$(I_{\Del,\tau})_{\tau=1,\ldots,\tau_{\Del}}$ of $\{1,\ldots,d\}$, and additional
integration points $(x_{\Del,\tau})_{\tau=1,\ldots,\tau_{\Del}}$ in $\Del$.

Such a choice defines uniquely a monomial    $G^{j,j}=G^{j,j}(\F^j;(I_{\ell}),(I'_{\ell});(I_{\Del,\tau}),(x_{\Del,\tau}))$ in the
fields
$\psi_i^j$, $i=1,\ldots,d$ taken at the cluster points $(x_{\ell})$, $(x'_{\ell})$ and the $t$-derivation points $(x_{\Del,\tau})$, namely,
 \BEQ G^{j,j}:=\lambda^{\kappa_0 v(\F^j;G^{j,j})} \left[ \prod_{\ell\in L(\F^j)}  \psi_{I_{\ell}}^j(x_{\ell})  \psi^j_{I'_{\ell}}(x'_{\ell})
 \right] \ \cdot\ \left[ \prod_{\Del\in\F^j} \prod_{\tau=1}^{\tau_{\Del}}
 \psi^j_{I_{\Del,\tau}}(x_{\Del,\tau})\right],
\label{eq:Gjj} \EEQ

where $v(\F^j;G^{j,j}):=2L(\F^j)+\sum_{\Del\in\F^j}\tau_{\Del}$ is the total number of vertices obtained from the horizontal cluster and
the $t$-derivations is (two per horizontal
cluster link, one per $t$-derivation acting on the exponential of the interaction).

Denote by $N_i^j(G^{j,j};\Del),i=1,\ldots,d$
the number
of fields $\psi^j_i(x)$, $x\in\Del$ occurring in $G^{j,j}$ if $\Del\in\F^j$, summing up
to $N^j(G^{j,j};\Del):=\sum_{i=1}^d N_i^j(G^{j,j};\Del)$,  and by
$N_i^j(G^{j,j}):=\sum_{\Del\in\F^j} N_i^j(G^{j,j};\Del)$ the total number of fields $\psi^j_i$
occurring
in $G^{j,j}$. Similarly, we denote by $N_i^j(\F^j)$ the number of half-propagators of
type $i$
in $\F^j$. Note that $N^j(G^{j,j};\Del)=O(n(\Del))$, see
(\ref{eq:Gjj}).

\begin{enumerate}
\item Let
\BEQ I^j_{{\mathrm{Gaussian}}}(\F^j;G^{j,j}):=\int d\mu(\psi^j)
\prod_{\ell\in L(\F^j)} C^j_{\psi}(i_{\ell},x_{\ell};i'_{\ell},x'_{\ell})
\ \cdot
 \ G^{j,j}. \EEQ
Then
\BEQ | I^j_{{\mathrm{Gaussian}}}(\F^j;G^{j,j})|\le K^{|\F^j|} \lambda^{\kappa_0 v(\F^j;G^{j,j})}
\prod_{i=1}^d M^{-j\beta_i(N_i^j(G^{j,j})+N_i^j(\F^j))}.\EEQ

\item

Fix the total number of vertices, $v=v(\F^j;G^{j,j})$, and fix one $M$-adic interval $\Del_1^j\in\D^j$.  Let ${\cal F}_{\Del_1^j}^j\subset{\cal F}^j$ be the subset of
{\em connected} horizontal cluster forests of scale $j$ containing $\Del_1^j$. Consider
the rescaled quantity
\BEQ I^{j,{\mathrm{rescaled}}}_{{\mathrm{Gaussian}}}(\F^j;G^{j,j}):= \prod_{i=1}^d M^{j\beta_i(N_i^j(G^{j,j})+N_i^j(\F^j))} I^j_{{\mathrm{Gaussian}}}(\F^j;G^{j,j})
\EEQ

and:

-- sum over all $\F^j\in{\cal F}_{\Del_1^j}^j$, $(I_{\ell})\subset \{1,\ldots,d\}\setminus
\{i_{\ell}\}$, $(I'_{\ell})\subset \{1,\ldots,d\}\setminus
\{i'_{\ell}\}$, $(I_{\Del,\tau})\subset\{1,\ldots,d\}$;

-- maximize for  $(x_{\ell}),(x'_{\ell}),(x_{\Del,\tau})$, each one ranging over
its associated interval in $\D^j$.

Call $I^{j,{\mathrm{rescaled}}}_{{\mathrm{Gaussian}}}(v;\Del_1^j)$ the result.
Then
\BEQ I^{j,{\mathrm{rescaled}}}_{{\mathrm{Gaussian}}}(v;\Del_1^j)\le (K\lambda^{\kappa_0})^v.
\label{eq:6.17} \EEQ

\end{enumerate}
\end{Lemma}

{\bf Proof.}

\begin{enumerate}
\item

The integral $I^j_{{\mathrm{Gaussian}}}(\F^j;G^{j,j})$ may be evaluated by using
Wick's lemma.
Each choice of contractions leads, using the numerator in the right-hand side of eq.
(\ref{eq:Cjpsi}), to some term with the correct homogeneity
factor,
$\lambda^{\kappa_0 v(\F^j;G^{j,j})} \prod_{i=1}^d M^{-j\beta_i(N_i^j(G^{j,j})+N^j_i(\F^j))}$. Consider
the rescaled fields $\psi^{j,{\mathrm{rescaled}}}_i:=M^{j\beta_i}\psi^j_i$. For reasons to be discussed presently, we shall apply
Corollary \ref{cor:Wick-bound}  to the rescaled fields
$n(\Del_x)^{-\gamma}\psi_i^{j,{\mathrm{rescaled}}}(x)$, for some power $\gamma\ge 1$.

The possibility to introduce this supplementary scaling factor $n(\Del)^{-\gamma}$
comes from the  following argument. Split $\frac{1}{(1+M^j|x-x'|)^r}$ into\\
$\frac{1}{(1+M^j|x-x'|)^{r'}} \frac{1}{(1+M^j|x-x'|)^{r''}}$, with $r=r'+r''$, $r',r''>0$.
 The product of propagators $\prod_{\ell\in L(\F^j)}
C^j_{\psi}(i_{\ell},x_{\ell};i'_{\ell},x'_{\ell})$
contributes,
see denominator in the right-hand side of eq. (\ref{eq:Cjpsi}), a convergence factor
  $K_{r''}^{|L(\F^j)|} \ \cdot\ \prod_{\Del\in\F^j}
\prod_{\Del'\sim\Del}
\left( d^j(\Del,\Del')\right)^{-r''/2}$, for some constant $K_{r''}$ depending
only on $r''$.
 Since the number of intervals $\Del'\in\D^j$
such that $d^j(\Del,\Del')\le n(\Del)/4$ is $\le n(\Del)/2$, this means
that at least
half of the intervals $\Del'\sim\Del$ are at a $d^j$-distance $>n(\Del)/4$
from
$\Del$, so that \footnote{In $D$ dimensions, $\frac{n(\Del)}{4}$ becomes $Kn(\Del)^{1/D}$ and eq. (\ref{eq:5.18})
holds, with different constants.}
\BEQ \prod_{\Del'\sim\Del} \left( d^j(\Del,\Del')\right)^{-r''/2}\le \left[
(n(\Del)/4)^{-r''/2}\right]^{n(\Del)/2}=K^{n(\Del)} n(\Del)^{-K'r''n(\Del)}. \label{eq:5.18}\EEQ

On the other hand, taking into account the $n(\Del_{\ell})^{\gamma}$, resp.
$n(\Del'_{\ell})^{\gamma}$ factors separated from the rescaled fields in cluster intervals
contributes $\prod_{\Del\in \F^j} n(\Del)^{\gamma  N^j(G^{j,j};\Del)}\le
\prod_{\Del\in \F^j} n(\Del)^{\gamma O(n(\Del))}$, a product of so-called {\em local factorials},  which is compensated by the
above convergence factor as soon as $r''$ is chosen large enough.

We may now apply Corollary \ref{cor:Wick-bound}, eq. (\ref{eq:Wick-bound}) to the rescaled fields, which yields, using once again
$N^j(G^{j,j};\Del)=O(n(\Del))$,
 \BEA &&  I^{j,{\mathrm{rescaled}}}_{{\mathrm{Gaussian}}}(\F^j;G^{j,j})\le
K^{\sum_{\Del\in\F^j} n(\Del)} \lambda^{\kappa_0 v(\F^j;G^{j,j})} \nonumber\\
&& \prod_{\Del\in\F^j}
\left[ 1+ (n(\Del))^{-\gamma} \sum_{\Del'\in\F^j}
(n(\Del'))^{-\gamma} \frac{O(n(\Del'))}{(1+d^j(\Del,\Del'))^{r'}}
\right]^{O(n(\Del))}.  \nonumber\\  \label{eq:bd1}\EEA

Now the sum $\sum_{\Del'\in\D^j}
\frac{1}{(1+d^j(\Del,\Del'))^{r'}}$ converges as
soon as $r'>D$.
 Hence each term between
square brackets is bounded by a constant. Since $\sum_{\Del\in\F^j} n(\Del)=2|\F^j|-2
=O(|\F^j|)$, one gets:

\BEQ I^{j,{\mathrm{rescaled}}}_{{\mathrm{Gaussian}}}\le K^{|\F^j|}\lambda^{\kappa_0
v(\Del^j;G^{j,j})}.\EEQ

\item Associate to a connected forest $\F^j\in {\cal F}_{\Del_1^j}^j$ and a monomial
$G^{j,j}$ as in (\ref{eq:Gjj}) its Wick expansion, represented as a sum over a set of connecting
pairings of $\D^j$ as in Corollary \ref{cor:Wick-bound-spatial} (1), except that
$N(\Del)\le d(n(\Del)+\tau_{\Del})+n(\Del)$ -- the number of fields and half-propagators
in the $M$-adic interval $\Del$ -- depends on $\F^j$ and $G^{j,j}$, and is unbounded since
$n(\Del)$ may be arbitrarily large. Hence (to get a finite bound for our sum of Gaussian
integrals) we shall use the $\F^j$-dependent rescaling  by $n(\Del)^{-\gamma}=O(N(\Del)^{-\gamma})$ of the fields defined in 1., at the
price of the  extension of the exploration procedure described in the proof of Corollary \ref{cor:Wick-bound-spatial} (2). Note however that the mapping
$(\F^j,G^{j,j})\mapsto$ connecting pairing is not one-to-one, since a link of the resulting
pairing $\F$ may come either from the links of $\F^j$ or from the pairings of $G^{j,j}$; this
contributes at most a factor 2 per pairing, hence at most $2^{dv/2}.$

Now, the factor $\sum_{\Del'\in\D^j} \sum_{i,i'=1}^d K((\Del,i),(\Del',i'))$ of
Corollary \ref{cor:Wick-bound-spatial} (2) associated to the rescaled fields defined in 1. is bounded up to a constant by $\sum_{\Del'\in\D^j}
\frac{1}{(1+d^j(\Del,\Del'))^r}<\infty$, hence the result.

\end{enumerate}
\hfill \eop

\bigskip

The above arguments  extend easily to single scale Mayer trees of polymers of scale $j$. The new rules are:
\begin{itemize}
\item[(i)] there may be some undetermined number of copies of each interval $\Del^j$, each with a different color;
\item[(ii)] fields in intervals with different colors are uncorrelated;
\item[(iii)] each cluster forest of a given color is connected; one of them (the red one, say) contains a fixed interval, $\Del_1^j$;
\item[(iv)] the different cluster forests are connected by Mayer links. These define a tree structure on the set of colors, and imply for each link between 2 colors, say, red and blue, an overlap between one red interval and one blue interval (chosen at random if they have several overlaps).
 \end{itemize}

The proof of Lemma \ref{lem:single-scale-bound} (2) is the same as before, except that the exploration procedure must now take into account Mayer links. Let $n_{Mayer}(\P')$ be the coordination number of a (red, say) polymer $\P'$ in the Mayer tree. The overlap constraint between $\P'$ and its neighbours $\P_1,\ldots,\P_{n_{Mayer}(\P')-1}$ in the tree splits into multiple overlaps of order $n_i=n_1,\ldots,n_c\ge 1$ between an interval $\Del_i$ in $\P'$ and $n_i$ intervals in $n_i$ neighbouring trees, $\P_{i,1},\ldots,\P_{i,n_i}$, with $n_1+\ldots+n_c=n_{Mayer}(\P')-1$. The exploration procedure at the red interval $\Del_i$ adds $n_i$ to $K_{\Del_i}$, see eq. (\ref{eq:6.12}), corresponding to the number of possible choices of neighbouring trees, but Cayley's theorem, see proof of Proposition \ref{prop:Mayer-bound}, yields a factor $\frac{1}{(n_{Mayer}(\P')-1)!}\le \frac{1}{n_1!\ldots n_c!}$. Summing over all possible values of $n_i$ leads to replacing $K_{\Del_i}$ by
$\sum_{n_i\ge 0} \frac{K_{\Del_i}+n_i}{n_i!}=O(1)$.

\medskip

The whole procedure must be slightly amended to take into account {\em rooted Mayer trees} (see \S 3.4)
connecting possibly an interval $\Del\in{\bf\Del}_{ext}(\P)$, where $\P\in{\cal P}^{j\to}$ is a polymer with $\ge N_{ext,max}$ external legs, to intervals {\em without external legs} of polymers of type 1.
Then one should associate some small power of $\lambda$, say $\lambda^{\kappa}$ with $\kappa\ll\kappa_0$
(at the price of reducing slightly $\kappa_0$ in eq. (\ref{eq:6.17})) to each interval with a
field lying in it, while the intervals $\Del$ of the above type and containing moreover no field
define intervals of a new type (of type 2, say), with no small factor attached to them. The whole
discussion is very similar to the Remark after Proposition \ref{prop:Mayer-bound}. For such intervals, $K_{\Del}\equiv 1$ must be replaced with $1+\sum_{n_i\ge 1} \lambda^{\kappa} \frac{n_i}{n_i!}=1+
O(\lambda^{\kappa})$. As explained at the end of this paragraph, such a factor is not a problem.

\bigskip

We may now give a more general, multiscale bound which takes into account
secondary fields and high-momentum fields.  We rescale the low-momentum fields by reference
to their {\em dropping scale} $j$, and not to their production scale $k$, see \S 3.3, which
leaves outside a supplementary spring factor that will be used to fix the horizontal motion of the polymers as in \S 4.1. As mentioned in the Remarks following Corollary \ref{cor:spring-factor} and Definition \ref{def:3.10}, we shall {\em not} split low-momentum fields $\psi_i^{\to j}$ into a sum (field average)+(secondary field) if $\beta_i<-D/2$. The following definition is valid in all cases:

\begin{Definition}[spring factors] \label{def:beta-tilde}
Let $\tilde{\del}^j\psi^h_i:=\del^j\psi^h_i$ if $\beta_i\ge -D/2$, with $\del^j$ defined
by means of a wavelet admitting $\lfloor \beta_i+\frac{D}{2}\rfloor$ vanishing moments,
 and $\tilde{\del}^j\psi^h_i:=\psi^h_i$  if $\beta_i<-D/2$, so that,
by Corollary \ref{cor:spring-factor},
\BEQ |\langle \tilde{\del}^{j,{\mathrm{rescaled}}} \tilde{\psi}_i^h(x)
 \tilde{\del}^{j',{\mathrm{rescaled}}} \psi_{i'}^h(x')\rangle|\le
K_r. \frac{M^{\tilde{\beta}_i(j-h)} M^{\tilde{\beta}_{i'}(j'-h)}}{(1+d^h(\Del_x^j,\Del_{x'}^{j'}))^r},\EEQ
where $\tilde{\del}^{j,{\mathrm{rescaled}}} \tilde{\psi}_i^h:=
M^{j\beta_i} \tilde{\del}^j \psi_i^h$ is the rescaled field, and
\BEQ \tilde{\beta}_i= \beta_i\ \ (\beta_i<-D/2),\quad \tilde{\beta}_i =\beta_i-\lfloor \beta_i+\frac{D}{2}
+1\rfloor\in[-1-D/2,-D/2)\ \  (\beta_i\ge -D/2).\EEQ

\end{Definition}

{\bf Example.} The $\sigma$-field in the $(\phi,\partial\phi,\sigma)$-model has $\beta_{\sigma}=-2\alpha>-1/2$.
Thus low-momentum fields are severed from their averages, and $\tilde{\beta}_{\sigma}=\beta_{\sigma}-1
=-1-2\alpha$.

\begin{Hypothesis}[high-momentum fields] \label{hyp}

Assume
\begin{itemize}
\item[(i)] either that all $\beta_i$, $i=1,\ldots,d$ are $<0$ \footnote{
which is the case of $\phi^4$-theory for $D>2$ for instance};
\item[(ii)] or, more generally, that there is a  scale constraint on ${\cal L}_{int}$ of the following form: rewriting
${\cal L}_{int}$ as
\BEQ {\cal L}_{int}(x)=\sum_{q\ge 2} \sum_{1\le i_1,\ldots,i_q\le d}  \sum_{j_1\le\ldots\le j_q} K_{i_1,\ldots,i_q}^{j_1,\ldots,j_q} \psi_{i_1}^{j_1}(x)\ldots
\psi_{i_q}^{j_q}(x),\EEQ
then
\BEQ \left( K_{i_1,\ldots,i_q}^{j_1,\ldots,j_q}\not=0\right)\Longrightarrow \left( \beta_{i_1}<0,\beta_{i_1}+\beta_{i_2}<0,\ldots,
\beta_{i_1}+\ldots+\beta_{i_q}<0\right). \label{eq:hypothesis} \EEQ
\end{itemize}

\end{Hypothesis}

This condition on the scales of the {\em low-momentum fields} is of course equivalent to
a condition on the scales of the {\em high-momentum fields} due to the homogeneity of
the vertices.

\medskip

Note that Hypothesis (\ref{eq:hypothesis}) holds true for our $(\phi,\partial\phi,\sigma)$-model, since splitting a vertex leads to
one  low-momentum field, either $\partial\phi$ or $\sigma$, with respective  scaling
exponents $\alpha-1,-2\alpha<0$.  In general, it has the following obvious consequence.

\begin{Lemma} \label{lem:beta-gamma}
Assume Hypothesis (\ref{eq:hypothesis}) holds, and fix $I_q=(i_1,\ldots,i_q)$.  Choose
$\eps>0$ such that
\BEQ \eps< \min\left( |\beta_{i_1}|,|\beta_{i_1}+\beta_{i_2}|,\ldots,
|\beta_{i_1}+\ldots+\beta_{i_q}|\right) \EEQ
 whenever there exists $j_1\le \ldots\le j_q$ such that $K_{i_1,\ldots,i_q}^{j_1,\ldots,j_q}\not=0$, and let
\BEQ \gamma_i:=\frac{\eps}{q}-\beta_i, \quad i=i_1,\ldots,i_q; \quad \gamma_I:=\sum_{i\in I_q} \gamma_i. \label{eq:gamma} \EEQ

Then $\beta_i+\gamma_i>0$ for all $i$, and $\gamma_{i_{q'}}+\ldots+\gamma_{i_q}<D$ for every $q'=1,\ldots,q$.
\end{Lemma}

{\bf Proof.} Since $\beta_{i_1}+\ldots+\beta_{i_q}=-D$,
\BEQ \gamma_{i_{q'}}+\ldots+\gamma_{i_q}<\eps+(D+\beta_{i_1}+\ldots+\beta_{i_{q'-1}})<D.\EEQ \hfill \eop

\bigskip

Under the above Hypothesis, one has a multiscale generalization of Lemma \ref{lem:single-scale-bound} by considering the contribution of {\em all} fields of some fixed scale $j$. We adopt the following
convention. If a vertex $\psi_{I_q}(x)$ is split into $i$ fields of scale $j$ and $|I_q|-i$
 fields of scale $\not=j$, then it contributes a (fractional number) of vertices, $\frac{i}{|I_q|}$,
 of scale $j$. This implies a small factor $\lambda^{\kappa_0 v}$ for the Gaussian bounds at scale $j$,
 where $v$ (the total number of vertices at scale $j$) is a fraction with bounded denominator.

\begin{Lemma}[multiscale generalization] \label{lem:multiscale-bound}

Assume $\psi_1,\ldots,\psi_d\subset\{(\partial^n \tilde{\psi}_1)_{n\ge 0},\ldots,
(\partial^n \tilde{\psi}_{\tilde{d}})_{n\ge 0}\}$ are derivatives of some independent
multiscale Gaussian fields $\tilde{\psi}_1,\ldots,\tilde{\psi}_{\tilde{d}}$.
 Fix
some constant $\kappa_0\in(0,1)$ as in the previous lemma, as well as some reference
scale $j_{min}\le j$.

\begin{enumerate}
\item For each $k\ge j$, consider a horizontal cluster forest $\F^k\in{\cal F}^k$ and
 associated cluster points $x_{\ell^k},x'_{\ell^k}$, and choose subsets $(I_{\ell^k}),
(I'_{\ell^k}),(I_{\Del^k,\tau}), (x_{\Del^k,\tau})$ as in the previous lemma.   Do the same for each $h<j$, and choose
for each point $x=x_{\ell^h}$, $x'_{\ell^h}$ or $x_{\Del^h,\tau}$ a {\em restriction interval}
$\Del^j=\Del^j_{\ell^h},\Del^{'j}_{\ell^h},\Del^j_{\Del^h,\tau}$ such that $\Del^j\subset \Del_{\ell^h},
\Del'_{\ell^h}$ or $\Del^h$ respectively.

Such as choice
defines uniquely a monomial $G^{j,j}$ as before, and one more monomial per different scale,
\BEQ G^{j,k}:= \lambda^{\kappa_0 v(\F^j;G^{j,k})} \left[\prod_{\ell^k\in L(\F^k)}
 \tilde{\del}^k\psi_{I_{\ell^k}}^j(x_{\ell^k})
\tilde{\del}^k\psi^j_{I'_{\ell^k}}(x'_{\ell^k})
 \right]  \left[ \prod_{\Del^k\in\F^k} \prod_{\tau=1}^{\tau_{\Del^k}}
\tilde{\del}^k\psi^j_{I_{\Del^k,\tau}}(x_{\Del^k,\tau})\right].
\label{eq:delGk} \EEQ
for $k>j$, and, see eq. (\ref{eq:gamma}) for notations,
\BEA &&  G^{j,h}:=\lambda^{\kappa_0 v(\F^j;G^{j,h})} \left[\prod_{\ell^h\in L(\F^h)} \sum_{\Del^j\subset\Del_{\ell^h}} \sum_{\Del^{'j}\subset
\Del'_{\ell^h}} \right. \nonumber\\
&& \left. \qquad M^{-\gamma_{I_{\ell^h}}(j-h)}
 Res^h_{\Del^j} \psi_{I_{\ell^h}}^j(x_{\ell^h}) M^{-\gamma_{I'_{\ell^h}}(j-h)}
Res^h_{\Del^{'j}}\psi^j_{I'_{\ell^h}}(x'_{\ell^h})
 \right] \ \cdot \nonumber\\
&&  \cdot\ \left[ \prod_{\Del^h\in\F^h}  \prod_{\tau=1}^{\tau_{\Del^h}} \sum_{\Del^j\subset\Del^h}
 M^{-\gamma_{I_{\Del^h,\tau}}(j-h)}
Res^h_{\Del^j}\psi^j_{I_{\Del^h,\tau}}(x_{\Del^h,\tau})\right]
\label{eq:Gh} \EEA
for $h<j$, where the $\Del^j$, $\Del^{'j}$ are restriction intervals.

Let $v(\F,G^j):=v(\F^j;G^{j,j})+\sum_{k>j} v(\F^k; G^{j,k})+\sum_{h<j} v(\F^j;G^{j,h})$ be the total
number of vertices. Let finally $G^j:=G^{j,j} \prod_{k>j} G^{j,k} \prod_{h<j} G^{j,h}$ and
\BEQ I^j_{{\mathrm{Gaussian}}}(\F;G^{j}):= \int d\mu(\psi^j)
\prod_{\ell\in L(\F^j)} C^j_{\psi}(i_{\ell^j},x_{\ell^j};
i'_{\ell^j},x'_{\ell^j})
   \ G^{j}. \EEQ
Then
\BEQ I^j_{{\mathrm{Gaussian}}}(\F;G^{j})\le K^{|\F|} M^{-(D+|\omega^*_{max}|)\# {\bf\Del}^{j_{min}\to}}
\prod_{a\ge j_{min}} \prod_{i=1}^d M^{-a\beta_i(N_i^a(G^{j,a})+N_i^a(\F^a))},\EEQ
where $a$ stands either for a low-momentum scale $h\le  j$ or a high-momentum scale $k>j$, and
$\omega_{max}<0$ is an in Definition \ref{def:diverging-graphs}.

\item Fix the total number of vertices, $v:=v(\F, G^{j})$. Define
${\cal F}^{j_{min}\to}(\vec{\Del}^{j_{min}\to})$ as in Corollary \ref{cor:Wick-bound-spatial} (3).  Consider, similarly to the previous lemma, the rescaled quantity
\BEA &&  I^{j,\mathrm{rescaled}}_{{\mathrm{Gaussian}}}(\F;G^{j}):= \nonumber\\
&& \qquad
 \prod_{i=1}^d M^{j\beta_i(N^j_i(G^{j,j})+N_i^j(\F^j))} \prod_{a\ge j_{min},a\not=j}
M^{a\beta_i N^j_i(G^{j,a})}  I^j_{{\mathrm{Gaussian}}}(\F; G^{j}). \nonumber\\
\EEA

and:

-- sum over all $\F\in{\cal F}^{j_{min}\to}(\vec{\Del}^{j_{min}\to})$, $(I_{\ell^k})\subset \{1,\ldots,d\}\setminus
\{i_{\ell^k}\}$, $(I'_{\ell^k})\subset \{1,\ldots,d\}\setminus
\{i'_{\ell^k}\}$, $(I_{\Del^k,\tau})\subset\{1,\ldots,d\}$ ($k\ge j$), and similarly for $h<j$;

-- sum over all possibly choices of the restriction intervals $\Del_j$;

-- maximize for  $(x_{\ell^k}),(x'_{\ell^k}),(x_{\Del^k,\tau})$, each one over
its associated interval in $\D^k$, $k\ge j$, and for $(x_{\ell^h}), (x'_{\ell^h}), (x_{\Del^h,\tau})$, $h<j$,
each one ranging over its associated restriction interval in $\D^j$ (and not $\D^h$!).

Call $I^{j,{\mathrm{rescaled}}}_{{\mathrm{Gaussian}}}(v;\vec{\Del}^{j_{min}\to})$ the result.
Then
\BEQ I^{j,{\mathrm{rescaled}}}_{{\mathrm{Gaussian}}}(v;\vec{\Del}^{j_{min}\to})\le
M^{-(D+|\omega^*_{max}|)\# {\bf\Del}^{j_{min}\to}}(K\lambda^{\kappa_0})^v.
\label{eq:multiscale-Gaussian-bound} \EEQ

\end{enumerate}
\end{Lemma}

{\bf Proof.}

Consider first the factor $M^{-(D+|\omega_{max}|)\# {\bf\Del}^{j_{min}\to}}$. It comes
from the part of the rescaling spring factors used for fixing horizontally the polymers (see
\S 4.1). Let $\P^{j\to}_c$ be one of the connected components of the multi-scale forest at scale $j$,
and ${\bf\Del}^j_c\subset{\bf\Del}^{j_{min}\to}\cap\D^j$ be its intervals of scale $j$ with external
legs. Then the rescaling spring factors of the corresponding low-momentum fields $(T\psi_{i_{n_{ext}}})^{\to(j-1)}(x_{n_{ext}})$, $n_{ext}=1,2,\ldots,N_{ext}(\P^{j\to}_c)$, yield
a factor $\prod_{n_{ext}=1}^{N_{ext}(\P^{j\to}_c)} M^{\beta^*_{n_{ext}}}$ when going down from scale $j$
to scale $j-1$, where $\beta^*_{n_{ext}}=\beta_{n_{ext}}$ or $\beta_{n_{ext}}-1$, see \S 4.1, are
such that $\sum_{n_{ext}=1}^{N_{ext}(\P^{j\to}_c)} \beta^*_{n_{ext}}\le -D-|\omega^*_{max}|$.

\medskip

Let us now prove 2. directly, since 1. is a weaker form of 2. Note that by the Cauchy-Schwarz type inequality $\int |fgh|\le (\int |f|^3 \int |g|^3
\int |h|^3)^{1/3}$, one may separate low-momentum fields from high-momentum fields and from the fields produced at scale $j$, to which the previous lemma applies \footnote{Of course, fields produced at scale $j$ may also be treated on an equal foot with high-momentum fields for instance.}.

Let us first consider {\em low-momentum fields}.
 We use the same rescaling as
in the proof of the previous lemma, namely, we consider the rescaled fields
$n(\Del^k_x)^{-\gamma}\tilde{\del}^{k,{\mathrm{rescaled}}} \psi_i^j,$ with $k>j$. The ${\bf\Del}^{j_{min}\to}$-connecting pairing associated to $(\F,G^j)$ has links of two types:

\begin{itemize}
\item[(i)] links due to pairings $\langle \psi_i^j \psi_{i'}^j\rangle$, or, more or less equivalently, cluster links $C^j_{\psi}(i_{\ell^j},x_{\ell^j}; i'_{\ell^j},x'_{\ell^j})$ {\em of scale} $j$;
\item[(ii)] {\em cluster links}  to the propagators   $C^a_{\psi}(i_{\ell^a},x_{\ell^a}; i'_{\ell^a},x'_{\ell^a})$ {\em of scale} $a\not=j$, $a=k>j$ in the specific case of low-momentum fields.
\end{itemize}

Cluster links of scale $k>j$ (or $h<j$) contribute a factor $\frac{1}{(1+d^k(\Del^k,(\Del')^k))^r}$, which is required both for the bound on
$I^j$ and for that on $I^k$. Since $r$ is arbitrary, one chooses it large enough and splits the above factor among the different scales of the vertices. On the other hand, the scaling of the propagators $C^k_{\psi}$ or
$C^h_{\psi}$ is left for the computation of $I^k$ or $I^h$. Note the possible existence of {\em chains}
of propagators of scale $k$ connecting two vertices with low-momentum fields of scale $j$; summing over
all possible chains yields the same factor of order $\frac{1}{(1+d^k(\Del^k,(\Del')^k))^r}$ as soon
as $r>D$.

By Definition  \ref{def:beta-tilde},  the term between square brackets in eq. (\ref{eq:multiscale-Wick-bound}) is bounded up to a constant (see proof of
Corollary \ref{cor:Wick-bound-spatial} (3))  by $A_{cluster}(\Del^k)+A^{low}(\Del^k)$, where
\BEQ A_{cluster}(\Del^k):=\sum_{i,i'=1}^d
 \sum_{k'=j}^{k} M^{\tilde{\beta}_i(k'-j)} M^{\tilde{\beta}_{i'}(k'-j)}
 \sum_{\Del'\in \D^{k'}}  \frac{1}{(1+d^{k'}(\Del^{k'},\Del'))^r} \EEQ
 where $\Del^{k'}$ is the unique interval of scale $k'$ such that $\Del^{k'}\supset\Del^k$, and
\BEQ A^{low}(\Del^{k}):=\sum_{i,i'=1}^d
 \sum_{k'=j}^{k} M^{\tilde{\beta}_i(k'-j)}  \sum_{k''=j}^{k'}
M^{\tilde{\beta}_{i'}(k''-j)}
\sum_{\Del''\in \D^{k''}} \frac{1}{(1+d^j(\Del^{k'},\Del''))^r}. \label{eq:5.31} \EEQ

The above sum, $\sum_{\Del^{k''}\in\D^{k''}} \frac{1}{(1+d^j(\Del^{k'},\Del^{k''}))^r}$, is of order
 $M^{D(k''-j)}$, which yields
\BEQ A^{low}(\Del^{k})\le K \sum_{i,i'=1}^d A^{low}_{\tilde{\beta}_i,\tilde{\beta}_{i'}}, \quad
 A^{low}_{\tilde{\beta}_i,\tilde{\beta}_{i'}}:=  \sum_{k'=j}^{\rho} M^{\tilde{\beta}_i(k'-j)} \sum_{k''=j}^{k'}
M^{(\tilde{\beta}_{i'}+D)(k''-j)}, \label{eq:AbetaD}\EEQ
while clearly $A_{cluster}(\Del^k)$ is finite since $\tilde{\beta}_i,\tilde{\beta}_{i'}<0$.

 \medskip

Let us finish the proof with the assumption that $D=1$. There are 3 different cases:

-- either $\beta_{i'}<-1$; then $\tilde{\beta}_{i'}=\beta_{i'}$ and $\sum_{k''=j}^{k'}
M^{(\beta_{i'}+1)(k''-j)}=O(1)$, $\sum_{k'=j}^{\rho} M^{\tilde{\beta}_i(k'-j)}=O(1)$ since
$\tilde{\beta}_i<0$;

-- or $-1\le \beta_{i'}<-\half$, resp. $0\le \beta_{i'}<\half$: then $\tilde{\beta}_{i'}=\beta_{i'}$, resp.
$\beta_{i'}-1$, and $\sum_{k''=j}^{k'} M^{(\tilde{\beta}_{i'}+1)(k''-j)}=O((k'-j)
M^{(\tilde{\beta}_{i'}+1)(k'-j)})$, $\sum_{k'=j}^{\rho} (k'-j) M^{(\tilde{\beta}_i+
\tilde{\beta}_{i'}+1)(k'-j)}=O(1)$ since $ \tilde{\beta}_i,\tilde{\beta}_{i'}<-\half$;

-- or $-\half\le \beta_{i'}<0$, resp. $\half\le \beta_{i'}<1$; then $\tilde{\beta}_{i'}=\beta_{i'}-1$,
resp. $\beta_{i'}-2$, and $\sum_{k''=j}^{k'} M^{(\tilde{\beta}_{i'}+1)(k''-j)}=O(1)$, while the sum
over $k'$  converges as in the first case.

\medskip

The simpler case $D\ge 2$ is left to the reader (there are only 2 subcases: $\beta<-D/2$, or $\beta\ge
-D/2$, the latter subcase to be  split according to the value of $\lfloor \beta+\frac{D}{2}\rfloor$).

\bigskip

Consider now {\em high-momentum fields}, produced at a scale $h<j$ (or $h\le j$). By the same method, one ends up with the following quantity to bound instead of
eq. (\ref{eq:5.31}):
\BEQ A^{high}(\Del^j):=\sum_{i,i'=1}^d \sum_{h'=j_{min}}^h M^{-(\beta_i+\gamma_i)(j-h')} \sum_{h''=j_{min}}^{h'}
M^{-(\beta_{i'}+\gamma_{i'})(j-h'')} \sum_{\Del^{'j}\in\D^j} \frac{1}{(1+d^j(\Del^j,\Del^{'j}))^r},\EEQ
where $\Del^j,\Del^{'j}$ range over restriction intervals, plus the finite  term
$A_{cluster}(\Del^h)$ due to cluster links of scale $h$ as before. Now $\sum_{\Del^{'j}\in\D^j}
\frac{1}{(1+d^j(\Del^j,\Del^{'j}))^r}<K$ and $\beta_i+\gamma_i,\beta_{i'}+\gamma_{i'}>0$ by Lemma \ref{lem:beta-gamma},
hence $A^{high}(\Del^j)$ is bounded by a constant.

\hfill\eop

Note that a small part of the rescaling spring factor $M^{\tilde{\beta}_i(k-j)}$ ($\tilde{\beta}_i<0$) for low-momentum
fields may be used to obtain a factor $M^{-\eps}<1$ for $\eps>0$ small enough per {\em interval}
$\Del$ belonging to a fixed polymer $\P$ -- in particular in {\em empty} intervals where there is no field. This simple remark
is essential in the sequel since various estimates yield a factor $1+O(\lambda^{\kappa})$ per
interval, which may be compensated by this (not so) small factor $M^{-\eps}$. On the other hand, to each
{\em vertex} or equivalently to each {\em non-empty interval} -- or to each {\em field} -- is associated
a factor of order $\lambda^{\kappa}$, which may be arbitrarily small. This is a general principle
for the bounds to come now.


\subsubsection{Gaussian bounds for polymers}


It remains to be seen how these Gaussian estimates, valid for each scale, combine to give  estimates for the (Mayer-extended)
polymer
evaluation functions defined in subsection 3.4. Note that
(rescaled) low-momentum field averages have been left out; the domination estimates in subsection 6.2 prove that it is
possible in the case of the $(\phi,\partial\phi,\sigma)$-model  to bound
them  while leaving a {\em small factor  per field},  $\lambda^{\kappa_0+\kappa'_0}$, $\kappa_0,\kappa'_0>0$,  negligible with respect to
 $\lambda^{\kappa_0}$.

We take this into account for our next Theorem by choosing a small enough exponent $\kappa_0$.

\begin{Theorem}  \label{th:6.1.3}

Fix some reference scale $j_{min}\ge 0$ and some exponent $\kappa>0$.

 Let ${\cal P}_0^{j_{min}\to}(\Del^{j_{min}})$ be the set of vacuum (Mayer-extended)
polymers down to scale $j_{min}$ containing some fixed interval $\Del^{j_{min}}$ of scale $j_{min}$.  Let
{\em Eval}$_{\kappa_0}(\P)$, $\P\in {\cal P}_0^{j_{min}\to}(\Del^{j_{min}})$  be the sum over all
 multi-scale splitting of vertices into cluster forests $\F^j$ extending over $\P$ and monomials $G^j$
of the product integrals $\prod_{j=j_{min}}^{\rho} \int d\mu(\psi^j) \prod_{\ell\in L(\F^j)} C^j_{\psi}(i_{\ell^j},x_{\ell^j};
i'_{\ell^j},x'_{\ell^j}) G^j$, all fields in $G^j$ being accompanied as before with the factor $\lambda^{\kappa_0}$.
Then
\BEQ \sum_{\P\in {\cal P}_0^{j_{min}\to}(\Del^{j_{min}});\ \#\{{\mathrm{vertices\ of\ }}\P=v\}}
  \left|  {\mathrm{Eval}}_{\kappa_0+\kappa'_0}(\P)\right|\le
(K\lambda^{\kappa_0})^v. \label{eq:Eval1} \EEQ
In particular, for $\lambda$ small enough,
\BEQ  \sum_{\P\in {\cal P}_0^{j_{min}\to}(\Del^{j_{min}})}
  \left|  {\mathrm{Eval}}_{\kappa_0+\kappa'_0}(\P)\right|<\infty.\EEQ

\end{Theorem}

{\bf Proof.}

 Let ${\cal P}_0^{j_{min}\to}({\bf\Del}^{j_{min}\to})$ be the set of  vacuum (Mayer-extended)
polymers down to scale $j_{min}$ with fixed set of intervals ${\bf \Del}^{j_{min}\to}$ as in the previous Lemma. Note that the horizontal fixing scaling factor $M^{-(D+|\omega^*_{max}|)\# {\bf\Del}^{j_{min}\to}}$ makes it possible to sum over all inclusion links of the polymers. Namely,
each inclusion link $\Del^j\subset\Del^{j-1}$ -- implying necessarily some $t$-derivative in $\Del^j$ -- produces a factor $M^D$ due to the choice of $\Del^j$ among
the $M^D$ intervals $\Del\in\D^j$ such that $\Del\subset\Del^{j-1}$, which is compensated by the factor
$M^{-(D+|\omega^*_{max}|)}$ attached to $\Del^j$.

\medskip

Hence it suffices to prove eq. (\ref{eq:Eval1}) for $\P$ ranging in the set ${\cal P}_0^{j_{min}\to}(
{\bf\Del}^{j_{min}\to})$.

\bigskip

 Let $A^{j_{min}\to}_v({\bf\Del}^{j_{min}\to}):=
 \sum_{\P\in {\cal P}_0^{j_{min}\to}({\bf \Del}^{j_{min}\to});
\#\{ {\mathrm{vertices\ of\ }}\P\}=v} \left| {\mathrm{Eval}}_{\kappa_0+\kappa'_0}(\P) \right|.$
 Split the small factor per vertex into $\lambda^{\kappa_0}\lambda^{\kappa'_0}$,
 and set apart $\lambda^{\kappa_0}$ to get a global homogeneity factor
$ \lambda^{\kappa_0 v}$. Then
\BEQ A^{j_{min}\to}_v({\bf\Del}^{j_{min}\to})\le \lambda^{\kappa_0 v}
\sum_{\P\in {\cal P}_0^{j_{min}\to}({\bf\Del}^{j_{min}\to})}
\left| {\mathrm{Eval}}_{\kappa'_0}(\P)\right|,  \EEQ  where the total
number of vertices is now unrestricted.

Let $\P\in {\cal P}_0^{j_{min}\to}({\bf\Del}^{j_{min}\to}) $. Consider some product of fields of
type $q$ produced in some interval $\Del^j$ of scale $j$, $\psi_{I_q}(x_{\Del^j,\tau})$ or $\psi_{I_q\setminus\{i_{\ell}\}}(x_{\ell})$,
interpreted as some pairing of $\psi_{I_q}(x_{\ell})$ with $\psi_{i'_{\ell}}(x'_{\ell})$,
and:

-- choose some non-empty subset $I^{high}=(i_{q'},\ldots,i_q)\subset I_q$;

-- choose some high-momentum scales $(k_i)_{i\in I^{high}}$, with
 $j\le k_{i_{q'}}\le\ldots \le k_{i_q}$ as in Hypothesis \ref{hyp}, and
restriction intervals $(\Del_i^{k_i})_{i\in I^{high}}$, $\Del_i^{k_i}\in\D^{k_i}$,
$\Del_i^{k_i}\subset \Del^j$;

-- letting $I^{low}:=I_q\setminus I^{high}$, choose some low-momentum scales
$(h_i)_{i\in I^{low}}$, $h_i<j$.

Then \BEQ \psi_{I_q}:=\sum_{I^{high}\subset I} \sum_{(k_i)} \sum_{\Del_i^{k_i}}
\sum_{(h_i)} \left[ \prod_{i\in I^{high}} Res^j_{\Del^{k_i}}\psi_i^{k_i}\right]
\ \cdot \ \left[  \prod_{i\in I^{low}} \psi_i^{h_i}\right] \EEQ
is the decomposition of $\psi_{I_q}$ into all possible splittings. Any given splitting of $\psi_{I_q}$ is supported on an
$M$-adic interval $\cap_{i\in I^{high}} \Del^{k_i}$ of size bounded by
\BEA M^{-k_{i_q} D}&=&M^{-jD}\cdot M^{-(k_{i_{q'}}-j)D}\ldots M^{-(k_{i_q}-k_{i_{q-1}})D} \nonumber\\
&\le & M^{-jD} M^{-(\gamma_{i_{q'}}+\ldots+\gamma_{i_q})(k_{i_{q'}}-j)} \ldots M^{-\gamma_{i_q}(k_{i_q}-k_{i_{q-1}})} \nonumber\\
 &=& M^{-jD} M^{-\sum_{i\in I^{high}} \gamma_i(k_i-j)} \EEA
 by Lemma \ref{lem:beta-gamma}.

  Fixing $k_i$, letting $j$ range over all scales $\le k_i$
and changing notations $(k_i,j)\to (j,h)$ yields the spring factors $M^{-\gamma_{I_{\ell^h}}(j-h)}, M^{-\gamma_{I'_{\ell^h}}(j-h)},M^{-\gamma_{I_{\Del^h,\tau}}(j-h)}$ of eq. (\ref{eq:Gh}). The
remaining factor $M^{-jD}=|\Del^j|$ may be rewritten as $M^{j\sum_{i\in I_q}\beta_i}$
which is distributed between the different fields, $\psi_i^{j_i}\to
\psi_i^{j_i,{\mathrm{rescaled}}}=M^{j\beta_i}\psi_i^{j_i}$ as in \S 6.1.2.  Recall from the end
of the preceding paragraph that each interval
of $\P$ comes with a factor $M^{-\eps}<1$.  Hence, letting $j_{max}$
be the maximal scale of a given  polymer $\P$
one has by  Lemma \ref{lem:multiscale-bound},
\BEA \sum_{\P\in {\cal P}_0^{j_{min}\to}(\vec{\Del}^{j_{min}\to})}
\left| {\mathrm{Eval}}_{\kappa'_0}(\P) \right|
 &\le&  \sum_{j_{max}=j_{min}}^{\rho}  \prod_{j=j_{min}}^{j_{max}}  \left(M^{-\eps}+\sum_{v\ge 1}
  I^{j,{\mathrm{rescaled}}}_{{\mathrm{Gaussian}}}(v;\vec{\Del}^{j_{min}\to})\right) \nonumber\\
&\le & \sum_{j_{max}=j_{min}}^{\rho} (M^{-\eps}+K'\lambda^{\kappa'_0})^{j_{max}-j_{min}+1} <\infty
\EEA
if $\lambda$ if chosen small enough so that in particular (with the constant $K$ of
eq. (\ref{eq:multiscale-Gaussian-bound})) $\sum_{v\ge 1} (K\lambda^{\kappa'_0})^v\le
K'\lambda^{\kappa'_0}<1-M^{-\eps}.$

\hfill \eop


\subsubsection{Combinatorial factors}

The last point -- for the Gaussian part of the final bounds  -- is to control the combinatorial
factors due to the horizontal and vertical cluster expansions. Let us show briefly how to do this.

As a general rule, differentiating a product of $n$ fields yields $n$ terms (by the Leibniz formula). This
implies supplementary combinatorial factors when estimating the polymer evaluation functions $F(\P)$, compared
to the estimates of Theorem \ref{th:6.1.3}. Consider the $O(n(\Del^j))$ derivations in a given
interval $\Del^j$ due to horizontal/vertical cluster expansion at scale $j$. A field produced by one
such derivation may be acted upon by another one in the same interval, yielding a local factorial
 $(O(n(\Del^j)))^{O(n(\Del^j))}$.
Otherwise, derivations in $\Del^j$ act on fields produced at an earlier stage in some interval
$\Del^k\subset\Del^j$ belonging to the same polymer $\P$ as $\Del^j$. Integrating over these, and borrowing
some small power of $\lambda$ from one of the differentiated fields,  yields
at most an averaged factor in $\Del^j$  of order
$K:=\frac{1}{|\Del^j|} \lambda^{\kappa} \sum_{k>j} \sum_{\Del^k\in\P\cap\D^k;\ \Del^k\subset\Del^j} \int_{\Del^k} dx $.
Once again, there are $O(n(\Del^j))$ such derivations. One may always gain a local factorial
$\frac{1}{n(\Del^j)!}$ to some arbitrary power (see proof of Lemma \ref{lem:single-scale-bound});
 using $\frac{K^{n(\Del^j)}}{n(\Del^j)!} \le e^K$, and
multiplying over all scales $j$ and all intervals $\Del^j\in\P\cap\D^j$, one obtains
\BEQ \exp \lambda^{\kappa} \sum_k \sum_{\Del^k\in\P\cap\D^k} \sum_{j<k} M^{-(k-j)},\EEQ
hence a factor of order $1+O(\lambda^{\kappa})$ per interval $\Del\in\P$, compensated by some factor $M^{-\eps}$ as
explained at the end of \S 6.1.2.


\subsection{Domination bounds}


Unlike Gaussian bounds, which are rather sophisticated, these are essentially based on the simple fact that
$|x|e^{-A|x|}=A^{-1} (A|x|)e^{-A|x|}\le KA^{-1}$ if $A>0$.

\begin{Lemma}[domination]  \label{lem:domination}

Let $\psi$ be a multiscale Gaussian field with scaling dimension $\beta$. Then
\BEA &&  |(T\psi)^{\to(k-1)}(\Del^k)|^n \exp -\lambda^{\kappa} M^{m\beta k}\ \cdot \
\frac{1}{|\Del^k|} \int_{\Del^k} \left( (T\psi)^{\to(k-1)}\right)^m(x)\ dx \nonumber\\
&& \qquad \qquad \qquad \qquad \qquad
\le K^n n^{n/m} \lambda^{-\kappa n/m} M^{-n\beta k}. \label{eq:dom} \EEA

\end{Lemma}

{\bf Proof.}

 Let $u:=(T\psi)^{\to(k-1)}(\Del^k)$ and $v:=\frac{1}{|\Del^k|} \int_{\Del^k} \left( (T\psi)^{\to(k-1)}\right)^m(x)\ dx$; by H\"older's inequality, $|u|\le v^{1/m}$,
so that
\BEA |u|^n e^{-\lambda^{\kappa} M^{m\beta k}v} &\le & (\frac{\lambda^{\kappa}}{n} M^{m\beta k})^{-n/m}
\ \cdot\ \left( (\frac{\lambda^{\kappa}}{n} M^{m\beta k}v)^{1/m}
 \exp -\frac{\lambda^k}{n} M^{m\beta k}v \right)^n \nonumber\\
&\le & K^n n^{n/m} \lambda^{-\kappa n/m} M^{-n\beta k}.\EEA

\hfill \eop

{\bf Example} ($(\phi,\partial\phi,\sigma)$-model). Lemma \ref{lem:domination} implies in particular the
  following four kinds of
low-momentum field domination:

\begin{itemize}

\item[(i)] Av$_{{\cal L}_4<\del{\cal L}_4}$ terms

Consider low-momentum fields $\sigma$ produced at a scale $k$ by letting some derivation  $\frac{\del}{\del \sigma}$
 or $\partial_t$
(due resp. to horizontal and vertical cluster expansions of scale $k$) act on ${\cal L}_4$.  When $\rho-k$ is large enough, they will be
dominated by the part of the counterterm $\del {\cal L}_4$ which is coupled to $b^{\rho-1}$. Eq. (\ref{eq:dom}) yields (using $1-t^2\ge (1-t)^2$ for $t\in[0,1]$)

\BEA &&  \frac{(1-t^k_{\Del})^n}{n!} \left| \lambda (T\sigma)^{\to(k-1)}(\Del^k) \right|^n e^{-\lambda^2
(1-(t^k_{\Del})^2)  M^{(1-4\alpha)k} \int_{\Del^k} |(T\sigma)^{\to(k-1)}
|^2(x) dx} \nonumber\\
&& \qquad \qquad \qquad  \le \left( KM^{2k\alpha} n^{-\half}\right)^n.\EEA
Note the factor $(1-t^k_{\Del})^n$ coming from the rest term in the Taylor-Lagrange expansion as in \S 3.3, with
$n=N_{ext,max}+O(n(\Del))$. This is precisely the reason why we chose to Taylor expand to order $N_{ext,max}+O(n(\Del))$ and
not simply to order $N_{ext,max}$. The other terms in the Taylor-Lagrange expansion have $t^k_{\Del}=0$.

Replacing in the exponential $\lambda^2 M^{(1-4\alpha)k}$ by the term $b^{\rho-1}\thickapprox \lambda^2
M^{(1-4\alpha)(\rho-1)}$, one gains a supplementary "petit facteur" $\left(M^{-\half(1-4\alpha)(\rho-1-k)}\right)^n$.

\item[(ii)] Av$_{{\cal L}_4<\del{\cal L}_{12}}$-terms

When $\rho-k$ is too small, the "petit facteur" is not small any more. One dominates by some fraction
(say, one tenth)  of the boundary term, $\frac{1}{10} \del{\cal L}_{12}$ instead. Again, eq. (\ref{eq:dom}) yields (using $1-(t_{\Del}^k)^6\ge (1-t_{\Del}^k)^6$)

\BEA &&    \frac{(1-t^k_{\Del})^n}{n!} \left|\lambda  (T\sigma)^{\to(k-1)}(\Del^k) \right|^n e^{-\frac{1}{10}\lambda^3
(1-(t^k_{\Del})^6)  M^{(1-12\alpha)k} \int_{\Del^k} |(T\sigma)^{\to(k-1)}
|^6(x) dx}  \nonumber\\
 &&\qquad \le  \frac{n^{n/6}}{n!} \left( K\lambda^{\half} M^{2\alpha k}\right)^n \le
\left(K'\lambda^{\half} M^{2\alpha k}\right)^n. \nonumber\\ \EEA

Replacing in the exponential $M^{(1-12\alpha)k}$ by the term $M^{(1-12\alpha)\rho}$ present in
$\del{\cal L}_{12}$ (note that $1-12\alpha<0$ !), one {\em loses} this time a supplementary large
factor $\left( M^{\frac{12\alpha-1}{6}(\rho-k)}\right)^n$. However, it is accompanied by a small
 factor $\lambda^{\half}$ per field.

Assume $\rho-k=0,1,\ldots,q$. We fix $q$ such that $\lambda^{\half} M^{q\frac{12\alpha-1}{6}}=
\lambda^{1/4}$, i.e. $q=\frac{3}{2(12\alpha-1)}\frac{\ln(1/\lambda)}{\ln M}.$ Thus the
 "petit facteur" in the preceding Av$_{{\cal L}_4<\del{\cal L}_4}$ case -- where $\rho-k>q$ -- is at most
$\lambda^{\kappa}$, where $\kappa=\frac{3(1-4\alpha)}{4(12\alpha-1)}$.

\item[(iii)] Av$_{\del{\cal L}_{12}<\del{\cal L}_{12}}$-terms

Consider low-momentum fields $\sigma$ produced from some vertex $\Del^{\rho'}$ by letting some derivation
 act on $\del{\cal L}_{12}$ this time. It produces  $i\le 5$ low-momentum $\sigma$-fields,
  accompanied by $\lambda^3$.
   Again, eq. (\ref{eq:dom}) yields, by dominating  by one tenth
  of the boundary term, $\frac{1}{10}\del{\cal L}_{12}$, and leaving out once and for all the $t$-coefficients,

\BEA && \frac{1}{n!} \left| \lambda^3
 ((T\sigma)^{\to(\rho'-1)}(\Del^{\rho'}))^i \right|^{n}
e^{-\frac{1}{10}\lambda^3
  M^{(1-12\alpha)\rho} \int_{\Del^{\rho'}} |(T\sigma)^{\to(\rho'-1)}
|^6(x) dx} \nonumber\\
 && \qquad \le
 \left(K'\lambda^{3/i-1/2} M^{2\alpha \rho'} M^{(\frac{12\alpha-1}{6})(\rho-\rho')}\right)^{in}. \nonumber\\ \EEA

  The corresponding vertex produces after
integration a factor $|\Del^{\rho'}|.M^{-(12\alpha-1)\rho}=M^{-12\alpha\rho}.M^{\rho-\rho'}$,
multiplied by $M^{2\alpha\rho'}=M^{2\alpha\rho}M^{-\frac{1}{6}(\rho-\rho')} M^{-\frac{12\alpha-1}{6}(\rho-\rho')}$ per high-momentum field, and by $M^{2\alpha \rho'}
M^{\frac{12\alpha-1}{6}(\rho-\rho')}\lambda^{3/i-1/2}= M^{2\alpha\rho}M^{-\frac{1}{6}(\rho-\rho')}\lambda^{3/i-1/2}$ per low-momentum
 field. Thus all together one has obtained a small factor $\le \lambda^{1/2} $  per vertex,
 to be shared between the six fields, and $M^{-\frac{12\alpha-1}{6}(\rho-\rho')}\le 1$ per high-momentum field.

\item[(iv)] Av$_{\del{\cal L}_4<\del{\cal L}_{4}}$

Consider low-momentum fields $\sigma$ produced in an interval $\Del^k$ of
scale $k$ by letting some derivation  act on $\del{\cal L}_4$ or on the counterterm
of scale $\rho$, $-\frac{b^{\rho}}{2}(t_x^{\rho})^2 ((T\sigma)^{\to\rho}(x))^2$.
It produces at most one low-momentum $\sigma$-field, accompanied by $b^{\rho'}=(b^{\rho'}/\lambda)
\cdot \lambda$
instead of $\lambda$.

Again, eq. (\ref{eq:dom}) yields, by dominating as in (i)
  by the part of the counterterm $\del {\cal L}_4$ which is coupled to $b^{\rho-1}$,
  a factor $O(M^{2k\alpha} \cdot M^{-\half(1-4\alpha)(\rho-k)})$ per field, alas multiplied by
  $b^{\rho'}/\lambda\thickapprox \lambda M^{\rho'(1-4\alpha)}$.  The rest of the argument goes as
  in subsection 6.1.4 -- a general argument called "aplatissement du fortement  connexe" in (colloquial)
  French. The factor $M^{-\half(1-4\alpha)(\rho-k)}$ may be simply bounded by $1$.

  Such terms may be produced  only in intervals $\Del^{\rho'}\subset\Del^k$ such that $\Del^{\rho'}\in\P$. Integrating
  over all such intervals -- and taking into account the $M^{-2k\alpha}$-scaling of the high-momentum
  field $\sigma^k$ left behind -- yields at most $\lambda K:=\lambda \sum_{\rho'\ge k}
  M^{-4\alpha(\rho'-k)}
  \#\{\P\cap\D^{\rho'}\}$ per low-momentum field produced, all together $(\lambda^{1/2})^n\cdot
  (\lambda^{\half}K)^n$. One may always gain a local factorial $\frac{1}{n!}$; using $\frac{K^n}{n!}
  \le e^K$ and multiplying over all scales $k$ and all intervals $\Del^k\in\P\cap\D^k$, one gets
  \BEQ \lambda^{v/2}\cdot \exp \lambda^{\half} \sum_{\rho'} \sum_{\Del^{\rho'}\in\P} \sum_{k\le\rho'}
  M^{-4\alpha(\rho'-k)}\EEQ
where $v$ is the  number of vertices where such low-momentum fields have been produced,
 hence a factor $\lambda^{\half}$ per vertex
and a factor of order $1+O(\lambda^{\half})$ per interval $\Del\in\P$.

\end{itemize}

\bigskip


\subsection{Final bounds}


The main Theorem is the following.

\begin{Theorem}  \label{th:6.2}

\begin{enumerate}
\item There exists a constant, positive two-by-two matrix  $K$ such that
\BEQ b^j=K\lambda^2 M^{j(1-4\alpha)}
(1+O(\lambda)).  \label{eq:bj-bound} \EEQ
\item There exists a constant $K'$ such that the Mayer bound of Proposition \ref{prop:Mayer-bound}
for the scale $j$ free energy,
\BEQ |f^{j\to\rho}(\lambda)|\le K'\lambda(1+O(\lambda)) \label{eq:Mayer-bound} \EEQ
 holds uniformly in $j$.
\end{enumerate}
\end{Theorem}

{\bf Proof.}

Let us first prove eq. (\ref{eq:bj-bound}) for $b^j$.

Consider a product (G-monomial)$\times$ (product of propagators) as in subsection 3.3, written generically as $GC$.
Multiply it by a product of averaged low-momentum fields of one of the four types
Av$_{{\cal L}_4<\del{\cal L}_4}$,  Av$_{{\cal L}_4<\del{\cal L}_{12}}$,
Av$_{\del{\cal L}_{12}<\del{\cal L}_{12}}$ or
Av$_{\del{\cal L}_4<\del{\cal L}_{4}}$, see subsection 6.2, generically written as $Av_{low}$.

\medskip

We make the following induction hypothesis:

\medskip

{\em Induction hypothesis. $\tilde{b}^k=K(1+O(\lambda))$ for all $k>j$ for some scale-independent constant $K$, where
$\tilde{b}^k:=\lambda^{-2} M^{-(1-4\alpha)k} b^k$ is the {\em rescaled} mass counterterm of scale $k$.}

\medskip

We shall soon see how to compute the constant $K$.
For the time being,
we must bound
\BEA && \int d\mu_{\vec{s}}(\phi) d\mu_{\vec{s}}(\sigma) \sum_{\P\in {\cal P}_{\sigma,\sigma}^{j\to}(\Del^j)}
\sum_{(G,C)}  Av_{low} GC e^{-\int_{|\P|} ({\cal L}_4+\del{\cal L}_4+\half\del{\cal L}_{12})(\ ;\ \vec{t})(x) dx} \ \cdot \nonumber\\
&& \qquad \cdot\ e^{-\int_{|\P|} \left[ \half \del{\cal L}_{12}^{\to\rho}(\ ;\vec{t})(x) -
\frac{b^{\rho}}{2} (t_x^{\rho})^2 ((T\sigma)^{\to\rho}(x))^2\right] dx}, \label{eq:dom-mass}\nonumber\\ \EEA
 where $|\P|$ is the support of the polymer $\P$ with two external $\sigma$-legs, in the notation of Proposition
\ref{prop:Mayer}.

\bigskip

{\em First step (domination of the mass counterterm).}

Note first that the term between square brackets in eq. (\ref{eq:dom-mass}) is {\em negative} when
$X:=\lambda M^{-2\rho\alpha} ||(T\sigma)^{\to\rho}||$ is small. Up to unessential coefficients,
it is equal to $M^{\rho}(\lambda^{-3}X^6-(t_x^{\rho})^2 X^2)$, which is minimal, of order $M^{\rho}\lambda^{9/2}$ --
 for $t_x^{\rho}\approx 1$, which is the worst case -- when $X$ is of order
$\lambda^{3/4}$. The factor $t_x^{\rho}$ in front of $(T\sigma)^{\to(\rho-1)}$ selects the intervals
in $\D^{\rho}$ belonging to the polymer. Hence the exponential in eq. (\ref{eq:dom-mass}) is
 bounded by $e^{K\int_{|\P\cap\D^{\rho}|} M^{\rho}\lambda^{9/2} dx}$, all together a factor
 of order $1+O(\lambda^{9/2})$ per interval $\Del\in\P\cap\D^{\rho}$ of the polymer with $t_{\Del}^{\rho}\not=0$.

\medskip

{\em Second step (domination of low-momentum fields).}  Split  $Av_{low} e^{-\int_{|\P|} ({\cal L}_4+\del{\cal L}_4+\half\del{\cal L}_{12})(\ ;\ \vec{t})(x) dx}$
from the expression in eq. (\ref{eq:dom-mass}).
As shown in subsection 6.2, these produce a small factor of order $\lambda^{\kappa}$, $\kappa>0$ per field. More precisely,
 any $\kappa<\inf(\frac{1}{4},\frac{3(1-4\alpha)}{4(12\alpha-1)})$ is suitable. There remains (by using the Cauchy-Schwarz inequality) to bound
$\sum_{\P\in {\cal P}_{\sigma,\sigma}^{j\to}(\Del^j)}
\sum_{(G,C)}   \left( \int d\mu_{\vec{s}}(\phi)d\mu_{\vec{s}}(\sigma) |GC|^2\right)^{\half}.$

\medskip

{\em Third step (computation of $b^j$).}

Let us now estimate $b^j$ by means of our induction hypothesis and of the Gaussian bounds of
 subsection 6.1. Consider for instance the diagonal term $b^j_{+,+}$. The terms with the fewest number of vertices are:

 -- the term with $0$ vertex obtained by applying twice $\frac{\partial}{\partial\sigma_+}$ to the counterterm
 $\del{\cal L}_4$, namely, $b^j_{+,+}$ (by definition);

\begin{figure}[h]
  \centering
   \includegraphics[scale=0.5]{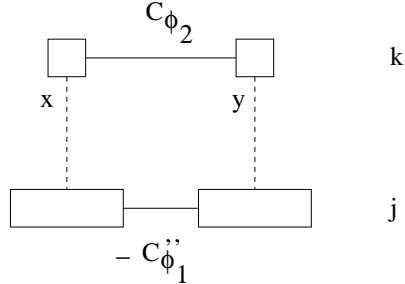}
   \caption{Main part of the mass counterterm of scale $j$.}
  \label{Fig-counterterm}
\end{figure}

 -- the polymers with two vertices, see Fig. \ref{Fig-counterterm}, which sum up to
 $$\lambda^2 \sum_{k\ge j}\int_V C_{\phi_2}^k(x,y) (-\partial^2)C_{\phi_1}^j(x,y)dy=:\lambda^2 M^{j(1-4\alpha)}
 K_{V},$$  where
 \BEA  && K_V:=M^{-j(1-4\alpha)}\int_V C_{\phi}^{j\to\rho}(x,y) (-\partial^2)C_{\phi}^j(x,y)dy
 \nonumber\\
 \qquad &&\to_{|V|\to\infty} K:=M^{-j(1-4\alpha)}\int |\xi|^{-4\alpha} \chi^j(\xi)\chi^{j\to}(\xi)d\xi=\int |\xi|^{-4\alpha}
 \chi^1(\xi)\chi^{1\to 2}(\xi) d\xi, \nonumber\\ \EEA
  a scale-independent quantity. As in section 1, the non-diagonal counterterm $b^j_{+,-}$ or $b^j_{-,+}$ may be computed in the same way, yielding in the end a scale-independent positive two-by-two matrix.

 More complicated polymers are of order $M^{j(1-4\alpha)} \lambda^2 g(\lambda^2 \tilde{b}^j,\lambda^2 (\tilde{b}^k)_{k>j},
\lambda)$,
 where $\tilde{b}^j:=\lambda^{-2}M^{-j(1-4\alpha)}b^j$ is the rescaled mass counterterm as in
 the induction hypothesis. The
 function $g$ is $C^{\infty}$ in a neighbourhood of $0$ and vanishes at $0$. Hence, by the implicit
 function theorem, $b^j=KM^{j(1-4\alpha)}(1+O(\lambda))$ as in Proposition \ref{prop:Mayer}.

\bigskip

The bound for the scale $j$ free energy $f^{j\to\rho}$ is now straightforward.
\hfill \eop

\bigskip

Bounds for $n$-point functions are easy generalizations of the preceding Theorem.
Consider for instance the $2$-point function $\langle |{\cal F}\phi_1(\xi)|^2
\rangle_{\lambda}$, with $M^j\le|\xi|\le M^{j+1}$. By momentum conservation,
and by definition of the Fourier partition of unity, see subsection 2.1, this is equal to the sum over $j_1,j_2=j,j\pm 1$ of  $\langle {\cal F}\phi^{j_1}_1(\xi) {\cal F}\phi^{j_2}_1(-\xi)
\rangle_{\lambda}$. The term of order $0$ in $\lambda$ is given by the Gaussian evaluation
$\esper \left[ {\cal F}\phi_1^{j_1}(\xi) {\cal F}\phi_1^{j_2}(-\xi)\right]$. Further terms
involve at least one $\sigma$-propagator with a small factor (see Lemma \ref{lem:2.13})
 of order $\inf(1,\frac{M^{j(1-4\alpha)}}{b^{\rho}})\le K \inf(1,\lambda^{-2} M^{-(\rho-j)(1-4\alpha)})$
 which goes to 0 when $\rho\to\infty$.



\end{document}